%% file: singlestable_24_arxiv.tex
\UseRawInputEncoding 

\documentclass{article}

\usepackage{latexsym}
\usepackage{amsfonts}
\usepackage{amssymb}
\usepackage{cite}

\usepackage{proof}

\input{header}

\title{
An ordinal analysis of a single stable ordinal
}
\author{Toshiyasu Arai
\\
Graduate School of Mathematical Sciences
\\
University of Tokyo
\\
3-8-1 Komaba, Meguro-ku,
Tokyo 153-8914, JAPAN
\\
tosarai@ms.u-tokyo.ac.jp
}

\date{}

\begin{document}

\maketitle

\begin{abstract}
In this paper we give an ordinal analysis of a set theory extending  ${\sf KP}\ell^{r}$ 
with an axiom stating that
`there exists a transitive set $M$ such that $M\prec_{\Sigma_{1}}V$'.
\end{abstract}

\section{Introduction}\label{sect:introduction}

In this paper we give an ordinal analysis of a set theory 
${\sf KP}\ell^{r}+(M\prec_{\Sigma_{1}}V)$ extending  ${\sf KP}\ell^{r}$ 
with an axiom stating that
`there exists a transitive set $M$ such that $M\prec_{\Sigma_{1}}V$'.
An ordinal analysis of an extension
${\sf KP}i+(M\prec_{\Sigma_{1}}V)$
 is given in M. Rathjen\cite{RathjenAFML2}.

$\Sigma^{1-}_{2}\mbox{-CA}+\Pi^{1}_{1}\mbox{-CA}_{0}$
 is a second order arithmetic obtained from 
$\Pi^{1}_{1}\mbox{-CA}_{0}$
 by adding the Comprehension Axiom for
 parameter free $\Sigma^{1}_{2}$-formulas.
 It is easy to see that $\Sigma^{1-}_{2}\mbox{-CA}+\Pi^{1}_{1}\mbox{-CA}_{0}$ is interpreted canonically to the set theory
${\sf KP}\ell^{r}+(M\prec_{\Sigma_{1}}V)$.

To obtain an upper bound of the proof-theoretic ordinal of 
${\sf KP}\ell^{r}+(M\prec_{\Sigma_{1}}V)$, 
we employ operator controlled derivations introduced
by W. Buchholz\cite{Buchholz}, 
in which a set theory ${\sf KP}i$ for recursively inaccessible universes is analyzed
proof-theoretically.
Our proof is an extension of \cite{KPPiN} 
in which a set theory ${\sf KP}\Pi_{N}$ of $\Pi_{N}$-reflection is analyzed, while
\cite{KPPiN} 
is an extension of M. Rathjen's analysis in \cite{Rathjen94} 
for 
$\Pi_{3}$-reflection.
A new ingredient is a use of an explicit Mostowski collapsing as in \cite{sneak}.

The set theory ${\sf KP}\ell^{r}$ in J\"ager's monograph\cite{J3} 
is obtained from the Kripke-Platek set theory ${\sf KP}\omega$ with
the axiom of Infinity, cf.\,\cite{Ba, J3}, 
by deleting $\Delta_{0}$-Collection on the universe,
restricting Foundation schema to $\Delta_{0}$-formulas, and
adding an axiom $(Lim)$, $\forall x\exists y(x\in y\land Ad(y))$, stating that
the universe is a limit of admissible sets, where $Ad$ is a unary predicate such that $Ad(y)$ is intended to designate that
`$y$ is a (transitive and) admissible set'.
Then a set theory ${\sf KP}\ell^{r}+(M\prec_{\Sigma_{1}}V)$ extends ${\sf KP}\ell^{r}$
 by adding an individual constant $M$ and
the axioms for the constant $M$:
$M$ is non-empty, transitive and stable, $M\prec_{\Sigma_{1}}V$ for the universe $V$:
$M\neq\emptyset$, $\forall x\in M\forall y\in x(y\in M)$ and
\begin{equation}\label{eq:stable}
 \varphi(u_{1},\ldots,u_{n}) \land \{u_{1},\ldots,u_{n}\}\subset M \to \varphi^{M}(u_{1},\ldots,u_{n})
\end{equation}
for each $\Sigma_{1}$-formula $\varphi$ in the set-theoretic language.

Note that $M$ is a model of ${\sf KP}\omega$ and the axiom $(Lim)$, i.e., a model of ${\sf KP}i$.

For positive integers $N$, let us define a subtheory $T_{N}$ of ${\sf KP}\ell^{r}+(M\prec_{\Sigma_{1}}V)$.
The intended model of $T_{N}$ is an admissible set $M_{N}$ in which there is a set $M$ with
$M\prec_{\Sigma_{1}}M_{N}$, and there are $(N-1)$ admissible sets $M_{n}\,(0<n<N)$
above $M$ such that $M\in M_{1}\in\cdots\in M_{N-1}\in M_{N}$.

\bdf\label{df:TN}
{\rm
Let $N>0$ be a positive integer.
$T_{N}$ denotes a set theory defined as follows.
The language of $T_{N}$ is $\{\in\}\cup\{M_{i}\}_{i<N}$ with individual constants $M_{i}$.
$T_{N}$ is obtained from the set theory ${\sf KP}\omega+(M\prec_{\Sigma_{1}}V)$
with $M:=M_{0}$ by 
adding axioms
$M_{n}\in M_{n+1}$ for $n+1<N$ and axioms stating that
each $M_{i}$ is a transitive admissible set for $n<N$.
}
\edf

\bprp\label{prp:reduction}
For each $\Sigma_{1}$-formula $\theta$,
if ${\sf KP}\ell^{r}+(M\prec_{\Sigma_{1}}V)\vdash\theta$, then there exists an $N$ such that
$T_{N}\vdash\theta$.
\eprp
\bprf
This is seen through a partial cut elimination and an asymmetric interpretation.
Note that each axiom (\ref{eq:stable}) is a $\Pi_{1}$-formula.
\eprf
\\

In the following theorems,
$\Omega=\omega_{1}^{CK}$ and $\psi_{\Omega}$ denotes a collapsing function such that $\psi_{\Omega}(\alpha)<\Omega$.
$\mathbb{S}$ is an ordinal term denoting a stable ordinal, and 
$\Omega_{\mathbb{S}+N}$ the $N$-th admissible ordinal above $\mathbb{S}$.
Our aim here is to show the following Theorem \ref{thm:2}, where
$\omega_{0}(\alpha)=\alpha$ and $\omega_{n+1}(\alpha)=\omega^{\omega_{n}(\alpha)}$.

\begin{theorem}\label{thm:2}
Suppose $T_{N}\vdash\theta^{L_{\Omega}}$
for a $\Sigma_{1}$-sentence $\theta$.
Then we can find an $n<\omega$ such that for $\alpha=\psi_{\Omega}(\omega_{n}(\Omega_{\mathbb{S}+N}+1))$,
$L_{\alpha}\models\theta$ holds.
\end{theorem}

Actually the bound is seen to be tight as the following Theorem \ref{th:wf} in \cite{singledistwfprf}
shows.
$OT$ denotes a computable notation system of ordinals,
and $OT_{N}$ a restriction of $OT$ such that $OT=\bigcup_{0<N<\omega}OT_{N}$ and
$\psi_{\Omega}(\varepsilon_{\Omega_{\mathbb{S}+N}+1})$ denotes
the order type of $OT_{N}\cap\Omega$.

\begin{theorem}\label{th:wf}
$\Sigma^{1-}_{2}\mbox{{\rm -CA}}+\Pi^{1}_{1}\mbox{{\rm -CA}}_{0}$
 proves that $(OT_{N},<)$ is a well ordering for {\rm each} $N$.
\end{theorem}

Thus the ordinal $\psi_{\Omega}(\Omega_{\mathbb{S}+\omega}):=\sup\{\psi_{\Omega}(\varepsilon_{\Omega_{\mathbb{S}+N}+1}): 0<N<\omega\}$ 
is the proof-theoretic ordinal of $\Sigma^{1-}_{2}\mbox{-CA}+\Pi^{1}_{1}\mbox{-CA}_{0}$ and of
${\sf KP}\ell^{r}+(M\prec_{\Sigma_{1}}V)$, where
$|{\sf KP}\ell^{r}+(M\prec_{\Sigma_{1}}V)|_{\Sigma_{1}^{\Omega}}$ denotes the 
$\Sigma_{1}^{\Omega}$-ordinal of ${\sf KP}\ell^{r}+(M\prec_{\Sigma_{1}}V)$, i.e., the ordinal
$\min\{\alpha\leq\omega_{1}^{CK}: \forall \theta\in\Sigma_{1}
\left(
{\sf KP}\ell^{r}+(M\prec_{\Sigma_{1}}V)\vdash\theta^{L_{\Omega}} \Rightarrow L_{\alpha}\models\theta
\right)\}$.

\begin{theorem}\label{th:main}
$\psi_{\Omega}(\Omega_{\mathbb{S}+\omega}) 
=|\Sigma^{1-}_{2}\mbox{{\rm -CA}}+\Pi^{1}_{1}\mbox{{\rm -CA}}_{0}|
= |{\sf KP}\ell^{r}+(M\prec_{\Sigma_{1}}V)|_{\Sigma_{1}^{\Omega}}.
$
\end{theorem}
Moreover
${\sf KP}\ell^{r}+(M\prec_{\Sigma_{1}}V)$ is seen to be
conservative over 
${\rm I}\Sigma_{1}+\{TI(\alpha,\Sigma^{0}_{1}(\omega)): \alpha<\psi_{\Omega}(\Omega_{\mathbb{S}+\omega})\}$
with respect to $\Pi^{0}_{2}(\omega)$-arithmetic sentences, where
$TI(\alpha,\Sigma^{0}_{1}(\omega))$ denotes the schema of transfinite induction up to $\alpha$,
and applied to $\Sigma^{0}_{1}$-arithmetic formulas in a language of the first-order arithmetic.
In particular each provably computable function in ${\sf KP}\ell^{r}+(M\prec_{\Sigma_{1}}V)$ is defined
by $\alpha$-recursion for an $\alpha<\psi_{\Omega}(\Omega_{\mathbb{S}+\omega})$, cf.\,Corollary \ref{cor:Pi2conservative}.

Let us mention the contents of this paper, and give a brief sketch of proofs.
In the next section \ref{sect:ordinalnotation} we define simultaneously
iterated Skolem hulls 
$\mathcal{H}_{\alpha}(X)$ of sets $X$ of ordinals, ordinals 
$\psi^{f}_{\kappa}(\alpha)$ for regular cardinals 
$\kappa$, and finite functions $f$,
and Mahlo classes $Mh^{\alpha}_{k}(\xi)$.
We invoke an existence of a shrewd cardinal introduced by M. Rathjen\cite{RathjenAFML2},
to justify the definitions.

Following W. Buchholz\cite{Buchholz},
operator controlled derivations are introduced,
 and inference rules for stability and reflections are eliminated from
derivations in section \ref{sect:controlledOme}.
Roughly the elimination procedure goes as follows.
First the axiom (\ref{eq:stable}) for stability is replaced by axioms of reflections
through a Mostwoski collapsing $\pi$.
Let $\exists x \varphi(x, u)$ be a $\Sigma_{1}$-formula with a parameter $u\in L_{\mathbb{S}}$ for a stable ordinal $\mathbb{S}$.
Let $u\in L_{\alpha}$ with an $\alpha<\mathbb{S}$.
Suppose $\varphi(v, u)$ is true for a $v$.
Since we are considering infinitary images of finite derivations, we can assume that
the set $v$ is constructed from some finite parameters $\{\beta_{i}\}_{i}$ of ordinals.
A Skolem hull of $L_{\alpha}\cup\{\beta_{i}\}_{i}\cup\{\Omega_{\mathbb{S}+n}: n\leq N\}$ under some functions such as
$\beta\mapsto \psi_{\Omega_{\mathbb{S}+n+1}}(\beta)$, where 
$\Omega_{\mathbb{S}+n}<\psi_{\Omega_{\mathbb{S}+n+1}}(\beta)<\Omega_{\mathbb{S}+n+1}$,
is collapsed down to an $L_{\rho}$ with $\alpha<\rho<\mathbb{S}$ through a Mostowski collapsing $\pi$.
Suppose that $\varphi(v,u)$ implies $\varphi(\pi(v),u)$. Then the axiom (\ref{eq:stable}) follows.
In the implication a reflection in a transfinite level involves
since $v\in L_{\beta}$ possibly with $\beta>\mathbb{S}$.
In other words, $\exists \beta<\Omega_{\mathbb{S}+N}\exists x\in L_{\beta}\varphi(x,u)$ should yield $\exists \beta_{0}<\mathbb{S}\exists x\in L_{\beta_{0}}\varphi(x,u)$,
where $\exists \beta<\Omega_{\mathbb{S}+N}\exists x\in L_{\beta}\varphi(x,u)$ is a $\Sigma_{\Omega_{\mathbb{S}+N}}$-sentence
over $L_{\mathbb{S}}$, so to speak.
To resolve such a transfinite reflection, we need collapsing functions $\psi_{\kappa}^{f}(\alpha)$
with finite functions $f$ indicating Mahlo degrees.
In \cite{KPPiN} 
we introduced
collapsing functions $\psi_{\kappa}^{\vec{\xi}}(\alpha)$ with finite sequences $\vec{\xi}$
of ordinals in length $N-2$ to analyze $\Pi_{N}$-reflection.
Here for transfinite reflections we need functions $f$ of finite supports.

IH denotes the Induction Hypothesis, MIH the Main IH and SIH the Subsidiary IH.
Throughout of this paper $N$ denotes a fixed positive integer.

\section{Ordinals for one stable ordinal}\label{sect:ordinalnotation}
In this section up to Lemma \ref{lem:welldefinedness.1},
we work in the set theory obtained from ${\sf ZFC}$
by adding the axiom stating that there exists a weakly inaccessible cardinal $\mathbb{S}$.
For ordinals $\alpha\geq\beta$,
$\alpha-\beta$ denotes the ordinal $\gamma$ such that $\alpha=\beta+\gamma$.
Let $\alpha$ and $\beta$ be ordinals.
$\alpha\dot{+}\beta$ denotes the sum $\alpha+\beta$
when $\alpha+\beta$ equals to the commutative (natural) sum $\alpha\#\beta$, i.e., when
either $\alpha=0$ or $\alpha=\alpha_{0}+\omega^{\alpha_{1}}$ with
$\omega^{\alpha_{1}+1}>\beta$.

$\mathbb{S}$ denotes a weakly inaccessible cardinal with $\omega_{\mathbb{S}}=\mathbb{S}$, and
$\omega_{\mathbb{S}+n}$ the $n$-th regular cardinal above $\mathbb{S}$ for $0<n<\omega$.
Let $\mathbb{K}=\omega_{\mathbb{S}+N}$ for a fixed positive integer $N$.

\bdf\label{df:Lam}
{\rm
For a positive integer $N$,
$\varphi_{b}(\xi)$ denotes the binary Veblen function on 
$\mathbb{K}=\omega_{\mathbb{S}+N}$ with $\varphi_{0}(\xi)=\omega^{\xi}$,
and
$\tilde{\varphi}_{b}(\xi):=\varphi_{b}(\Lambda\cdot \xi)$ for an epsilon number  $\Lambda<\mathbb{K}$.
Each ordinal $\varphi_{b}(\xi)$ is a fixed point of the function $\varphi_{c}$ for $c<b$
in the sense that $\varphi_{c}(\varphi_{b}(\xi))=\varphi_{b}(\xi)$.
The same holds for $\tilde{\varphi}_{b}$ and $\tilde{\varphi}_{c}$ with $c<b$.

Let $b,\xi<\Lambda^{+}$.
$\theta_{b}(\xi)$ 
[$\tilde{\theta}_{b}(\xi)$] 
denotes
a $b$-th iterate of $\varphi_{0}(\xi)=\omega^{\xi}$.
[of $\tilde{\varphi}_{0}(\xi)=\Lambda^{\xi}$], resp.
Specifically ordinals
$ \tilde{\theta}_{b}(\xi)$ 
are defined by recursion on $b$ as follows.
$\theta_{0}(\xi)
=\tilde{\theta}_{0}(\xi)
=\xi$, $\theta_{\omega^{b}}(\xi)=\varphi_{b}(\xi)$, and
$\tilde{\theta}_{\omega^{b}}(\xi)=\tilde{\varphi}_{b}(\xi)$, and
$\theta_{c\dot{+}\omega^{b}}(\xi)=\theta_{c}(\theta_{\omega^{b}}(\xi))$.
$\tilde{\theta}_{c\dot{+}\omega^{b}}(\xi)=\tilde{\theta}_{c}(\tilde{\theta}_{\omega^{b}}(\xi))$.

$\alpha>0$ is a strongly critical number if $\forall b,\xi<\alpha(\varphi_{b}(\xi)<\alpha)$.
$\Gamma_{a}$ denotes the $a$-th strongly critical number, and
$\Gamma(a)$ the next strongly critical number above $a$, while
$\varepsilon_{a}$ denotes the $a$-th epsilon number, and
$\varepsilon(a)$ the next epsilon number above $a$.
}
\edf

Let us define a normal form of non-zero ordinals $\xi<\Gamma(\Lambda)$.
Let $\xi=\Lambda^{\zeta}$.
If $\zeta<\Lambda^{\zeta}$, then $\tilde{\theta}_{1}(\zeta)$ is the normal form of $\xi$,
denoted by
$\xi=_{NF}\tilde{\theta}_{1}(\zeta)$.
Assume $\zeta=\Lambda^{\zeta}$, and let $b>0$ be the maximal ordinal such that
there exists an ordinal $\eta$ with $\xi=\tilde{\varphi}_{b}(\eta)$.
Then $\xi=\tilde{\varphi}_{b}(\eta)=_{NF}\tilde{\theta}_{\omega^{b}}(\eta)$.

Let $\xi=\tilde{\theta}_{b_{m}}(\xi_{m})\cdot a_{m}+\cdots+\tilde{\theta}_{b_{0}}(\xi_{0})\cdot a_{0}$, where
$\tilde{\theta}_{b_{i}}(\xi_{i})>\xi_{i}$,
$\tilde{\theta}_{b_{m}}(\xi_{m})>\cdots>\tilde{\theta}_{b_{0}}(\xi_{0})$, 
$b_{i}=\omega^{c_{i}}<\Lambda$, and
$0<a_{0},\ldots,a_{m}<\Lambda$.
Then
$\xi=_{NF}\tilde{\theta}_{b_{m}}(\xi_{m})\cdot a_{m}+\cdots+\tilde{\theta}_{b_{0}}(\xi_{0})\cdot a_{0}$.

\bdf\label{df:Lam2}
{\rm
Let $\xi<\Gamma(\Lambda)$ be a non-zero ordinal with its normal form\footnote{The normal form in (\ref{eq:CantornfLam}) is slightly extended from \cite{singledistwfprf}.}:
\begin{equation}\label{eq:CantornfLam}
\xi=\sum_{i\leq m}\tilde{\theta}_{b_{i}}(\xi_{i})\cdot a_{i}=_{NF}
\tilde{\theta}_{b_{m}}(\xi_{m})\cdot a_{m}+\cdots+\tilde{\theta}_{b_{0}}(\xi_{0})\cdot a_{0}
\end{equation}
$SC_{\Lambda}(\xi)=\bigcup_{i\leq m}(\{b_{i},a_{i}\}\cup SC_{\Lambda}(\xi_{i}))$.
$\tilde{\theta}_{b_{0}}(\xi_{0})$ is said to be the \textit{tail} of $\xi$, denoted 
$\tilde{\theta}_{b_{0}}(\xi_{0})=tl(\xi)$, and
$\tilde{\theta}_{b_{m}}(\xi_{m})$ the \textit{head} of $\xi$, denoted 
$\tilde{\theta}_{b_{m}}(\xi_{m})=hd(\xi)$.

\begin{enumerate}
\item\label{df:Exp2.3}
 $\zeta$ is a \textit{part} of $\xi$
 iff there exists an $n\, (0\leq n\leq m+1)$
 such that
 $\zeta=_{NF}\sum_{i\geq n}\tilde{\theta}_{b_{i}}(\xi_{i})\cdot a_{i}=
 \tilde{\theta}_{b_{m}}(\xi_{m})\cdot a_{m}+\cdots+\tilde{\theta}_{b_{n}}(\xi_{n})\cdot a_{n}$
 for $\xi$ in (\ref{eq:CantornfLam}).

\item\label{df:thtminus}
Let $\zeta=_{NF}\tilde{\theta}_{b}(\xi)$ with $\tilde{\theta}_{b}(\xi)>\xi$ and $b=\omega^{b_{0}}$,
and $c$ be ordinals.
An ordinal $\tilde{\theta}_{-c}(\zeta)$ is defined recursively as follows.
If $b\geq c$, then $\tilde{\theta}_{-c}(\zeta)=\tilde{\theta}_{b-c}(\xi)$.
Let $c>b$.
If $\xi>0$, then
$\tilde{\theta}_{-c}(\zeta)=\tilde{\theta}_{-(c-b)}(\tilde{\theta}_{b_{m}}(\xi_{m}))$ for the head term 
$hd(\xi)=\tilde{\theta}_{b_{m}}(\xi_{m})$ of 
$\xi$ in (\ref{eq:CantornfLam}).
If $\xi=0$, then let $\tilde{\theta}_{-c}(\zeta)=0$.

\end{enumerate}
}
\edf

A `Mahlo degree' $m(\pi)$ of ordinals $\pi$
with higher reflections is defined to be a finite function 
$f:\Lambda\to\Gamma(\Lambda)$.

\bdf\label{df:finitefunction}
{\rm
  \begin{enumerate}
  \item
 A function $f:\Lambda\to\Gamma(\Lambda)$ with a \textit{finite} support
${\rm supp}(f)=\{c<\Lambda: f(c)\neq 0\}\subset \varphi_{\Lambda}(0)$ is said to be a \textit{finite function}
if
$\forall i>0(a_{i}=1)$ and $a_{0}=1$ when $b_{0}>1$
$f(c)=_{NF}\tilde{\theta}_{b_{m}}(\xi_{m})\cdot a_{m}+\cdots+\tilde{\theta}_{b_{0}}(\xi_{0})\cdot a_{0}$
for any $c\in{\rm supp}(f)$.

It is identified with the finite function $f\!\upharpoonright\! {\rm supp}(f)$.
When $c\not\in {\rm supp}(f)$, let $f(c):=0$.
$SC_{\Lambda}(f):=\bigcup\{\{c\}\cup SC_{\Lambda}(f(c)): c\in {\rm supp}(f)\}$.
$f,g,h,\ldots$ range over finite functions.

For an ordinal $c$, $f_{c}$ and $f^{c}$ are restrictions of $f$ to the domains
${\rm supp}(f_{c})=\{d\in{\rm supp}(f): d< c\}$ and ${\rm supp}(f^{c})=\{d\in{\rm supp}(f): d\geq c\}$.
$g_{c}*f^{c}$ denotes the concatenated function such that
${\rm supp}(g_{c}*f^{c})={\rm supp}(g_{c})\cup {\rm supp}(f^{c})$, 
$(g_{c}*f^{c})(a)=g(a)$ for $a<c$, and
$(g_{c}*f^{c})(a)=f(a)$ for $a\geq c$.

\item\label{df:Exp2.5}
Let $f$ be a finite function and $c,\xi$ ordinals.
A relation $f<^{c}\xi$ is defined by induction on the
cardinality of the finite set $\{d\in {\rm supp}(f): d>c\}$ as follows.
If $f^{c}=\emptyset$, then $f<^{c}\xi$ holds.
For $f^{c}\neq\emptyset$,
$f<^{c}\xi$ iff
there exists a part $\mu$ of $\xi$ such that
$f(c)< \mu$
and 
$f<^{c+d} \tilde{\theta}_{-d}(tl(\mu))$ 
for $d=\min\{c+d\in {\rm supp}(f): d>0\}$.

\end{enumerate}

}
\edf

The following Proposition \ref{prp:tht4} is shown in \cite{singledistwfprf}.

\bprp\label{prp:tht4}
\begin{enumerate}
\item\label{prp:tht4.1}
$\zeta\leq\xi \Rightarrow \tilde{\theta}_{-c}(\zeta)\leq\tilde{\theta}_{-c}(\xi)$.

\item\label{prp:tht4.2}
$\tilde{\theta}_{c}(\tilde{\theta}_{-c}(\zeta)) \leq\zeta$.
\end{enumerate}
\eprp

\bprp\label{prp:idless}
$f<^{c}\xi\leq\zeta \Rightarrow f<^{c}\zeta$.
\eprp
\bprf
By induction on the cardinality $n$
of the finite set
$\{d\in {\rm supp}(f): d> c\}=\{c<c+d_{1}<\cdots<c+d_{n}\}$.
If $n=0$, then $f(c)<\xi\leq\zeta$ yields $f<^{c}\zeta$.
Let $n>0$.
We have $f(c)<\mu$, 
and $f<^{c+d_{1}}\tilde{\theta}_{-d_{1}}(tl(\mu))$ 
for a part $\mu$ of $\xi$.
We show the existence of a part $\lambda$ of $\zeta$ such that 
$\mu\leq\lambda$, 
and $\tilde{\theta}_{-d_{1}}(tl(\mu)) \leq \tilde{\theta}_{-d_{1}}(tl(\lambda))$.
Then IH yields $f<^{c+d_{1}}\tilde{\theta}_{-d_{1}}(tl(\lambda))$, and $f<^{c}\zeta$ follows.

If $\mu$ is a part of $\zeta$, then $\lambda=\mu$ works.
Otherwise $\xi<\zeta$ and there exists a part $\lambda$ of $\zeta$ such that $\mu<\lambda$,
and
$tl(\mu)<tl(\lambda)$.
We obtain $\tilde{\theta}_{-d_{1}}(tl(\mu))\leq\tilde{\theta}_{-d_{1}}(tl(\lambda))$.
\eprf

$u,v,w,x,y,z,\ldots$ range over sets in the universe,
$a,b,c,\alpha,\beta,\gamma,\delta,\ldots$ range over ordinals$<\varepsilon(\Lambda)$,
$\xi,\zeta,\eta,\ldots$ range over ordinals$<\Gamma(\Lambda)$,
and $\pi,\kappa,\rho,\sigma,\tau,\lambda,\ldots$ range over regular ordinals.

\subsection{Skolem hulls and collapsing functions}\label{subsec:Skolemh}
In this subsection Skolem hulls $\mathcal{H}_{a}(\alpha)$, collapsing functions $\psi_{\pi}^{f}(\alpha)$ and
Mahlo classes $Mh^{a}_{c}(\xi)$ are defined.

\bdf
{\rm
\begin{enumerate}

\item\label{df:Lam3}
Let $A\subset\mathbb{S}$ be a set, and $\alpha\leq\mathbb{S}$ a limit ordinal.
\[
\alpha\in M(A) :\Leftrightarrow A\cap\alpha \mbox{ is stationary in }
\alpha
\Leftrightarrow \mbox{ every club subset of } \alpha \mbox{ meets } A.
\]

 \item\label{df:Lam4}
 $\kappa^{+}$ denotes the next regular ordinal above $\kappa$.
 For $n<\omega$, $\kappa^{+n}$ is defined recursively by
 $\kappa^{+0}=\kappa$ and $\kappa^{+(n+1)}=\left(\kappa^{+n}\right)^{+}$.

\item
$\Omega_{\alpha}:=\omega_{\alpha}$ for $\alpha>0$, $\Omega_{0}:=0$, and
$\Omega=\Omega_{1}$.
$\mathbb{S}$ is a weakly inaccessible cardinal, and
$\Omega_{\mathbb{S}}=\mathbb{S}$.
$\Omega_{\mathbb{S}+n}=\mathbb{S}^{+n}$ is the $n$-th cardinal above $\mathbb{S}$.
\end{enumerate}
}
\edf

In the following Definition \ref{df:Cpsiregularsm}, 
$\varphi\alpha\beta=\varphi_{\alpha}(\beta)$ denotes the binary Veblen function on $\mathbb{K}^{+}=\omega_{\mathbb{S}+N+1}$ with $\varphi_{0}(\beta)=\omega^{\beta}$,
$\tilde{\theta}_{b}(\xi)$ the function defined in Definition \ref{df:Lam}
for $\Lambda<\mathbb{K}=\omega_{\mathbb{S}+N}$ with the positive integer $N$.

For 
$a<\varepsilon(\mathbb{K})$,
$c<\Lambda$, and
$\xi<\Gamma(\Lambda)$, 
define simultaneously 
classes $\mathcal{H}_{a}(X)\subset\Gamma(\mathbb{K})$,
$Mh^{a}_{c}(\xi)\subset(\mathbb{S}+1)$, and 
ordinals $\psi_{\kappa}^{f}(a)\leq\kappa$ by recursion on ordinals $a$ as follows.
We see that these are 
$\Delta_{1}$-definable in {\sf ZFC},
 cf.\,Proposition \ref{prp:definability}.

\bdf\label{df:Cpsiregularsm}
{\rm
Let
$\mathbb{K}=\Omega_{\mathbb{S}+N}$, 
$\mathcal{H}_{a}[Y](X):=\mathcal{H}_{a}(Y\cup X)
$ for sets $Y\subset \Gamma(\mathbb{K})$.
Let $a<\varepsilon(\mathbb{K})$ and $X\subset\Gamma(\mathbb{K})$.

\begin{enumerate}
\item\label{df:Cpsiregularsm.1}
(Inductive definition of $\mathcal{H}_{a}(X)$.)

\begin{enumerate}
\item\label{df:Cpsiregularsm.10}
$\{0,\Omega_{1},\mathbb{S}\}\cup\{\Omega_{\mathbb{S}+n}: 0<n\leq N\}\cup X\subset\mathcal{H}_{a}(X)$.

\item\label{df:Cpsiregularsm.11}
If $x, y \in \mathcal{H}_{a}(X)$,
then $x+y\in \mathcal{H}_{a}(X)$,
and 
$\varphi xy\in \mathcal{H}_{a}(X)$.

\item\label{df:Cpsiregularsm.12}
Let $\alpha\in\mathcal{H}_{a}(X)\cap\mathbb{S}$. Then for each $0<k\leq N$, 
$\Omega_{\alpha+k}\in\mathcal{H}_{a}(X)$.

\item\label{df:Cpsiregularsm.1345}
Let $\alpha=\psi_{\pi}^{f}(b)$ with $\{\pi,b\}\subset\mathcal{H}_{a}(X)$, 
$\pi\in\{\mathbb{S}\}\cup\Psi_{N}$,
$b<a$, and a finite function $f$ such that
$SC(f)\subset\mathcal{H}_{a}(X)\cap\mathcal{H}_{b}(\alpha)$.

Then $\alpha\in\mathcal{H}_{a}(X)\cap\Psi_{N}$.

\end{enumerate}

\item\label{df:Cpsiregularsm.2}
 (Definitions of $Mh^{a}_{c}(\xi)$ and $Mh^{a}_{c}(f)$)
\\
The classes $Mh^{a}_{c}(\xi)$ are defined for $c<\Lambda$, $\xi< \Gamma(\Lambda)$,
and ordinals $a<\varepsilon(\mathbb{K})$.
Let $\pi$ be a regular ordinal$\leq \mathbb{S}$. Then 
by main induction on ordinals $\pi\leq\mathbb{S}$
with subsidiary induction on $c<\Lambda$ 
we define $\pi\in Mh^{a}_{c}(\xi)$ iff 
$\{a,c,\xi\}\subset\mathcal{H}_{a}(\pi)$ and
{\small
\begin{equation}\label{eq:dfMhkh}
 \forall f<^{c}\xi 
 \forall g \left(
 SC_{\Lambda}(f)\cup SC_{\Lambda}(g) \subset\mathcal{H}_{a}(\pi) 
 \,\&\, 
\pi\in Mh^{a}_{0}(g_{c})
 \Rightarrow \pi\in M(Mh^{a}_{0}(g_{c}*f^{c}))
 \right)
\end{equation}
}
where $f, g$ vary through 
finite
 functions from $\Lambda$ to $\varphi_{\Lambda}(0)$,
and 
\[
Mh^{a}_{c}(f)  :=  \bigcap\{Mh^{a}_{d}(f(d)): d\in {\rm supp}(f^{c})\}
=
\bigcap\{Mh^{a}_{d}(f(d)): c\leq d\in {\rm supp}(f)\}.
\]
In particular
$Mh^{a}_{0}(g_{c})=\bigcap\{Mh^{a}_{d}(g(d)): d\in {\rm supp}(g_{c})\}
=\bigcap\{Mh^{a}_{d}(g(d)): c> d\in {\rm supp}(g)\}$.
When $f=\emptyset$ or $f^{c}=\emptyset$, let $Mh^{a}_{c}(\emptyset):=\mathbb{K}$.

\item\label{df:Cpsiregularsm.3}
 (Definition of $\psi_{\pi}^{f}(a)$)
\\
 Let $a<\varepsilon(\mathbb{K})$ 
 be an ordinal, $\pi$ a regular ordinal and
 $f$ a finite function.
Then let
{\small
\begin{equation}\label{eq:Psivec}
\psi_{\pi}^{f}(a)
 :=  \min(\{\pi\}\cup\{\kappa\in Mh^{a}_{0}(f)\cap\pi:   \mathcal{H}_{a}(\kappa)\cap\pi\subset\kappa ,
   \{\pi,a\}\cup SC(f)\subset\mathcal{H}_{a}(\kappa)
\})
\end{equation}
}
For the empty function $\emptyset$,
$\psi_{\pi}(a):=\psi_{\pi}^{\emptyset}(a)$.

\item
For classes $A\subset(\mathbb{S}+1)$, let
$\alpha\in M^{a}_{c}(A)$ iff $\alpha\in A$ and
\begin{equation}\label{eq:Mca}
\forall g
[
\alpha\in Mh_{0}^{a}(g_{c}) \,\&\, SC(g_{c})\subset\mathcal{H}_{a}(\alpha) \Rightarrow
\alpha\in M\left( Mh_{0}^{a}(g_{c}) \cap A \right)
]
\end{equation}
\end{enumerate}
}

\edf

\bprp\label{prp:Mhless}
Assume $\pi\in Mh^{a}_{c}(\zeta)$ and $\xi<\zeta$ with $SC_{\Lambda}(\xi)\subset\mathcal{H}_{a}(\pi)$.
Then 
$\pi\in Mh^{a}_{c}(\xi)\cap M^{a}_{c}(Mh^{a}_{c}(\xi))$.
\eprp
\bprf
Proposition \ref{prp:idless} yields
$\pi\in Mh^{a}_{c}(\xi)$. 
$\pi\in M^{a}_{c}(Mh^{a}_{c}(\xi))$ is seen from the function $f$ such that
$f<^{c}\zeta$ with ${\rm supp}(f)=\{c\}$ and $f(c)=\xi$.
\eprf

\bprp\label{prp:MMh}
Suppose $\pi\in Mh^{a}_{c}(\xi)$.
\begin{enumerate}
\item\label{prp:MMh.1}
Let $f<^{c}\xi$ with $SC_{\Lambda}(f) \subset\mathcal{H}_{a}(\pi)$.
Then
$\pi\in M_{c}^{a}(Mh^{a}_{c}(f^{c}))$.

\item\label{prp:MMh.2}
Let $\pi\in M^{a}_{d}(A)$ for $d>c$ and $A\subset\mathbb{S}$.
Then $\pi\in M_{c}^{a}(Mh^{a}_{c}(\xi)\cap A)$.
\end{enumerate}
\eprp
\bprf
\ref{prp:MMh}.\ref{prp:MMh.1}.
Let $g$ be a function such that $\pi\in Mh_{0}^{a}(g_{c})$ with $SC_{\Lambda}(g_{c})\subset\mathcal{H}_{a}(\pi)$.
By the definition (\ref{eq:dfMhkh}) of $\pi\in Mh^{a}_{c}(\xi)$ we obtain
$\pi\in M\left( Mh_{0}^{a}(g_{c}) \cap Mh^{a}_{c}(f^{c}) \right)$.
\\
\ref{prp:MMh}.\ref{prp:MMh.2}.
Let $\pi\in M^{a}_{d}(A)$ for $d>c$. Then $\pi\in Mh^{a}_{c}(\xi)\cap A$.
Let $g$ be a function such that $\pi\in Mh_{0}^{a}(g_{c})$ with $SC_{\Lambda}(g_{c})\subset\mathcal{H}_{a}(\pi)$.
We obtain by (\ref{eq:Mca}) and $d>c$ with the function $g_{c}*h$,
$\pi\in M\left( Mh_{0}^{a}(g_{c}) \cap Mh^{a}_{c}(\xi)\cap A \right)$, where
${\rm supp}(h)=\{c\}$ and $h(c)=\xi$.
\eprf
\\

$T$ denotes the extension of ${\sf ZFC}$ by the axiom stating that
$\mathbb{S}$ is a weakly inaccessible cardinal.

\bprp\label{prp:definability}
Each of 
$x\in\mathcal{H}_{a}(y)\, (a<\varepsilon(\mathbb{K}),y<\Gamma(\mathbb{K}))$,
$x\in Mh^{a}_{c}(f)\, (c<\Lambda)$ and $x=\psi^{f}_{\kappa}(a)$
is a $\Delta_{1}$-predicate in $T$.
\eprp
\bprf 
An inspection of Definition \ref{df:Cpsiregularsm} shows that
$x\in\mathcal{H}_{a}(y)$, $\psi^{f}_{\kappa}(a)$ and $x\in Mh^{a}_{c}(f)$ 
are simultaneously defined by recursion on 
$a<\varepsilon(\mathbb{K})$, in which $x\in Mh^{a}_{c}(f)$ is defined by recursion on ordinals $x\leq\mathbb{S}$
with subsidiary recursion on $c<\Lambda$.
\eprf
\\

\textit{Shrewd cardinals} are introduced by M. Rathjen\cite{RathjenAFML2}.
A cardinal $\kappa$ is \textit{shrewd} iff for any $\eta>0$, $P\subset V_{\kappa}$,
and formula $\varphi(x,y)$, if
$V_{\kappa+\eta}\models\varphi[P,\kappa]$, then there are $0<\kappa_{0},\eta_{0}<\kappa$ such that
$V_{\kappa_{0}+\eta_{0}}\models\varphi[P\cap V_{\kappa_{0}},\kappa_{0}]$.

Let $\tilde{T}$ denote the extension of $T$ by the axiom stating that
$\mathbb{S}$ is a shrewd cardinal.

\blem\label{lem:welldefinedness.1}
$\tilde{T}$ proves that
$\mathbb{S}\in Mh^{a}_{c}(\xi)\cap M(Mh^{a}_{c}(\xi))$
for every $a<\varepsilon(\mathbb{K})$, $c,\xi<\mathbb{K}$ such that
$\{a,c,\xi\}\subset\mathcal{H}_{a}(\mathbb{S})$.
\elem 
\bprf
We show the lemma by induction on $\xi<\mathbb{K}$.

Let $\{a,c,\xi\}\cup SC_{\Lambda}(f)\subset\mathcal{H}_{a}(\mathbb{S})$ and $f<^{c}\xi$.
We show $\mathbb{S}\in M^{a}_{c}(Mh^{a}_{c}(f^{c}))$, 
which yields $\mathbb{S}\in Mh^{a}_{c}(\xi)$.
For each $d\in{\rm supp}(f^{c})$ we obtain $f(d)<\xi$ by 
$\tilde{\theta}_{-e}(\zeta)\leq\zeta$.
IH yields $\mathbb{S}\in Mh^{a}_{c}(f^{c})$.
By the definition (\ref{eq:Mca}) it suffices to show that
$
\forall g
[
\mathbb{S}\in Mh_{0}^{a}(g_{c}) \,\&\, SC_{\Lambda}(g_{c})\subset\mathcal{H}_{a}(\mathbb{S}) \Rightarrow
\mathbb{S}\in M\left( Mh_{0}^{a}(g_{c}) \cap Mh^{a}_{c}(f^{c}) \right)
]
$.

Let $g$ be a function such that $SC_{\Lambda}(g_{c})\subset\mathcal{H}_{a}(\mathbb{S})$ and
$\mathbb{S}\in Mh_{0}^{a}(g_{c})$.
We have to show $\mathbb{S}\in M(A\cap B)$ for $A=Mh_{0}^{a}(g_{c})\cap\mathbb{S}$ and 
$B=Mh^{a}_{c}(f^{c})\cap\mathbb{S}$.
Let $C$ be a club subset of $\mathbb{S}$.

We have $\mathbb{S}\in Mh_{0}^{a}(g_{c})\cap Mh^{a}_{c}(f^{c})$,
and
$\{a,c\}\cup SC_{\Lambda}(g_{c},f^{c})\subset\mathcal{H}_{a}(\mathbb{S})$.
Pick a $b<\mathbb{S}$ so that $\{a,c\}\cup SC_{\Lambda}(g_{c},f^{c})\subset\mathcal{H}_{a}(b)$.
Since the cardinality of the set $\mathcal{H}_{a}(\mathbb{S})$ is equal to $\mathbb{S}$, pick a bijection 
$F:\mathbb{S}\to \mathcal{H}_{a}(\mathbb{S})$.
Each $\alpha
\in\mathcal{H}_{a}(\mathbb{S})\cap\Gamma(\mathbb{K})$ is identified with its code, denoted
by $F^{-1}(\alpha)$.
Let $P$ be the class
$P=\{(\pi,d,\alpha)\in\mathbb{S}^{3} : \pi\in Mh^{a}_{F(d)}(F(\alpha))\}$,
where $F(d), F(\xi)<\mathbb{K}$ with $\{F(d),F(\alpha)\}\subset\mathcal{H}_{a}(\pi)$.
For fixed $a$, the set
$\{(d,\eta)\in \mathbb{K}\times\mathbb{K}   : \mathbb{S}\in Mh^{a}_{d}(\eta)\}$
is defined from the class $P$ by recursion on ordinals $d<\mathbb{K}$.

Let $\varphi$ be a formula such that $V_{\mathbb{S}+\mathbb{K}}\models\varphi[P, C,\mathbb{S},b]$ iff
$\mathbb{S}\in Mh_{0}^{a}(g_{c})\cap Mh^{a}_{c}(f^{c})$ and $C$ is a club subset of $\mathbb{S}$.
Since $\mathbb{S}$ is shrewd, pick $b<\mathbb{S}_{0}<\mathbb{K}_{0}<\mathbb{S}$ such that
$V_{\mathbb{S}_{0}+\mathbb{K}_{0}}\models\varphi[P\cap \mathbb{S}_{0},C\cap\mathbb{S}_{0},\mathbb{S}_{0},b]$.
We obtain $\mathbb{S}_{0}\in A\cap B\cap C$.

Therefore $\mathbb{S}\in Mh^{a}_{c}(\xi)$ is shown.
$\mathbb{S}\in M(Mh^{a}_{c}(\xi))$ is seen from the shrewdness of 
$\mathbb{S}$.
\eprf

\bcor\label{cor:welldefinedness.1}
$\tilde{T}$ proves that
$\forall a<\varepsilon(\mathbb{K}) \forall c<\mathbb{K}[\{a,c,\xi\}\subset\mathcal{H}_{a}(\mathbb{S}) \to \psi_{\mathbb{S}}^{f}(a)<\mathbb{S})]$
for every $\xi<\mathbb{K}$ and finite functions 
$f$ such that ${\rm supp}(f)=\{c\}$, $c<\mathbb{K}$ and
$f(c)=\xi$.
\ecor
\bprf
By Lemma \ref{lem:welldefinedness.1} we obtain
$\mathbb{S}\in M(Mh^{a}_{c}(\xi))$.
Now suppose $\{a,c,\xi\}\subset\mathcal{H}_{a}(\mathbb{S})$.
The set
$C=\{\kappa<\mathbb{S}: 
\mathcal{H}_{a}(\kappa)\cap\mathbb{S}\subset\kappa, 
\{a,c,\xi\}\subset\mathcal{H}_{a}(\kappa)\}$
 is a club subset of the regular cardinal $\mathbb{S}$,
 and $Mh^{a}_{c}(\xi)$ is stationary in $\mathbb{S}$.
This shows the existence of a $\kappa\in Mh_{c}^{a}(\xi)\cap C\cap\mathbb{S}$, and hence
$\psi_{\mathbb{S}}^{f}(a)<\mathbb{S}$ by the definition (\ref{eq:Psivec}).
\eprf

\subsection{$\psi$-functions}\label{subsec:psif}

\blem\label{lem:stepdown}
Assume $\mathbb{S}\geq\pi\in Mh^{a}_{d}(\xi)\cap Mh^{a}_{c}(\xi_{0})$, $\xi_{0}\neq 0$,
and $d<c$.
Moreover let $\xi_{1}\in\mathcal{H}_{a}(\pi)$ for 
$\xi_{1}\leq\tilde{\theta}_{c-d}(\xi_{0})$,
and $tl(\xi)\geq \xi_{1}$ when $\xi\neq 0$.
Then
$\pi\in Mh^{a}_{d}(\xi+\xi_{1})\cap M^{a}_{d}(Mh^{a}_{d}(\xi+\xi_{1}))$.
\elem
\bprf
$\pi\in M^{a}_{d}(Mh^{a}_{d}(\xi+\zeta))$ follows from $\pi\in Mh^{a}_{d}(\xi+\zeta)$
 and $\pi\in Mh^{a}_{c}(\xi_{0})\subset M^{a}_{c}(Mh^{a}_{c}(\emptyset))$
 by Proposition \ref{prp:MMh}.\ref{prp:MMh.1}.

Let $f$ be a finite function such that 
$SC_{\Lambda}(f)\subset\mathcal{H}_{a}(\pi)$, and
$f<^{d}\xi+\xi_{1}$.
We show $\pi\in M^{a}_{d}(Mh^{a}_{d}(f^{d}))$ by main induction on 
the cardinality of the finite set $\{e\in {\rm supp}(f): e>d\}$
 with subsidiary induction on $\xi_{1}$.

First let $f<^{d}\mu$ for a part $\mu$ of $\xi$. By Proposition \ref{prp:Mhless} we obtain
$\pi\in Mh^{a}_{d}(\mu)$ and $\pi\in M^{a}_{d}(Mh^{a}_{d}(f^{d}))$.

In what follows let $f(d)=\xi+\zeta$ with $\zeta<\xi_{1}$.
By SIH we obtain $\pi\in Mh^{a}_{d}(f(d))\cap M^{a}_{d}(Mh^{a}_{d}(f(d)))$.
If $\{e\in {\rm supp}(f): e>d\}=\emptyset$, then $Mh^{a}_{d}(f^{d})=Mh^{a}_{d}(f(d))$, and we are done.
Otherwise let $e=\min\{e\in {\rm supp}(f): e>d\}$.
By SIH we can assume $f<^{e}\tilde{\theta}_{-(e-d)}(tl(\xi_{1}))$.
By $\xi_{1}\leq\tilde{\theta}_{c-d}(\xi_{0})$, Propositions \ref{prp:idless} and 
\ref{prp:tht4}.\ref{prp:tht4.1},
we obtain $f<^{e}\tilde{\theta}_{-(e-d)}(\tilde{\theta}_{c-d}(\xi_{0}))=\tilde{\theta}_{-e}(\tilde{\theta}_{c}(\xi_{0}))$.
We claim that $\pi\in M^{a}_{c_{0}}(Mh_{c_{0}}^{a}(f^{c_{0}}))$ for $c_{0}=\min\{c,e\}$.
If $c=e$, then the claim follows from the assumption $\pi\in Mh_{c}^{a}(\xi_{0})$ and $f<^{e}\xi_{0}$.
Let $e=c+e_{0}>c$. Then $\tilde{\theta}_{-e}(\tilde{\theta}_{c}(\xi_{0}))=\tilde{\theta}_{-e_{0}}(hd(\xi_{0}))$, and
$f<^{c}\xi_{0}$ with $f(c)=0$ yields the claim.
Let $c=e+c_{1}>e$. Then $\tilde{\theta}_{-e}(\tilde{\theta}_{c}(\xi_{0}))=\tilde{\theta}_{c_{1}}(\xi_{0})$.
MIH yields the claim.

On the other hand we have
$Mh_{d}^{a}(f^{d})=Mh^{a}_{d}(f(d))\cap Mh_{c_{0}}^{a}(f^{c_{0}})$.
$\pi\in Mh^{a}_{d}(f(d))\cap M^{a}_{c_{0}}(Mh_{c_{0}}^{a}(f^{c_{0}}))$ with $d<c_{0}$ yields by 
Proposition \ref{prp:MMh}.\ref{prp:MMh.2},
$\pi\in M^{a}_{d}(Mh^{a}_{d}(f(d))\cap Mh_{c_{0}}^{a}(f^{c_{0}}))$, i.e.,
$\pi\in M^{a}_{d}(Mh^{a}_{d}(f^{d}))$.
\eprf

\bdf\label{df:nfform2}
{\rm
For finite functions $f$ and $g$,
\[
 Mh^{a}_{0}(g)\prec Mh^{a}_{0}(f)
 :\Leftrightarrow 
\forall\pi\in Mh^{a}_{0}(f)
\left(
SC_{\Lambda}(g)\subset\mathcal{H}_{a}(\pi) \Rightarrow \pi\in M(Mh^{a}_{0}(g))
\right)
.
\]
}
\edf

\bcor\label{cor:stepdown}
Let $f,g$ be finite functions and $c\in{\rm supp}(f)$.
Assume  that 
there exists an ordinal
$d<c$ 
such that
$(d,c)\cap {\rm supp}(f)=(d,c)\cap {\rm supp}(g)=\emptyset$, 
$g_{d}=f_{d}$, 
$g(d)<f(d)+\tilde{\theta}_{c-d}(f(c))\cdot\omega$,
and
$g<^{c}f(c)$.

Then
$Mh^{a}_{0}(g)\prec Mh^{a}_{0}(f)$ holds.
In particular if $\pi\in Mh^{a}_{0}(f)$ and
$SC_{\Lambda}(g)\subset\mathcal{H}_{a}(\pi)$, then
$\psi_{\pi}^{g}(a)<\pi$.
\ecor
\bprf
Let $\pi\in Mh^{a}_{0}(f)=\bigcap\{Mh^{a}_{e}(f(e)): e\in {\rm supp}(f)\}$ and
$SC_{\Lambda}(g)\subset\mathcal{H}_{a}(\pi)$.
Lemma \ref{lem:stepdown} with 
$\pi\in Mh^{a}_{d}(f(d))\cap Mh^{a}_{c}(f(c))$ yields
$\pi\in Mh^{a}_{d}(g(d))\cap M^{a}_{c}(Mh^{a}_{c}(g^{c}))$.
On the other hand we have $\pi\in Mh^{a}_{0}(g_{d})=\bigcap\{Mh^{a}_{e}(f(e)): e\in {\rm supp}(f)\cap d\}$.
Hence $\pi\in M(Mh^{a}_{0}(g))$.

Now suppose $SC_{\Lambda}(g)\subset\mathcal{H}_{a}(\pi)$.
The set
$C=\{\kappa<\pi: \mathcal{H}_{a}(\kappa)\cap\pi\subset\kappa, \{\pi,a\}\cup SC_{\Lambda}(g)\subset\mathcal{H}_{a}(\kappa)\}$
 is a club subset of the regular cardinal $\pi$,
 and $Mh^{a}_{0}(g)$ is stationary in $\pi$.
This shows the existence of a $\kappa\in Mh_{0}^{a}(g)\cap C\cap\pi$, and hence
$\psi_{\pi}^{g}(a)<\pi$ by the definition (\ref{eq:Psivec}).
\eprf
\\

Assume that $g_{c}=f_{c}$ and $g<^{c}f(c)$ for a $c>0$. 
Then there exists a $d<c$ for which the assumption in Corollary \ref{cor:stepdown}
is met.

\subsection{Normal forms in ordinal notations}\label{subsec:normalf}
In this subsection we introduce an \textit{irreducibility} of finite functions,
which is needed to define a normal form in ordinal notations.

\bprp\label{prp:nfform}
Let $f$ be a finite function such that
$\{a\}\cup SC_{\Lambda}(f)\subset\mathcal{H}_{a}(\pi)$.
Assume $tl(f(c))\leq\tilde{\theta}_{d}(f(c+d))$ 
holds for some $c, c+d\in {\rm supp}(f)$ with 
$d>0$.
Then 
$\pi\in Mh^{a}_{0}(f) \Leftrightarrow \pi\in Mh^{a}_{0}(g)$ holds,
where
$g$ is a finite function such that
$g(c)=f(c)-tl(f(c))$ and $g(e)=f(e)$ for every $e\neq c$.
\eprp
\bprf
First assume 
$\pi\in Mh^{a}_{0}(f)$.
We obtain 
$\pi\in Mh^{a}_{0}(g)$ by Proposition \ref{prp:idless}.
Let $\pi\in Mh^{a}_{0}(g)$, and $tl(f(c))\leq\tilde{\theta}_{d}(f(c+d))$.
On the other hand we have $\pi\in Mh^{a}_{c+d}(f(c+d))$.
By Lemma \ref{lem:stepdown} and $\pi\in Mh^{a}_{c}(g(c))$
we obtain $\pi\in Mh^{a}_{c}(f(c))$ for $f(c)=g(c)+tl(f(c))$.
Hence $\pi\in Mh^{a}_{0}(f)$.
\eprf

\bdf\label{df:irreducible}
{\rm

An \textit{irreducibility} of finite functions $f$ is defined by induction on the cardinality
$n$ of the finite set ${\rm supp}(f)$.
If $n\leq 1$, $f$ is defined to be irreducible.
Let $n\geq 2$ and $c<c+d$ be the largest two elements in ${\rm supp}(f)$, and let $g$ be 
a finite function
such that ${\rm supp}(g)={\rm supp}(f_{c})\cup\{c\}$, $g_{c}=f_{c}$ and
$g(c)=f(c)+\tilde{\theta}_{d}(f(c+d))$.

Then $f$ is irreducible iff 
$tl(f(c))>\tilde{\theta}_{d}(f(c+d))$ and
$g$ is irreducible.

}
\edf

\bdf\label{df:lx}
 {\rm 
 Let  $f,g$ be irreducible finite functions, and $b$ an ordinal.
Let us define a relation $f<^{b}_{lx}g$
by induction on the cardinality $\#\{e\in{\rm supp}(f)\cup{\rm supp}(g): e\geq b\}$ as follows.
$f<^{b}_{lx}g$ holds iff $f^{b}\neq g^{b}$ and
for the ordinal $c=\min\{c\geq b : f(c)\neq g(c)\}$,
one of the following conditions is met:

\begin{enumerate}

\item\label{df:lx.23}
$f(c)<g(c)$ and let $\mu$ be the shortest part of $g(c)$ such that $f(c)<\mu$.
Then for any $c<c+d\in{\rm supp}(f)$,  
if $tl(\mu)\leq\tilde{\theta}_{d}(f(c+d))$, then 
$f<_{lx}^{c+d}g$ holds.

\item\label{df:lx.24}
$f(c)>g(c)$ and let $\nu$ be the shortest part of $f(c)$ such that $\nu>g(c)$.
Then there exist a $c<c+d\in {\rm supp}(g)$ such that
$f<_{lx}^{c+d}g$ and
$tl(\nu)\leq \tilde{\theta}_{d}(g(c+d))$.

\end{enumerate}

}
\edf

In \cite{singledistwfprf} 
the following Lemma \ref{lem:psinucomparison} is shown.
\blem\label{lem:psinucomparison}
If $f<^{0}_{lx}g$, then
$Mh^{a}_{0}(f)\prec Mh^{a}_{0}(g)$, cf.\,Definition \ref{df:nfform2}.
\elem

\bprp\label{prp:psicomparison}
Let $f,g$ be 
irreducible finite functions, and assume that
$\psi_{\pi}^{f}(b)<\pi$ and $\psi_{\kappa}^{g}(a)<\kappa$.
Then $\psi_{\pi}^{f}(b)<\psi_{\kappa}^{g}(a)$ iff one of the following cases holds:
\begin{enumerate}
\item\label{prp:psicomparison.0}
$\pi\leq \psi_{\kappa}^{g}(a)$.

\item\label{prp:psicomparison.1}
$b<a$, $\psi_{\pi}^{f}(b)<\kappa$ and 
$SC_{\Lambda}(f)\cup\{\pi,b\}\subset\mathcal{H}_{a}(\psi_{\kappa}^{g}(a))$.

\item\label{prp:psicomparison.2}
$b>a$ and $SC_{\Lambda}(g)\cup\{\kappa,a\}\not\subset\mathcal{H}_{b}(\psi_{\pi}^{f}(b))$.

\item\label{prp:psicomparison.25}
$b=a$, $\kappa<\pi$ and $\kappa\not\in\mathcal{H}_{b}(\psi_{\pi}^{f}(b))$.

\item\label{prp:psicomparison.3}
$b=a$, $\pi=\kappa$, $SC_{\Lambda}(f)\subset\mathcal{H}_{a}(\psi_{\kappa}^{g}(a))$, and
$f<^{0}_{lx}g$.

\item\label{prp:psicomparison.4}
$b=a$, $\pi=\kappa$, 
$SC_{\Lambda}(g)\not\subset\mathcal{H}_{b}(\psi_{\pi}^{f}(b))$.

\end{enumerate}

\eprp
\bprf
This is seen as in Proposition 2.19 of \cite{KPPiN} 
using Lemma \ref{lem:psinucomparison}.
\eprf
\\

In \cite{singledistwfprf} 
a computable notation system $OT_{N}$
is defined from Proposition \ref{prp:psicomparison} for each positive integer $N$.
Constants are $0$ and $\mathbb{S}$, and constructors are $+,\varphi, \Omega$ and $\psi$.
$\Omega$-terms $\Omega_{\alpha}\in OT_{N}$ if
$\alpha\in\{1\}\cup\{\kappa+n: 0<n\leq N\}$ for $\kappa\in\{\mathbb{S}\}\cup\Psi_{N}$
with $\Omega:=\Omega_{1}$.
Let us spell out clauses for $\psi$-terms, the set $\Psi_{N}$ and $m(\rho)$.

\bdf\label{df:notationsystem}
{\rm
$E_{\mathbb{S}}(\alpha)$ denotes the set of subterms $\beta$ of ordinal terms $\alpha$
such that $\beta<\mathbb{S}$, and
the length $\ell\alpha$  of $\alpha$ is the total number of occurrences of symbols in $\alpha$.
\begin{enumerate}

 \item\label{df:notationsystem.9}
 
Let $\alpha=\psi_{\pi}(a)$ with $\{\pi,a\}\subset OT_{N}$, 
$\{\pi,a\}\subset\mathcal{H}_{a}(\alpha)$,
and if $\pi=\Omega_{\kappa+n}$ with $\kappa\in\Psi_{N}$ and $0<k\leq N$, then 
$a<\Gamma(\Omega_{\kappa+N})$.
Then
$\alpha\in OT_{N}$.
Let
$m(\alpha)=\emptyset$.

 \item\label{df:notationsystem.10}
Let $\alpha=\psi_{\mathbb{S}}^{f}(a)$, where $\xi,a,c\in OT_{N}$,
$\xi>0$, $c<\mathbb{K}$,
$\{\xi,a,c\}\subset\mathcal{H}_{a}(\alpha)$, ${\rm supp}(f)=\{c\}$ and $f(c)=\xi$.
Then
$\alpha\in \Psi_{N}$ and $f=m(\alpha)$.

 \item\label{df:notationsystem.11}
Let $\{a,d\}\subset OT_{N}$, $\pi\in\Psi_{N}$,
$f=m(\pi)$,
 $d<c\in {\rm supp}(f)$,
and $(d,c)\cap {\rm supp}(f)=\emptyset$.
Let $g$ be an irreducible function such that 
$SC_{\Lambda}(g)\subset OT_{N}$,
$g_{d}=f_{d}$, $(d,c)\cap {\rm supp}(g)=\emptyset$,
$g(d)<f(d)+\tilde{\theta}_{c-d}(f(c))\cdot\omega$, 
and $g<^{c}f(c)$.

Then 
$\alpha=\psi_{\pi}^{g}(a)\in \Psi_{N}$ 
if 
$\{\pi,a\}\cup SC_{\Lambda}(f)\cup SC_{\Lambda}(g)\subset\mathcal{H}_{a}(\alpha)$ and

\begin{equation}\label{eq:notationsystem.11}
E_{\mathbb{S}}(SC_{\Lambda}(g))<\alpha
\end{equation}

\end{enumerate}
}
\edf
The restriction (\ref{eq:notationsystem.11}) is needed to prove the well-foundednesss 
of $OT_{N}$ in \cite{singledistwfprf}.

In what follows by ordinals we mean ordinal terms in $OT_{N}$ for a fixed positive integer $N$.

\subsection{A Mostowski collapsing}

In this subsection we define a Mostowski collapsing $\alpha\mapsto\alpha[\rho/\mathbb{S}]$,
which is needed to replace inference rules for stability by ones of reflections.

\bprp\label{prp:EmS}
Let $\rho=\psi_{\kappa}^{f}(a)\in \Psi_{N}\cap\kappa$ with 
$\mathcal{H}_{\gamma}(\kappa)\cap\mathbb{S}\subset\kappa$ for $\gamma\leq a$.
Then $\mathcal{H}_{\gamma}(\rho)\cap\mathbb{S}\subset\rho$.
\eprp
\bprf
If $\kappa=\mathbb{S}$, then $\mathcal{H}_{\gamma}(\rho)\cap\mathbb{S}\subset \mathcal{H}_{a}(\rho)\cap\mathbb{S}\subset\rho$
for $\gamma\leq a$.
Let $\kappa=\psi_{\pi}^{g}(b)<\mathbb{S}$. 
We have $\kappa\in \mathcal{H}_{a}(\rho)$ by (\ref{eq:Psivec}),
and hence $b<a$ by $\kappa>\rho$.
We obtain
$\mathcal{H}_{\gamma}(\rho)\cap\mathbb{S}\subset \mathcal{H}_{\gamma}(\kappa)\cap\mathbb{S}\subset\kappa$.
$\gamma\leq a$
yields
$\mathcal{H}_{\gamma}(\rho)\cap\mathbb{S}\subset \mathcal{H}_{\gamma}(\rho)\cap\kappa\subset
\mathcal{H}_{a}(\rho)\cap\kappa\subset\rho$.
\eprf

\bdf{\rm
Let
$\alpha\in E^{\mathbb{S}}_{\rho}$ iff $E_{\mathbb{S}}(\alpha)\subset\rho$ for $\alpha\in OT_{N}$.
}
\edf

\bprp\label{prp:Erho}
$SC_{\Lambda}(m(\rho))\subset E^{\mathbb{S}}_{\rho}$ for $\rho\in\Psi_{N}$.
\eprp
\bprf
This is seen from Definition \ref{df:notationsystem}.\ref{df:notationsystem.10} and
 (\ref{eq:notationsystem.11})
in Definition \ref{df:notationsystem}.\ref{df:notationsystem.11}.
\eprf

\bprp\label{prp:EK}
Let $\rho\leq\delta<\mathbb{S}$. Then for $\alpha\in E^{\mathbb{S}}_{\rho}\cap OT_{N}$,
$\alpha\in \mathcal{H}_{\gamma}(\rho)$ iff $\alpha\in \mathcal{H}_{\gamma}(\delta)$.
\eprp
\bprf
We show
$\alpha\in \mathcal{H}_{\gamma}(\rho)$ iff $\alpha\in \mathcal{H}_{\gamma}(\delta)$ by induction on the lengths $\ell\alpha$
of $\alpha\in E^{\mathbb{S}}_{\rho}\cap OT_{N}$.
First let $\alpha=\psi_{\sigma}(a)>\mathbb{S}$. Then $\sigma=\Omega_{\mathbb{S}+n}$ and for $\mathbb{S}>\beta\in\{\rho,\delta\}$,
$\alpha\in \mathcal{H}_{\gamma}(\beta)$ iff $a\in \mathcal{H}_{\gamma}(\beta)\cap\gamma$, and
$E_{\mathbb{S}}(\alpha)=E_{\mathbb{S}}(a)$.
IH yields $a\in \mathcal{H}_{\gamma}(\rho)$ iff $a\in \mathcal{H}_{\gamma}(\delta)$.
Next let $\alpha=\psi_{\pi}^{f}(a)<\mathbb{S}$. Then $\alpha<\rho\leq\delta$ by $\alpha\in E^{\mathbb{S}}_{\rho}$.
Hence $\alpha\in \mathcal{H}_{\gamma}(\rho)\cap \mathcal{H}_{\gamma}(\delta)$.
Other cases are seen from IH.
\eprf

\bdf\label{df:Mostwskicollaps}
{\rm
Let $\alpha\in E^{\mathbb{S}}_{\rho}$ with $\rho\in \Psi_{N}$.
We define an ordinal $\alpha[\rho/\mathbb{S}]$ recursively as follows.

\begin{enumerate}
\item
$\alpha[\rho/\mathbb{S}]:=\alpha$ when $\alpha<\mathbb{S}$.
In what follows assume $\alpha\geq\mathbb{S}$.

\item
$\mathbb{S}[\rho/\mathbb{S}]:=\rho$.
$\Omega_{\mathbb{S}+n}[\rho/\mathbb{S}]:=\Omega_{\rho+n}$.

\item
The map commutes with $+$ and $\varphi$, i.e.,
$(\varphi\alpha\beta)[\rho/\mathbb{S}]=\varphi(\alpha[\rho/\mathbb{S}])(\beta[\rho/\mathbb{S}])$,
 and
$(\alpha_{1}+\cdots+\alpha_{n})[\rho/\mathbb{S}]=\alpha_{1}[\rho/\mathbb{S}]+\cdots+\alpha_{n}[\rho/\mathbb{S}]$.

\item
$\left(\psi_{\Omega_{\mathbb{S}+n}}(a)\right)[\rho/\mathbb{S}]=\psi_{\Omega_{\rho+n}}(a[\rho/\mathbb{S}])$ for $0<n\leq N$.

\end{enumerate}
}
\edf

Note that $\alpha[\rho/\mathbb{S}]\in\mathcal{H}_{\mathbb{S}}(E_{\mathbb{S}}(\alpha)\cup\{\rho\})$.

\blem\label{lem:Mostowskicollaps}
For $\rho\in\Psi_{N}$,
$\{\alpha[\rho/\mathbb{S}]:\alpha\in E^{\mathbb{S}}_{\rho}\}$ is a transitive collapse of $E^{\mathbb{S}}_{\rho}$:
$\beta<\alpha\Leftrightarrow\beta[\rho/\mathbb{S}]<\alpha[\rho/\mathbb{S}]$,
$\gamma>\mathbb{S} \Rightarrow 
\left(
\beta\in\mathcal{H}_{\alpha}(\gamma)\Leftrightarrow 
\beta[\rho/\mathbb{S}]\in\mathcal{H}_{\alpha[\rho/\mathbb{S}]}(\gamma[\rho/\mathbb{S}]))
\right)$
and
$OT_{N}\cap\alpha[\rho/\mathbb{S}]=\{\beta[\rho/\mathbb{S}]:\beta\in E^{\mathbb{S}}_{\rho}\cap\alpha\}$
for $\alpha,\beta,\gamma\in E^{\mathbb{S}}_{\rho}$.

\elem
\bprf
Simultaneously we show first
$\beta<\alpha \Leftrightarrow \beta[\rho/\mathbb{S}]<\alpha[\rho/\mathbb{S}]$,
and second
$\beta\in\mathcal{H}_{\alpha}(\gamma)\Leftrightarrow \beta[\rho/\mathbb{S}]\in\mathcal{H}_{\alpha[\rho/\mathbb{S}]}(\gamma[\rho/\mathbb{S}]))$
if $\gamma>\mathbb{S}$
by induction on the sum $2^{\ell\alpha}+2^{\ell\beta}$ of lengths for
$\alpha,\beta,\gamma\in E^{\mathbb{S}}_{\rho}$.
We see easily that $\mathbb{S}>\Gamma_{\mathbb{K}[\rho]+1}>\alpha[\rho/\mathbb{S}]>\rho$ when $\alpha>\mathbb{S}$, where $\mathbb{K}[\rho/\mathbb{S}]=\Omega_{\rho+N}$.
Also $\alpha[\rho/\mathbb{S}]\leq\alpha$.

Let $\mathbb{S}<\beta=\psi_{\pi}(b)<\psi_{\kappa}(a)=\alpha$ with
$\{b,a\}\subset E^{\mathbb{S}}_{\rho}\cap OT_{N}$, 
where $\pi=\Omega_{\mathbb{S}+n}$, $\kappa=\Omega_{\mathbb{S}+m}$
for $0<n,m\leq N$,
$b\in\mathcal{H}_{b}(\beta)$ and $a\in\mathcal{H}_{a}(\alpha)$.
We have
$\beta=\psi_{\pi}(b)<\psi_{\kappa}(a)=\alpha$ iff either $\pi<\kappa$, i.e., $n<m$, or $\pi=\kappa\,\&\, b<a$,
and similarly for $\beta[\rho/\mathbb{S}]=\psi_{\Omega_{\rho}+n}(b[\rho/\mathbb{S}])<\psi_{\Omega_{\rho+m}}(a[\rho/\mathbb{S}])$.
From IH we see that $b\in\mathcal{H}_{b}(\beta) \Leftrightarrow b[\rho/\mathbb{S}]\in\mathcal{H}_{b[\rho]}(\beta[\rho/\mathbb{S}])$ and
$b<a \Leftrightarrow b[\rho/\mathbb{S}]<a[\rho/\mathbb{S}]$.
Hence 
$\beta<\alpha \Leftrightarrow \beta[\rho/\mathbb{S}]<\alpha[\rho/\mathbb{S}]$.
Other cases are seen by IH.
Next suppose $\beta\in\mathcal{H}_{\alpha}(\gamma)$ for $\gamma>\mathbb{S}$.
Then $\beta[\rho/\mathbb{S}]\in\mathcal{H}_{\alpha[\rho/\mathbb{S}]}(\gamma[\rho/\mathbb{S}])$ is seen from 
the first assertion using the fact $\gamma[\rho/\mathbb{S}]>\rho$.

Finally let
$\beta\in OT_{N}\cap\alpha[\rho/\mathbb{S}]$ for $\alpha\in E^{\mathbb{S}}_{\rho}$.
We show by induction on $\ell\beta$ that there exists a 
$\gamma\in E^{\mathbb{S}}_{\rho}$
such that $\beta=\gamma[\rho/\mathbb{S}]$.
If $\beta<\rho$, then $\beta[\rho/\mathbb{S}]=\beta$．
Also $\rho=\mathbb{S}[\rho/\mathbb{S}]$.
Let $\Gamma_{\mathbb{K}[\rho]+1}>\beta=\psi_{\pi}(b)>\rho$ with 
$b\in\mathcal{H}_{b}(\beta)$.
Then we see
$\pi=\Omega_{\rho+n}$ for an $n\leq N$ and 
$b<\Gamma_{\mathbb{K}[\rho]+1}$.
By IH there is a $c\in E^{\mathbb{S}}_{\rho}$ such that
$c[\rho/\mathbb{S}]=b$.
Then $\beta=\psi_{\Omega_{\rho+n}}(c[\rho/\mathbb{S}])=\gamma[\rho/\mathbb{S}]$ with
$\gamma=\psi_{\Omega_{\mathbb{S}+n}}(c)$,
$c\in\mathcal{H}_{c}(\gamma)$ and $E_{\mathbb{S}}(\gamma)=E_{\mathbb{S}}(c)$.
Hence $\gamma\in E^{\mathbb{S}}_{\rho}$.
Other cases are seen by IH.
\eprf

\bprp\label{prp:EK2}
$\mathcal{H}_{\gamma}(E_{\rho}^{\mathbb{S}})\subset E_{\rho}^{\mathbb{S}}$ for $\rho\in\Psi_{N}$.
\eprp
\bprf
Let $\mathcal{H}_{\gamma}(\rho)\cap\mathbb{S}\subset\rho$.
We show $\alpha\in E^{\mathbb{S}}_{\rho}$ by induction on $\ell\alpha$ for
$\alpha\in \mathcal{H}_{\gamma}(E^{\mathbb{S}}_{\rho})\cap OT_{N}$.
Let $\{\kappa,a\}\cup SC_{\Lambda}(g)\subset\mathcal{H}_{\gamma}(E^{\mathbb{S}}_{\rho})$ be such that
$a<\gamma$ and $\{\kappa,a\}\cup SC_{\Lambda}(g)\subset\mathcal{H}_{a}(\alpha)\cap OT_{N}$ with
$\alpha=\psi_{\kappa}^{g}(a)\in\mathcal{H}_{\gamma}(E^{\mathbb{S}}_{\rho})$.
We need to show $\alpha<\rho$. Suppose $\rho\leq\alpha<\mathbb{S}$.
IH yields $\{\kappa,a\}\cup SC_{\Lambda}(g)\subset E^{\mathbb{S}}_{\rho}\cap OT_{N}$.
Proposition \ref{prp:EK} yields 
$\{\kappa,a\}\cup SC_{\Lambda}(g)\subset\mathcal{H}_{a}(\rho)\subset\mathcal{H}_{\gamma}(\rho)$
by $a<\gamma$.
Hence $\alpha\in\mathcal{H}_{\gamma}(\rho)\cap\mathbb{S}\subset\rho$.
Other cases are seen from IH.
\eprf

\section{Upperbounds}\label{sect:controlledOme}

Operator controlled derivations are introduced by W. Buchholz\cite{Buchholz}.
In this section except otherwise stated,
 $\alpha,\beta,\gamma,\ldots,a,b,c,d,\ldots$ range over ordinals in $OT_{N}$,
$\xi,\zeta,\nu,\mu,\ldots$ range over ordinals in $\mathbb{K}$,
$f,g,h,\ldots$ range over finite functions on $\mathbb{K}$,
and $\pi,\kappa,\rho,\sigma,\tau,\lambda,\ldots$ range over regular ordinals in $OT_{N}$.
$Reg$ denotes the set of regular ordinals$\leq\mathbb{K}=\Omega_{\mathbb{S}+N}$
with a positive integer $N$.

\subsection{Classes of sentences}\label{subsec:classformula}

Following Buchholz\cite{Buchholz} 
let us introduce a language for ramified set theory $RS$.

\bdf
{\rm
$RS$\textit{-terms} and their \textit{levels} are inductively defined.
\begin{enumerate}
\item
For each $\alpha\in OT_{N}\cap \mathbb{K}$, $\mathsf{L}_{\alpha}$ is an $RS$-term of level $\alpha$.

\item
For a set-theoretic formula $\phi(x,y_{1},\ldots,y_{n})$ in the language $\{\in\}$ and 
$RS$-terms $a_{1},\ldots,a_{n}$ of levels$<\! \alpha\in OT_{N}\cap\mathbb{K}$, 
$[x\in \mathsf{L}_{\alpha}:\phi^{\mathsf{L}_{\alpha}}(x,a_{1},\ldots,a_{n})]$ 
is an $RS$-term of level $\alpha$.
\end{enumerate}
}
\edf

\bdf\label{df:RS}
{\rm
\begin{enumerate}
\item
$|u|$ denotes the level of $RS$-terms $u$, and $Tm(\alpha)$ the set of $RS$-terms of level$<\alpha\in OT_{N}\cap(\mathbb{K}+1)$.
$Tm=Tm(\mathbb{K})$ is then the set of $RS$-terms, which are denoted by $u,v,w,\ldots$

\item\label{df:RS.1}
$RS$-\textit{formulas} 
are constructed from \textit{literals}
$u\in v, u\not\in v$
by propositional connectives $\lor,\land$, bounded quantifiers $\exists x\in u, \forall x\in u$ and
unbounded quantifiers $\exists x,\forall x$.
Unbounded quantifiers $\exists x,\forall x$ are denoted by $\exists x\in L_{\mathbb{K}},\forall x\in L_{\mathbb{K}}$, resp.
It is convenient for us not to restrict propositional connectives $\lor,\land$ to binary ones.
Specifically when $A_{i}$ are $RS$-formulas for $i<n<\omega$,
$A_{0}\lor\cdots\lor A_{n-1}$ and $A_{0}\land\cdots\land A_{n-1}$ are $RS$-formulas.
Even when $n=1$, $A_{0}\lor\cdots\lor A_{0}$ is understood to be different from the formula $A_{0}$.

\item\label{df:RS.2}
For $RS$-terms and $RS$-formulas $\iota$,
$\mathsf{k}(\iota)$ denotes the set of ordinal terms $\alpha$ such that the 
constant $L_{\alpha}$ occurs in $\iota$.
$|\iota|=\max(\mathsf{k}(\iota)\cup\{0\})$.

\item\label{df:RS.3}
$\Delta_{0}$-formulas, $\Sigma_{1}$-formulas and $\Sigma$-formulas are defined as in \cite{Ba}.
Specifically if $\psi$ is a $\Sigma$-formula, then so is the formula $\forall y\in z\psi$.
$\theta^{(a)}$ denotes a $\Delta_{0}$-formula obtained from a $\Sigma$-formula $\theta$
by restricting each unbounded existential quantifier to $a$.

\item\label{df:RS.4}
For a set-theoretic $\Sigma_{1}$-formula $\psi(x_{1},\ldots,x_{m})$
and 
$u_{1},\ldots,u_{m}\in Tm(\kappa)$ with $\kappa\leq\mathbb{K}$,
$\psi^{(L_{\kappa})}(u_{1},\ldots,u_{m})$
 is a $\Sigma_{1}(\kappa)$\textit{-formula}.
 $\Delta_{0}(\kappa)$-formulas and $\Sigma(\kappa)$-formulas are defined similarly

\item\label{df:RS.5}
For $\theta\equiv\psi^{(\mathsf{L}_{\kappa})}(u_{1},\ldots,u_{m})\in\Sigma(\kappa)$ 
and $\lambda<\kappa$,
$\theta^{(\lambda,\kappa)}:\equiv \psi^{(\mathsf{L}_{\lambda})}(u_{1},\ldots,u_{m})$.

\item\label{df:RS.6}
Let $\rho\leq\mathbb{S}$, and
$\iota$ an $RS$-term or an $RS$-formula such that
$\mathsf{k}(\iota)\subset E^{\mathbb{S}}_{\rho}$ with $E^{\mathbb{S}}_{\mathbb{S}}=\mathbb{K}$.
Then
$\iota^{[\rho/\mathbb{S}]}$ denotes the result of replacing each unbounded quantifier
$Qx$ by $Qx\in \mathsf{L}_{\mathbb{K}[\rho/\mathbb{S}]}$,
and each ordinal term $\alpha\in \mathsf{k}(\iota)$
by $\alpha[\rho/\mathbb{S}]$ for the Mostowski collapse in
Definition \ref{df:Mostwskicollaps}.
$\iota^{[\rho/\mathbb{S}]}$ is defined recursively as follows.

 \begin{enumerate}
 \item
 $(\mathsf{L}_{\alpha})^{[\rho/\mathbb{S}]}\equiv \mathsf{L}_{\alpha[\rho/\mathbb{S}]}$ with $\alpha\in E^{\mathbb{S}}_{\rho}$.
 When $\{\alpha\}\cup\bigcup\{\mathsf{k}(u_{i}):i\leq n\}\subset E^{\mathbb{S}}_{\rho}$,
 $\left([x\in \mathsf{L}_{\alpha}:\phi^{\mathsf{L}_{\alpha}}(x,u_{1},\ldots,u_{n})]\right)^{[\rho/\mathbb{S}]}$
 is defined to be the $RS$-term
 $[x\in \mathsf{L}_{\alpha[\rho/\mathbb{S}]}:\phi^{\mathsf{L}_{\alpha[\rho/\mathbb{S}]]}}(x, (u_{1})^{[\rho/\mathbb{S}]},\ldots,(u_{n})^{[\rho/\mathbb{S}]})]$.
 
 \item
 $(\lnot A)^{[\rho/\mathbb{S}]}\equiv\lnot A^{[\rho/\mathbb{S}]}$.
 $(u\in v)^{[\rho/\mathbb{S}]}\equiv\left(u^{[\rho/\mathbb{S}]}\in v^{[\rho/\mathbb{S}]}\right)$.
 $(A_{0}\lor \cdots\lor A_{n-1})^{[\rho/\mathbb{S}]}\equiv (A_{0})^{[\rho/\mathbb{S}]}\lor\cdots\lor(A_{n-1})^{[\rho/\mathbb{S}]}$.
 $(\exists x\in u A)^{[\rho/\mathbb{S}]}\equiv (\exists x\in u^{[\rho/\mathbb{S}]} A^{[\rho/\mathbb{S}]})$.
 $(\exists x A)^{[\rho/\mathbb{S}]}\equiv(\exists x\in \mathsf{L}_{\mathbb{K}[\rho/\mathbb{S}]} A^{[\rho/\mathbb{S}]})$.
 
 \end{enumerate}

\end{enumerate}
}
\edf

\bprp\label{prp:levelcollaps}
Let $\rho\in Reg\cap(\mathbb{S}+1)$.

\begin{enumerate}
\item\label{prp:levelcollaps1}
Let $v$ be an $RS$-term with $\mathsf{k}(v)\subset E^{\mathbb{S}}_{\rho}$,
and $\alpha=|v|$.
Then 
$v^{[\rho/\mathbb{S}]}$ is an $RS$-term of level $\alpha[\rho/\mathbb{S}]$,
$\left| v^{[\rho/\mathbb{S}]} \right|=\alpha[\rho/\mathbb{S}]$ and
$\mathsf{k}(v^{[\rho/\mathbb{S}]})=\left(\mathsf{k}(v)\right)^{[\rho/\mathbb{S}]}$.

\item\label{prp:levelcollaps2}
Let $\alpha\leq\mathbb{K}$ be such that $\alpha\in E^{\mathbb{S}}_{\rho}$. Then
$
\left(Tm(\alpha)\right)^{[\rho/\mathbb{S}]}
:= 
\{v^{[\rho/\mathbb{S}]}: v\in Tm(\alpha), \mathsf{k}(v)\subset E^{\mathbb{S}}_{\rho}\}
= Tm(\alpha[\rho/\mathbb{S}])$.

\item\label{prp:levelcollaps3}
Let $A$ be an $RS$-formula with $\mathsf{k}(A)\subset \mathcal{H}_{\gamma}(\rho)$, and assume
$\mathcal{H}_{\gamma}(\rho)\cap\mathbb{S}\subset\rho$.
Then $A^{[\rho/\mathbb{S}]}$ is an $RS$-formula such that
$\mathsf{k}(A^{[\rho/\mathbb{S}]})\subset\{\alpha[\rho/\mathbb{S}]: \alpha\in\mathsf{k}(A)\}\cup\{\mathbb{K}[\rho/\mathbb{S}]\}$.
\end{enumerate}
\eprp
\bprf
\ref{prp:levelcollaps}.\ref{prp:levelcollaps1}.
We see easily that 
$v^{[\rho/\mathbb{S}]}$ is an $RS$-term of level $\alpha[\rho/\mathbb{S}]$.
\\
\ref{prp:levelcollaps}.\ref{prp:levelcollaps2}.
We see 
$\left(Tm(\alpha)\right)^{[\rho/\mathbb{S}]}\subset Tm(\alpha[\rho/\mathbb{S}])$ from 
Proposition \ref{prp:levelcollaps}.\ref{prp:levelcollaps1}.
Conversely let $u$ be an $RS$-term with $\mathsf{k}(u)=\{\beta_{i}: i<n\}$ 
and $\max\{\beta_{i}: i<n\}=|u|<\alpha[\rho/\mathbb{S}]$.
By 
Lemma \ref{lem:Mostowskicollaps}
there are ordinal terms $\gamma_{i}\in OT_{N}$ such that
$\gamma_{i}\in E^{\mathbb{S}}_{\rho}$ and
$\gamma_{i}[\rho/\mathbb{S}]=\beta_{i}$.
Let $v$ be an $RS$-term obtained from $u$ by replacing each constant $\mathsf{L}_{\beta_{i}}$ by $\mathsf{L}_{\gamma_{i}}$.
We obtain $v^{[\rho/\mathbb{S}]}\equiv u$, $v\in Tm(\alpha)$, and $\mathsf{k}(v)=\{\gamma_{i}: i<n\}\subset E^{\mathbb{S}}_{\rho}$.
This means $v\in\left(Tm(\alpha)\right)^{[\rho/\mathbb{S}]}$.
\eprf
\\

In what follows we need to consider \textit{sentences}.
Sentences are denoted $A,C$ possibly with indices.

For each sentence $A$, either a disjunction is assigned as $A\simeq\bigvee(A_{\iota})_{\iota\in J}$, or
a conjunction is assigned as $A\simeq\bigwedge(A_{\iota})_{\iota\in J}$.
In the former case $A$ is said to be a \textit{$\bigvee$-formula}, and in the latter
$A$ is a \textit{$\bigwedge$-formula}.
It is convenient for us, cf.\,Recapping \ref{mlem:singlemainl.1}, to modify
the assignment of disjunctions and conjunctions to sentences from \cite{Buchholz} 
such that
if $A\simeq\bigvee(A_{\iota})_{\iota\in J}$ is a $\bigvee$-formula, then each 
$A_{\iota}$ is a $\bigwedge$-formula, and similarly for $\bigwedge$-formula $A$.
\bdf\label{df:bigveewedge}
{\rm
If $A$ is a $\bigvee$-formula, then let $A^{\lor}:\equiv A$.
Otherwise let
$A^{\lor}:\equiv \left(\bigvee_{i\leq 0}B_{i}\right)$ with $B_{0}\equiv A$.
If $A$ is a $\bigwedge$-formula, then let $A^{\land}:\equiv A$.
Otherwise let
$A^{\land}:\equiv \left(\bigwedge_{i\leq 0}B_{i}\right)$ with $B_{0}\equiv A$.
}
\edf

\bdf\label{df:assigndc}
{\rm
Let
$[\rho]Tm(\alpha):=\{u\in Tm(\alpha) : \mathsf{k}(u)\subset E^{\mathbb{S}}_{\rho}\}$.

\begin{enumerate}
\item
For $v,u\in Tm(\mathbb{K})$ with $|v|<|u|$, let
\[
(v\dot{\in} u) :\equiv
\left\{
\begin{array}{ll}
A(v) & \mbox{{\rm if }} u\equiv[x\in \mathsf{L}_{\alpha}: A(x)]
\\
v\not\in \mathsf{L}_{0} & \mbox{{\rm if }} u\equiv \mathsf{L}_{\alpha}
\end{array}
\right.
\]
and 
$(u=v):\equiv(\forall x\in u(x\in v)\land \forall x\in v(x\in u))$.

\item\label{df:assigndc0}
For $v,u\in Tm(\mathbb{K})$, 
let $[\rho]J:=[\rho]Tm(|u|)$ with $J=Tm(|u|)$. Then
$(v\in u)  :\simeq  \bigvee(A_{w,0}\land A_{w,1}\land A_{w,2})_{w\in J}$, 
and
$(v\not\in u)  :\simeq  \bigwedge(\lnot A_{w,0}\lor \lnot A_{w,1}\lor \lnot A_{w,2})_{w\in J}$,
where
$A_{w,0}\equiv (w\varepsilon u)^{\lor}$, $A_{w,1}\equiv(\forall x\in w(x\in v))^{\lor}$ and 
$A_{w,2}\equiv(\forall x\in v(x\in w))^{\lor}$.

\item
$(A_{0}\lor\cdots\lor A_{n-1}):\simeq \bigvee(A_{\iota}^{\land})_{\iota\in J}$ 
and
$(A_{0}\land\cdots\land A_{n-1}):\simeq \bigwedge(A_{\iota}^{\lor})_{\iota\in J}$ for 
$J:=n$.

\item
For $u\in Tm(\mathbb{K})\cup\{\mathsf{L}_{\mathbb{K}}\}$,
$\exists x\in u\, A(x):\simeq \bigvee(A_{v})_{v\in J}$
and
$\forall x\in u\, \lnot A(x):\simeq \bigwedge(\lnot A_{v})_{v\in J}
$
for 
$A_{v}:\equiv ((v\dot{\in} u)^{\lor} \land (A(v))^{\lor})$,
$[\rho]J:=[\rho]Tm(|u|)$ with $J=Tm(|u|)$.

\end{enumerate}

}
\edf

\bprp\label{lem:assigncollaps}
Let $\rho\in \Psi_{N}\cup\{\mathbb{S}\}$.
For $RS$-formulas $A$, let
$A\simeq \bigvee(A_{\iota})_{\iota\in J}$ and assume $\mathsf{k}(A)\subset E^{\mathbb{S}}_{\rho}$.
Then
$A^{[\rho/\mathbb{S}]}\simeq \bigvee\left((A_{\iota})^{[\rho/\mathbb{S}]}\right)_{\iota\in [\rho]J}$.
The case $A\simeq \bigwedge(A_{\iota})_{\iota\in J}$ is similar.
\eprp
\bprf
This is seen from Proposition \ref{prp:levelcollaps}.\ref{prp:levelcollaps2}.
\eprf
\\

The rank $\mathrm{rk}(\iota)$ of 
sentences or terms $\iota$ 
is defined as in \cite{Buchholz} so that the following Proposition \ref{lem:rank} holds.
For completeness let us reproduce it.

\bdf\label{df:rank}
{\rm

$\mathrm{rk}(\lnot A):=\mathrm{rk}(A)$.
$\mathrm{rk}(\mathsf{L}_{\alpha})=\omega\alpha$.
$\mathrm{rk}([x\in \mathsf{L}_{\alpha}: A(x)])=\max\{\omega\alpha+1, \mathrm{rk}(A(\mathsf{L}_{0}))+3\}$.
$\mathrm{rk}(v\in u)=\max\{\mathrm{rk}(v)+7,\mathrm{rk}(u)+2\}$.
$\mathrm{rk}(A_{0}\lor \cdots\lor A_{n-1})=\max(\{0\}\cup \{\mathrm{rk}((A_{i})^{\land})+1:i<n\})$.
 $\mathrm{rk}(\exists x\in u\, A(x))=\max\{\mathrm{rk}(u), \mathrm{rk}(A(\mathsf{L}_{0}))+3\}$ for $u\in Tm(\mathbb{K})\cup\{\mathsf{L}_{\mathbb{K}}\}$.

For finite sets $\Delta$ of sentences, let
$\mathrm{rk}(\Delta)=\max(\{0\}\cup\{\mathrm{rk}(\delta):\delta\in\Delta\})$.
}

\edf

\bprp\label{lem:rank}
Let $A$ be a sentence with 
$A\simeq\bigvee(A_{\iota})_{\iota\in J}$ or $A\simeq\bigwedge(A_{\iota})_{\iota\in J}$.
\begin{enumerate}

\item\label{lem:rank0}
$\mathrm{rk}(A)<\mathbb{K}+\omega$.

\item\label{lem:rank1}
$\omega |u|\leq \mathrm{rk}(u)\in\{\omega |u|+i: i\in\omega\}$, and
$|A|\leq \mathrm{rk}(A)\in\{\omega |A|+i  : i\in\omega\}$.

\item\label{lem:rank2}
$\forall\iota\in J(\mathrm{rk}(A_{\iota})<\mathrm{rk}(A))$.

\item\label{lem:rank3}
Let $\alpha$ be an epsilon number. Then
$\mathrm{rk}(\exists x\in \mathsf{L}_{\alpha}B)=\alpha$ for $B\in\Delta_{0}(\alpha)$.
Conversely if $\mathrm{rk}(A)=\alpha$ for a $\bigvee$-formula $A$, then $A\in\Sigma_{1}(\alpha)$.

\item\label{lem:rank4}
Let $\rho\in Reg$ and $\mathsf{k}(\iota)\subset E^{\mathbb{S}}_{\rho}$.
Then
$\mathrm{rk}(\iota^{[\rho/\mathbb{S}]})=\left(\mathrm{rk}(\iota)\right)[\rho/\mathbb{S}]$.
\end{enumerate}

\eprp
\bprf
These are shown in \cite{Buchholz} except Proposition \ref{lem:rank}.\ref{lem:rank4}, which 
is seen from 
the facts
$(\omega\alpha)[\rho/\mathbb{S}]=\omega(\alpha[\rho/\mathbb{S}])$ and $(\alpha+1)[\rho/\mathbb{S}]=\alpha[\rho/\mathbb{S}]+1$ when $\alpha\in E^{\mathbb{S}}_{\rho}$.
We see that $\mathrm{rk}(\iota)\in E^{\mathbb{S}}_{\rho}$
from Proposition \ref{lem:rank}.\ref{lem:rank1}.
\eprf

\subsection{Sets $E_{m}(\alpha)$}

In this subsection 
sets $E_{m}(\alpha)\subset\Omega_{\mathbb{S}+N-m-1}$ of ordinals are defined for ordinals $\alpha$ so as to have 
Proposition \ref{prp:Em}.
The proposition states that for $\alpha\in E^{\mathbb{S}}_{\sigma}$,
$\alpha[\sigma/\mathbb{S}]\in\mathcal{H}_{\gamma}[E_{m}(\alpha)\cup\{\sigma\}]$ holds when $\gamma\geq\mathbb{S}$, and
conversely $\alpha\in\mathcal{H}_{\gamma}[E_{m}(\alpha[\sigma/\mathbb{S}])]$ holds if
$\alpha\in\mathcal{H}_{\gamma}(\psi_{\Omega_{\mathbb{S}+N-m}}(\gamma))$.
This means that if we have an ordinal $\sigma\in\Theta$ and $\iota\in[\sigma]J$,
then
$\mathsf{k}(\iota)\subset\mathcal{H}_{\gamma}[\Theta\cup E_{m}(\iota^{[\sigma/\mathbb{S}]})]$
iff $\mathsf{k}(\iota^{[\sigma/\mathbb{S}]})\subset\mathcal{H}_{\gamma}[\Theta\cup E_{m}(\iota)]$
under a mild condition, cf.\,Tautology \ref{lem:tautology.cap}.\ref{lem:tautology.cap3}.
We need $E_{m}(\alpha)$ to be a set of ordinals smaller than $\Omega_{\mathbb{S}+N-m-1}$
in eliminating inferences $(\Sigma(\Omega_{\mathbb{S}+N-m})\mbox{{\rm -rfl}})$, cf.\,Collapsing \ref{lem:picollpase}.

\bdf\label{df:Em}
{\rm
For ordinal terms $\alpha\in OT_{N}$ and $0\leq m<N$,
we define recursively finite sets $E_{m}(\alpha)\subset OT_{N}$ as follows.
\benu
\item
If $\alpha$ is not strongly critical, then $E_{m}(\alpha)=\bigcup\{E_{m}(\beta):\beta\in SC(\alpha)\}$.

\item 
$E_{m}(\alpha)=\emptyset$ for $\alpha\in \{\Omega\}\cup\{\Omega_{\mathbb{S}+n} : 0\leq n\leq N\}$.

\item
$E_{m}(\Omega_{\sigma+n})=\{\sigma\}$ for $\sigma\in\Psi_{N}$ and $0\leq n\leq N$.

\item
Let $\alpha=\psi_{\pi}(a)$ with $\pi=\Omega_{\mathbb{S}+n}$. Then
let
$E_{m}(\alpha)=E_{m}(a)$.

Let $\sigma\in\Psi_{N}$ be an ordinal such that $\alpha\in E^{\mathbb{S}}_{\sigma}$. Then define
\[
E_{m}(\alpha[\sigma/\mathbb{S}])=\left\{
\begin{array}{ll}
\{\alpha[\sigma/\mathbb{S}], \sigma\}\cup E_{m}(a[\sigma/\mathbb{S}]) & \mbox{ if } n\geq N-m
\\
\{\alpha, \sigma\}
\cup E_{m}(a[\sigma/\mathbb{S}])
& \mbox{ if } n<N-m
\end{array}
\right.
\]

\eenu
}
\edf

\bprp\label{prp:Em}
Let $\rho\in\Psi_{N}\cup\{\mathtt{u}\}$ with $E^{\mathbb{S}}_{\mathtt{u}}=OT_{N}$.
\benu
\item\label{prp:Em.1}
$E_{m}(\alpha)\subset\Omega_{\mathbb{S}+N-m-1}$.

\item\label{prp:Em.2}
Let $\alpha\in E^{\mathbb{S}}_{\sigma}\cap E^{\mathbb{S}}_{\rho}$ and $\gamma\geq\mathbb{S}$.
Then $\alpha[\sigma/\mathbb{S}]\in\mathcal{H}_{\gamma}[(E_{m}(\alpha)\cup\{\sigma\})\cap E^{\mathbb{S}}_{\rho}]$.

\item\label{prp:Em.3}
Let $\alpha\in E^{\mathbb{S}}_{\sigma}\cap E^{\mathbb{S}}_{\rho}\cap \mathcal{H}_{\gamma}(\psi_{\Omega_{\mathbb{S}+N-m}}(\gamma))$.
Then $\alpha\in\mathcal{H}_{\gamma}[E_{m}(\alpha[\sigma/\mathbb{S}]) \cap E^{\mathbb{S}}_{\rho}]$.

\item\label{prp:Em.4}
Let $\alpha\in E^{\mathbb{S}}_{\rho}\cap \mathcal{H}_{\gamma}(\psi_{\Omega_{\mathbb{S}+N-m}}(\gamma))$ and $\gamma\geq\mathbb{S}$.
Then $\alpha\in\mathcal{H}_{\gamma}[E_{m}(\alpha)\cap E^{\mathbb{S}}_{\rho}]$.

\item\label{prp:Em.5}
Let $\alpha\in E^{\mathbb{S}}_{\sigma}\cap  \mathcal{H}_{\gamma}(\psi_{\Omega_{\mathbb{S}+N-m}}(\gamma))$ and $\gamma\geq\mathbb{S}$.
Then $\mathcal{H}_{\gamma}[(E_{m}(\alpha)\cup\{\sigma\})\cap E^{\mathbb{S}}_{\rho}]=\mathcal{H}_{\gamma}[(E_{m}(\alpha[\sigma/\mathbb{S}])\cup\{\sigma\})\cap E^{\mathbb{S}}_{\rho}]$.

\item\label{prp:Em.6}
If $\alpha\in E^{\mathbb{S}}_{\rho}$, then $E_{m}(\alpha)\subset E^{\mathbb{S}}_{\rho}$.

\item\label{prp:Em.7}
Let $m\leq N$, $\alpha\in  \mathcal{H}_{\gamma}(\psi_{\Omega_{\mathbb{S}+N-m-1}}(\gamma))$ and $\gamma\geq\mathbb{S}$.
Then $\mathcal{H}_{\gamma}[E_{m+1}(\alpha)\cap E^{\mathbb{S}}_{\rho}]=\mathcal{H}_{\gamma}[E_{m}(\alpha)\cap E^{\mathbb{S}}_{\rho}]$.

\eenu
\eprp
\bprf
Each is seen by induction on the lengths $\ell\alpha$ of ordinal terms $\alpha$.
Let $\alpha=\psi_{\pi}(a)$ with $\pi=\Omega_{\mathbb{S}+n}$.
\\
\ref{prp:Em}.\ref{prp:Em.2}.
IH with $a\in E^{\mathbb{S}}_{\sigma}\cap E^{\mathbb{S}}_{\rho}$ yields 
$\gamma\geq\mathbb{S}>a[\sigma/\mathbb{S}]\in\mathcal{H}_{\gamma}[(E_{m}(a)\cup\{\sigma\})\cap E^{\mathbb{S}}_{\rho}]$,
where
$E_{m}(a)=E_{m}(\alpha)$. 
We obtain
$\alpha[\sigma/\mathbb{S}]=\psi_{\Omega_{\sigma+n}}(a[\sigma/\mathbb{S}])\in\mathcal{H}_{\gamma}[(E_{m}(\alpha)\cup\{\sigma\})\cap E^{\mathbb{S}}_{\rho}]$.
\\
\ref{prp:Em}.\ref{prp:Em.3}.
We may assume $n\geq N-m$.
Then $a<\gamma$ by $\alpha\in\mathcal{H}_{\gamma}(\psi_{\Omega_{\mathbb{S}+N-m}}(\gamma))$.
On the other, IH yields $a\in\mathcal{H}_{\gamma}[E_{m}(a[\sigma/\mathbb{S}]) \cap E^{\mathbb{S}}_{\rho}]$ with $E_{m}(a[\sigma/\mathbb{S}])\subset E_{m}(\alpha[\sigma/\mathbb{S}])$.
We obtain $\alpha\in\mathcal{H}_{\gamma}[E_{m}(\alpha[\sigma/\mathbb{S}]) \cap E^{\mathbb{S}}_{\rho}]$.
\\
\ref{prp:Em}.\ref{prp:Em.5} follows from Propositions \ref{prp:Em}.\ref{prp:Em.2}, \ref{prp:Em}.\ref{prp:Em.3} and \ref{prp:Em}.\ref{prp:Em.4}.
\\
\ref{prp:Em}.\ref{prp:Em.6}.
Let $\alpha=\psi_{\pi}(a)\in  E^{\mathbb{S}}_{\sigma}$ with $\pi=\Omega_{\mathbb{S}+n}$ and $\sigma<\rho$.
We obtain $a\in E^{\mathbb{S}}_{\sigma}\subset E^{\mathbb{S}}_{\rho}$, and $E_{m}(a[\sigma/\mathbb{S}])\subset E^{\mathbb{S}}_{\rho}$ by IH.
Furthermore we have $\{\alpha[\sigma/\mathbb{S}],\alpha,\sigma\}\subset E^{\mathbb{S}}_{\rho}$.
\\
\ref{prp:Em}.\ref{prp:Em.7}.
We obtain $\alpha\in  \mathcal{H}_{\gamma}(\psi_{\Omega_{\mathbb{S}+N-m}}(\gamma))$.
Let $\alpha=\psi_{\pi}(a)\in  E^{\mathbb{S}}_{\sigma}$ with $\pi=\Omega_{\mathbb{S}+N-m-1}$.
Then $E_{m}(\alpha[\sigma/\mathbb{S}])=\{\alpha,\sigma\}\cup E_{m}(a[\sigma/\mathbb{S}])$, and
$E_{m+1}(\alpha[\sigma/\mathbb{S}])=\{\alpha[\sigma/\mathbb{S}],\sigma\}\cup E_{m+1}(a[\sigma/\mathbb{S}])$.
We obtain
$\mathcal{H}_{\gamma}[E_{m+1}(a[\sigma/\mathbb{S}])\cap E^{\mathbb{S}}_{\rho}]=\mathcal{H}_{\gamma}[E_{m}(a[\sigma/\mathbb{S}])\cap E^{\mathbb{S}}_{\rho}]$
by IH.
Proposition \ref{prp:Em}.\ref{prp:Em.3} yields 
$\{\alpha\}\cap E^{\mathbb{S}}_{\rho}\subset \mathcal{H}_{\gamma}[E_{m+1}(\alpha[\sigma/\mathbb{S}])\cap E^{\mathbb{S}}_{\rho}]$ by
$\alpha\in  \mathcal{H}_{\gamma}(\psi_{\Omega_{\mathbb{S}+N-m-1}}(\gamma))$, while
$\{\alpha[\sigma/\mathbb{S}]\}\cap E^{\mathbb{S}}_{\rho}\subset \mathcal{H}_{\gamma}[E_{m}(\alpha[\sigma/\mathbb{S}])\cap E^{\mathbb{S}}_{\rho}]$ 
by Proposition \ref{prp:Em}.\ref{prp:Em.4} and
$\gamma\geq\mathbb{S}$.
\eprf

\bdf\label{df:Emk}
{\rm
$E_{m}(\iota):=\bigcup\{E_{m}(\alpha): \alpha\in\mathsf{k}(\iota)\}$
for $RS$-terms and $RS$-formulas $\iota$.
}
\edf

\bprp\label{prp:Emk}
Let $\iota$ be $RS$-terms and $RS$-formulas.
\benu
\item\label{prp:Emk.1}
$E_{m}(\iota)\subset\mathbb{S}$.

\item\label{prp:Emk.234}
Let $\mathsf{k}(\iota)\subset E^{\mathbb{S}}_{\sigma}\cap \mathcal{H}_{\gamma}(\psi_{\Omega_{\mathbb{S}+N-m}}(\gamma))$
and $\gamma\geq\mathbb{S}$.
Then $\mathcal{H}_{\gamma}[(\mathsf{k}(\iota)\cup E_{m}(\iota)\cup \{\sigma\})\cap E^{\mathbb{S}}_{\rho}]=
\mathcal{H}_{\gamma}[(\mathsf{k}(\iota^{[\sigma/\mathbb{S}]})\cup E_{m}(\iota{[\sigma/\mathbb{S}]})\cup\{\sigma\})\cap E^{\mathbb{S}}_{\rho}]$.

\item\label{prp:Emk.6}
If $\mathsf{k}(\iota)\subset E^{\mathbb{S}}_{\rho}$, then $E_{m}(\iota)\subset E^{\mathbb{S}}_{\rho}$.

\eenu
\eprp
\bprf
\ref{prp:Emk}.\ref{prp:Emk.234} follows from Propositions \ref{prp:Em}.\ref{prp:Em.2}, \ref{prp:Em}.\ref{prp:Em.3} and \ref{prp:Em}.\ref{prp:Em.5}.
\\
\ref{prp:Emk}.\ref{prp:Emk.6} follows from Proposition \ref{prp:Em}.\ref{prp:Em.6}.
\eprf

\subsection{Operator controlled $*$-derivations}\label{subsec:*derivation}
Inference rules are formulated in one-sided sequent calculi.
Let
$\mathcal{H}_{\gamma}[\Theta]:=\mathcal{H}_{\gamma}(\Theta)$, and
$\mathcal{H}_{\gamma}:=\mathcal{H}_{\gamma}(\emptyset)$.
We define a derivability relation $(\mathcal{H}_{\gamma}, \Theta)
\vdash^{* a}_{c}\Gamma$, 
where
$c$ is a bound of ranks of cut formulas and of initial sequents $({\rm stbl})$.
The derivability relation is designed to do the following job.
An infinitary derivation $\mathcal{D}_{0}$ in the relation $\vdash^{*}$ arises from a finitary proof,
in which initial sequents $({\rm stbl})$ occur to prove an axiom for stability, 
cf.\,Lemma \ref{th:embedreg}.

Inferences $(\Sigma\mbox{{\rm -rfl}})$ of $\mathbb{K}=\Omega_{\mathbb{S}+N}$ are
removed by collapsing,
ranks of formulas in the derivations are lowered to ordinals less than $\mathbb{K}$, and
the initial sequents $({\rm stbl})$ are replaced by inferences $({\rm rfl}(\rho,e,f))$
by putting a cap $\rho$ on formulas to get a derivation $\mathcal{D}_{1}$
 in another
derivability relation 
$(\mathcal{H}_{\gamma},\Theta, \mathtt{ Q}
)\vdash^{a}_{c,d,0,\Lambda_{0},\gamma_{0}}\Gamma$, cf.\,Capping \ref{lem:capping} in 
subsection \ref{subsec:derivationcap}.

\bdf\label{df:controlder*}
{\rm
$(\mathcal{H}_{\gamma},\Theta)\vdash^{* a}_{c} \Gamma$ holds
for a set
$\Gamma$ of formulas
if 
\begin{equation}\label{eq:controlder*1}
\{\gamma,a,c\}\cup E_{0}(\{\gamma,a,c\}) \cup \mathsf{k}(\Gamma)
\cup E_{0}(\Gamma)
\subset\mathcal{H}_{\gamma}[
\Theta
]
\end{equation}

and one of the following cases holds
:

\begin{description}

\item[$(\bigvee)$]\footnote{The condition $|\iota|< a$ is absent in the inference $(\bigvee)$.}
There exist 
$A\simeq\bigvee\{A_{\iota}: \iota\in J\}$,
$a(\iota)<a$ and an $\iota\in J$ such that $A\in\Gamma$ and
$(\mathcal{H}_{\gamma},\Theta)\vdash^{* a(\iota)}_{c}\Gamma,
A_{\iota}$.

\item[$(\bigwedge)$]
There exist 
$A\simeq\bigwedge\{A_{\iota}: \iota\in J\}$, 
ordinals $a(\iota)<a$ 
for each $\iota\in J$ such that
$A\in\Gamma$ and
$(\mathcal{H}_{\gamma},\Theta\cup \mathsf{k}(\iota)
\cup E_{0}(\iota)
)
\vdash^{* a(\iota)}_{c}\Gamma,
A_{\iota}$.

\item[$(cut)$]
There exist $a_{0}<a$ and $C$
such that 
$\mathrm{rk}(C)<c$,
$(\mathcal{H}_{\gamma},\Theta)\vdash^{* a_{0}}_{c}\Gamma,\lnot C$
and
$(\mathcal{H}_{\gamma},\Theta)\vdash^{* a_{0}}_{c}C, \Gamma$.

\item[$(\Sigma(\pi)\mbox{{\rm -rfl}})$]
There exist 
$a_{\ell}, a_{r}<a$ and a formula $C\equiv(\forall x\in u\, B(x))\in\Sigma(\pi)$ for a 
$\pi\in\{\Omega\}\cup\{\Omega_{\mathbb{S}+n+1} : 0\leq n<N\}$, a $u\in Tm(\pi)$ and a $B(\mathsf{L}_{0})\in\Sigma_{1}(\pi)$
such that 
$c>\pi$,
$(\mathcal{H}_{\gamma},\Theta)
\vdash^{* a_{\ell}}_{c}\Gamma,C$
and
$(\mathcal{H}_{\gamma},\Theta)
\vdash^{* a_{r}}_{c}
\lnot \exists x<\pi\,C^{(x,\pi)}, \Gamma$.

\item[({\rm stbl})]
There exist a $\bigwedge$-formula 
$B(\mathsf{L}_{0})\in\Delta_{0}(\mathbb{S})$
and a term $u\in Tm(\mathbb{K})$ 
such that
$\mathrm{rk}(B(u))<c$ and
$\{\lnot B(u),\exists x\in\mathsf{L}_{\mathbb{S}}B(x)\}\subset\Gamma$.

\end{description}
}
\edf

\blem\label{lem:tautology*}{\rm (Tautology)}
$(\mathcal{H}_{0},\Theta\cup\mathsf{k}(A)\cup E_{0}(A))
\vdash^{* 2d}_{0}
\lnot A, A$
 holds for
$d=\mathrm{rk}(A)$.
\elem
\bprf
By induction on $d$.
\eprf

\blem\label{lem:equality*}
{\rm (Equality)}
$(\mathcal{H}_{0},\Theta\cup\mathsf{k}(A,u,v)\cup E_{0}(A,u,v) )
\vdash^{* \omega(|u|\#|v|)\# 2d}_{0}
u\neq v,\lnot A(u), A(v)$ holds for $d=\mathrm{rk}(A(\mathsf{L}_{0}))$.

\elem
\bprf
By induction on $d$, cf.\cite{Buchholz,OA}.
\eprf

\blem\label{lem:inversionreg*}{\rm (Inversion)}
Let  $A\simeq \bigwedge(A_{\iota})_{\iota\in J}$,
$(\mathcal{H}_{\gamma},\Theta)\vdash^{* a}_{c}\Gamma,A$ and $\mathrm{rk}(A)\geq\mathbb{K}$.
Then for each $\iota\in J$, 
$(\mathcal{H}_{\gamma},\Theta\cup\mathsf{k}(\iota)\cup E_{0}(\iota)
)
\vdash^{* a}_{c}
\Gamma,A_{\iota}$ holds.
\elem
\bprf
By induction on $a$.
Since $\mathrm{rk}(A)\geq\mathbb{K}$, $A$ is not a major formula of any of  (stbl).
\eprf

\blem\label{lem:reduction*}{\rm (Reduction)}
Let $(\mathcal{H}_{\gamma},\Theta)\vdash^{*a}_{c}\Gamma_{0},\lnot C$ and
$(\mathcal{H}_{\gamma},\Theta)\vdash^{*b}_{c}C,\Gamma_{1}$ for
$C\simeq\bigvee(C_{\iota})_{\iota\in J}$ with $\mathbb{K}\leq\mathrm{rk}(C)\leq c$ and $a\geq b$.
Then
$(\mathcal{H}_{\gamma},\Theta)\vdash^{* a+b}_{c}\Gamma_{0},\Gamma_{1}$ holds.
\elem
\bprf 
By induction on $b$.
Since $\mathrm{rk}(C)\geq\mathbb{K}$, $C$ is not a major formula of any  (stbl).
If $b=0$, then $(\mathcal{H}_{\gamma},\Theta)\vdash^{*b}_{c}C,\Gamma_{1}$ follows from a void $(\bigwedge)$
with a major formula in $\Gamma_{1}$.
Let $b>0$.
If $\mathrm{rk}(C)<c$, then a $(cut)$ yields the lemma by $b\leq a<a+b$.
Consider the case when the last inference in $(\mathcal{H}_{\gamma},\Theta)\vdash^{*b}_{c}C,\Gamma_{1}$ is a $(\bigvee)$ with the
major formula $C$.
We have $(\mathcal{H}_{\gamma},\Theta)\vdash^{*b_{0}}_{c}C_{\iota},C,\Gamma_{1}$ for $\iota\in J$ and $b_{0}<b$.
IH yields $(\mathcal{H}_{\gamma},\Theta)\vdash^{* a+b_{0}}_{c}C_{\iota},\Gamma_{0},\Gamma_{1}$.
On the other hand we have
$(\mathcal{H}_{\gamma},\Theta\cup\mathsf{k}(\iota)\cup E_{0}(\iota))\vdash^{*a}_{c}\Gamma_{0},\lnot C_{\iota}$ by Inversion \ref{lem:inversionreg*}.
Assuming $\mathsf{k}(\iota)\subset\mathsf{k}(C_{\iota})$, we obtain $\mathsf{k}(\iota)\cup E_{0}(\iota)\subset\mathcal{H}_{\gamma}[\Theta]$ by (\ref{eq:controlder*1}), and
$(\mathcal{H}_{\gamma},\Theta)\vdash^{*a}_{c}\Gamma_{0},\lnot C_{\iota}$.
We obtain $\mathrm{rk}(C_{\iota})<\mathrm{rk}(C)$ by Proposition  \ref{lem:rank}.\ref{lem:rank2}, and
$(\mathcal{H}_{\gamma},\Theta)\vdash^{*a+b}_{c}\Gamma_{0},\Gamma_{1}$ by a $(cut)$.
\eprf

\blem\label{lem:predcereg*}{\rm (Cut-elimination)}\\
Suppose
$(\mathcal{H}_{\gamma},\Theta)\vdash^{* a}_{\mathbb{K}+1+m}\Gamma$ for $m<\omega$.
Then $(\mathcal{H}_{\gamma},\Theta)\vdash^{\omega_{m}(a)}_{\mathbb{K}+1}\Gamma$ holds.
\elem
\bprf
By main induction on $m$ with subsidiary induction on $a$.
\eprf

\blem\label{th:embedreg}{\rm (Embedding of Axioms)}
For each axiom $A$ in $T_{N}$, there is an $m<\omega$ such that
 $(\mathcal{H}_{0},\emptyset)\vdash^{* \mathbb{K}\cdot 2}_{\mathbb{K}+m} A$
holds.
\elem
\bprf
We show that the axiom $A(v) \land v\in \mathsf{L}_{\mathbb{S}} \to A^{(\mathbb{S},\mathbb{K})}(v)\,(A\in\Sigma_{1})$ follows 
by an inference $({\rm stbl})$.
Let $B(\mathsf{L}_{0})\in\Delta_{0}(\mathbb{S})$ be a $\bigwedge$-formula 
and $u\in Tm(\mathbb{K})$.
Then $\mathrm{rk}(B(u))<\mathbb{K}$.
Let $\mathsf{k}_{0}=\mathsf{k}(B(\mathsf{L}_{0}))$ and $\mathsf{k}_{u}=\mathsf{k}(u)$.
We obtain
$(\mathcal{H}_{0},\mathsf{k}_{0}\cup\mathsf{k}_{u})
\vdash^{* 0}_{\mathbb{K}}\lnot B(u),\exists x\in \mathsf{L}_{\mathbb{S}}B(x)$ by a $({\rm stbl})$.
A $(\bigwedge)$ yields
$
(\mathcal{H}_{0},\mathsf{k}_{0})
\vdash^{*1}_{\mathbb{K}}\lnot \exists x\, B(x),\exists x\in \mathsf{L}_{\mathbb{S}}B(x)
$.
\eprf

\blem\label{th:embedregthm}{\rm (Embedding)}
If $T_{N}\vdash \Gamma$ for sets $\Gamma$ of sentences, 
there are $m,k<\omega$ such that 
 $(\mathcal{H}_{0},\emptyset) \vdash_{\mathbb{K}+m}^{*\mathbb{K}\cdot 2+k}\Gamma$ 
holds.
\elem

\blem\label{lem:Kcollpase*}{\rm (Collapsing)}
Assume
$\Theta\subset
\mathcal{H}_{\gamma}(\psi_{\mathbb{K}}(\gamma))$
 and
$
(\mathcal{H}_{\gamma},\Theta
)\vdash^{*a}_{\mathbb{K}+1}\Gamma
$ with $\Gamma\subset\Sigma(\mathbb{K})$.
Then 
$(\mathcal{H}_{\hat{a}+1},\Theta)
\vdash^{* \beta}_{\beta}
\Gamma^{(\beta,\mathbb{K})}$ holds for $\hat{a}=\gamma+\omega^{a}$ and $\Lambda_{0}=\beta=\psi_{\mathbb{K}}(\hat{a})$.
\elem
\bprf
This is seen as in \cite{Buchholz}
by induction on $a$.
We have
$\{\gamma,a\}\cup E_{0}(\{\gamma,a\})\cup\mathsf{k}(\Gamma)\cup E_{0}(\Gamma)\subset\mathcal{H}_{\gamma}[\Theta]$ 
by (\ref{eq:controlder*1}), and $E_{0}(\beta)=E_{0}(\{\gamma,a\})$.
We obtain
$\{\beta\}\cup E_{0}(\beta)\subset\mathcal{H}_{\hat{a}+1}[\Theta]$ for (\ref{eq:controlder*1}).
\\
\textbf{Case 1}.
First consider the case when the last inference is a $({\rm stbl})$: 
We have 
a $\bigwedge$-formula
$B(\mathsf{L}_{0})\in\Delta_{0}(\mathbb{S})$, and
a term $u\in Tm(\mathbb{K})$
such that
$\mathsf{k}(B(u))\subset \mathcal{H}_{\gamma}[\Theta]\cap\mathbb{K}$ by (\ref{eq:controlder*1}).
The assumption $\Theta\subset
\mathcal{H}_{\gamma}(\psi_{\mathbb{K}}(\gamma))$ yields
$\mathsf{k}(B(u))\subset
\mathcal{H}_{\gamma}(\psi_{\mathbb{K}}(\gamma))\cap\mathbb{K}=\psi_{\mathbb{K}}(\gamma)$, and
$\mathrm{rk}(B(u))<\psi_{\mathbb{K}}(\gamma)\leq\beta$.
\\
\textbf{Case 2}. 
Second consider the case when the last inference is a $(\Sigma(\pi)\mbox{{\rm -rfl}})$ on 
$\pi\leq\mathbb{K}$:
There exist ordinals 
$a_{\ell}, a_{r}<a$ and a formula $C\in\Sigma(\pi)$
such that 
$(\mathcal{H}_{\gamma},\Theta)
\vdash^{* a_{\ell}}_{\mathbb{K}+1}\Gamma,C$,
and
$(\mathcal{H}_{\gamma},\Theta
)\vdash^{* a_{r}}_{\mathbb{K}+1}
\lnot \exists x<\pi\, C^{(x,\pi)},\Gamma$.
Let $\pi=\mathbb{K}$.
IH yields
$(\mathcal{H}_{\widehat{a}+1},\Theta)
\vdash^{*\beta_{\ell}}_{\beta}
\Gamma^{(\beta,\mathbb{K})}, C^{(\beta_{\ell},\mathbb{K})}$ for
$\beta_{\ell}=\psi_{\mathbb{K}}(\widehat{a_{\ell}})$ with 
$\widehat{a_{\ell}}=\gamma+\omega^{a_{\ell}}$, where
$\beta_{\ell}\in\mathcal{H}_{\widehat{a_{\ell}}+1}[\Theta]$ and
$E_{0}(\beta_{\ell})=E_{0}(\widehat{a_{\ell}})\subset E_{0}(\{\gamma,a_{\ell}\})\subset\mathcal{H}_{\gamma}[\Theta]$.

On the other hand we have
$(\mathcal{H}_{\widehat{a_{\ell}}+1},\Theta\cup\{\beta_{\ell}\}\cup E_{0}(\beta_{\ell})
)\vdash^{* a_{r}}_{\mathbb{K}+1}
\lnot C^{(\beta_{\ell},\mathbb{K})},\Gamma$
by Inversion \ref{lem:inversionreg*}, and
$(\mathcal{H}_{\widehat{a_{\ell}}+1},\Theta
)\vdash^{* a_{r}}_{\mathbb{K}+1}
\lnot C^{(\beta_{\ell},\mathbb{K})},\Gamma$.

Let
$\beta_{r}=\psi_{\mathbb{K}}(\widehat{a_{r}})<\beta$ with 
$\widehat{a_{r}}=\widehat{a_{\ell}}+1+\omega^{a_{r}}=\gamma+\omega^{a_{\ell}}+\omega^{a_{r}}<\gamma+\omega^{a}=\hat{a}$.
IH yields
$(\mathcal{H}_{\hat{a}+1},\Theta)
\vdash^{* \beta_{r}}_{\beta}
\lnot C^{(\beta_{\ell},\mathbb{K})},\Gamma^{(\beta,\mathbb{K})}$.
A $(cut)$ yields
$(\mathcal{H}_{\hat{a}+1},\Theta)
\vdash^{* \beta}_{\beta}
\Gamma^{(\beta,\mathbb{K})}$ for
$\mathrm{rk}(C^{(\beta_{\ell},\mathbb{K})})<\beta$.
When $\pi<\mathbb{K}$, IH followed by a $(\Sigma(\pi)\mbox{{\rm -rfl}})$
yields the lemma.
\\
\textbf{Case 3}.
Third consider the case when the last inference is a 
$(cut)$:
There exist $a_{0}<a$ and a $\bigvee$-formula $C$
such that $\mathrm{rk}(C)\leq\mathbb{K}$,
$(\mathcal{H}_{\gamma},\Theta)
\vdash^{* a_{0}}_{\mathbb{K}+1}
\Gamma,\lnot C$
and
$(\mathcal{H}_{\gamma},\Theta)
\vdash^{* a_{0}}_{\mathbb{K}+1}C,\Gamma$.
We obtain $C\in\Sigma_{1}(\mathbb{K})\cup \Delta_{0}(\mathbb{K})$ by Proposition \ref{lem:rank}.\ref{lem:rank3}.
IH yields
$(\mathcal{H}_{\widehat{a_{0}}+1},\Theta)
\vdash^{* \beta_{0}}_{\beta}C^{(\beta_{0},\mathbb{K})},\Gamma^{(\beta,\mathbb{K})}$
for $\beta_{0}=\psi_{\mathbb{K}}(\widehat{a_{0}})\in\mathcal{H}_{\widehat{a_{0}}+1}[\Theta]$ with $\widehat{a_{0}}=\gamma+\omega^{a_{0}}$,
and $E_{0}(\beta_{0})\subset\mathcal{H}_{\gamma}[\Theta]$.
On the other, we obtain
$(\mathcal{H}_{\widehat{a_{0}}+1},\Theta\cup\{\beta_{0}\}\cup E_{0}(\beta_{0})
)\vdash^{* a_{0}}_{\mathbb{K}+1}
\lnot C^{(\beta_{0},\mathbb{K})},\Gamma$
by Inversion \ref{lem:inversionreg*}, and
$(\mathcal{H}_{\widehat{a_{0}}+1},\Theta
)\vdash^{* a_{0}}_{\mathbb{K}+1}
\lnot C^{(\beta_{0},\mathbb{K})},\Gamma$.
Let $\beta_{1}=\psi_{\mathbb{K}}(\widehat{a_{1}})<\beta$ with $\widehat{a_{1}}=\widehat{a_{0}}+\omega^{a_{0}}=\gamma+\omega^{a_{0}}\cdot 2<\hat{a}$.
IH yields
$(\mathcal{H}_{\hat{a}+1},\Theta
)\vdash^{* \beta_{1}}_{\beta}
\lnot C^{(\beta_{0},\mathbb{K})},\Gamma^{(\beta,\mathbb{K})}$.
We obtain
$(\mathcal{H}_{\hat{a}+1},\Theta
)\vdash^{* \beta}_{\beta}
\Gamma^{(\beta,\mathbb{K})}$
by a $(cut)$ with $\mathrm{rk}(C^{(\beta_{0},\mathbb{K})})<\beta$.
\\
\textbf{Case 4}.  
Fourth consider the case when the last inference is a $(\bigvee)$:
A $\bigvee$-formula with
$A\in\Gamma$ is introduced.
Let $A_{i}\equiv A\simeq\bigvee\left(A_{\iota}\right)_{\iota\in J}$.
There are an $\iota\in J$, an ordinal
 $a(\iota)<a$
such that
$(\mathcal{H}_{\gamma},\Theta)\vdash^{* a(\iota)}_{\mathbb{K}+1}
\Gamma,A_{\iota}$.
We may assume $\mathsf{k}(\iota)\subset\mathsf{k}(A_{\iota})$.
We obtain 
by (\ref{eq:controlder*1}),
$\mathsf{k}(\iota)\subset\mathcal{H}_{\gamma}[\Theta]\cap\mathbb{K}
\subset \mathcal{H}_{\gamma}(\psi_{\mathbb{K}}(\gamma))\cap\mathbb{K}\subset
\psi_{\mathbb{K}}(\gamma)\subset\beta$.
IH yields $(\mathcal{H}_{\hat{a}+1},\Theta)
\vdash^{*\beta(\iota)}_{\beta}
\Gamma^{(\beta,\mathbb{K})},
A_{\iota}^{(\beta,\mathbb{K})}$
for $\beta(\iota)=\psi_{\mathbb{K}}(\widehat{a(\iota)})$ with
$\widehat{a(\iota)}=\gamma+\omega^{a(\iota)}$.
$(\mathcal{H}_{\hat{a}+1},\Theta
\vdash^{*\beta}_{\beta}
\Gamma^{(\beta,\mathbb{K})}$ 
follows from a $(\bigvee)$.
\\
\textbf{Case 5}.
Fifth consider the case when the last inference is a $(\bigwedge)$:
A formula $A\in\Sigma(\mathbb{K})$ with $A\simeq\bigwedge\left(A_{\iota}\right)_{\iota\in J}$ 
is introduced in $\Gamma$.
For every
$\iota\in J$ there exists an $a(\iota)<a$ 
such that
$(\mathcal{H}_{\gamma},\Theta\cup\mathsf{k}(\iota)\cup E_{0}(\iota)
)
\vdash^{* a(\iota)}_{\mathbb{K}+1}\Gamma, A_{\iota}$.
We see $\mathsf{k}(\iota)\cup E_{0}(\iota)\subset\psi_{\mathbb{K}}(\gamma)$ as follows.
First $E_{0}(\iota)\subset\Omega_{\mathbb{S}+N-1}<\psi_{\mathbb{K}}(\gamma)$ by Proposition \ref{prp:Emk}.\ref{prp:Emk.1}.
Second for example let $A\equiv(\forall x\in u\, B(x))$. Then 
$u\in Tm(\mathbb{K})$ by 
$A\in\Sigma(\mathbb{K})$, and $J=Tm(|u|)$. We obtain 
$\mathsf{k}(u)\subset\mathcal{H}_{\gamma}[\Theta]\cap\mathbb{K}\subset\psi_{\mathbb{K}}(\gamma)$.
Hence $|\iota|<|u|<\psi_{\mathbb{K}}(\gamma)$ for $\iota\in J$.
IH yields 
$(\mathcal{H}_{\hat{a}+1},\Theta\cup\mathsf{k}(\iota)\cup E_{0}(\iota)
)
\vdash^{* \beta(\iota)}_{\beta}
\Gamma,
A_{\iota}$
for each $\iota\in J$, where
$\beta(\iota)=\psi_{\mathbb{K}}(\widehat{a(\iota)})$ with
$\widehat{a(\iota)}=\gamma+\omega^{a(\iota)}$.
We obtain
$(\mathcal{H}_{\hat{a}+1},\Theta)\vdash^{* \beta}_{\beta}
\Gamma$ by a $(\bigwedge)$ and $\beta(\iota)<\beta$.
\eprf

\subsection{Stepping-down}\label{subsec:stepping}
The following Definition \ref{df:hstepdown} is mainly needed in subsection \ref{subsec:elimrfl}, but
also in Definition \ref{df:resolvent}, cf.\,the beginning of subsection \ref{subsec:elimrfl}.
Let $\Lambda_{0}<\mathbb{K}$ be the ordinal in Collapsing \ref{lem:Kcollpase*}, and $\Lambda=\Gamma(\Lambda_{0})$.

\bdf\label{df:special}
{\rm
Let
$s(f)=\max(\{0\}\cup{\rm supp}(f))$ for finite function $f$, and
$s(\rho)=s(m(\rho))$.
Let $f$ be a non-empty and irreducible finite function.
Then $f$ is said to be \textit{special} if there exists an ordinal $\alpha$
such that $f(s(f))=\alpha+\Lambda$
for the base $\Lambda$ of the $\tilde{\theta}$-function.
For a special finite function $f$, $f^{\prime}$ denotes a finite function such that
${\rm supp}(f^{\prime})={\rm supp}(f)$,
$f^{\prime}(c)=f(c)$ for $c\neq s(f)$, and
$f^{\prime}(s(f))=\alpha$ with $f(s(f))=\alpha+\Lambda$.
}
\edf
A special function $g$ has a room $\Lambda$ on its top $g(s(g))=\alpha+\Lambda$.
A stepping-down $h^{b}(g;a)$ of a special function $g$ to an ordinal $b$
is introduced in Definition \ref{df:hstepdown} by replacing the room $\Lambda$
by an ordinal $a$.
Such a stepping-down is needed to analyze inference rules for reflections in subsection
\ref{subsec:elimrfl}.

\bdf\label{df:hstepdown}
{\rm
Let $f,g:\Lambda\to \Gamma(\Lambda)$ be special finite functions and 
$\tilde{\theta}_{b}$ denotes the $b$-th iterate of the function $\tilde{\theta}_{1}(\xi)=\Lambda^{\xi}$ in
Definition \ref{df:Lam}.

\begin{enumerate}
\item\label{df:hstepdown.1}
For 
ordinals $a<\Lambda$, $b\leq s(g)$, 
let us define a finite function $h=h^{b}(g;a):\Lambda\to\Gamma(\Lambda)$ 
as follows.
$s(h)=b$, and
$h_{b}=g_{b}$.
To define $h(b)$,
let $\{b=b_{0}<b_{1}<\cdots<b_{n}=s(g)\}=\{b,s(g)\}\cup\left((b,s(g))\cap {\rm supp}(g)\right)$.
Define recursively ordinals $\alpha_{i}$ by
$\alpha_{n}=\alpha+a$ with $g(s(g))=\alpha+\Lambda$.
$\alpha_{i}=g(b_{i})+\tilde{\theta}_{c_{i}}(\alpha_{i+1})$ for
$c_{i}=b_{i+1}-b_{i}$.
Finally let $h(b)=\alpha_{0}+\Lambda$.

\item\label{df:hstepdown.4}
Let $b=s(f)$. Then
$f*g^{b+1}$ denotes a special finite function $h$ defined as follows.
If $g^{b+1}=\emptyset$, then $h=f$.
Let $g^{b+1}\neq\emptyset$. Then
${\rm supp}(h)={\rm supp}(f)\cup{\rm supp}(g^{b+1})$, 
$h(c)=f(c)$ for $c<b$, $h(b)=f^{\prime}(b)$ and
$h(c)=g(c)$ for $c>b$.

\end{enumerate}
}
\edf

\bprp\label{prp:hstepdown}
Let $a\leq\Lambda$, and
$f, g:\Lambda\to\Gamma(\Lambda)$ be special finite functions with a strongly critical number $\Lambda<\mathbb{K}$ such that
$f_{d}=g_{d}$ and
$f<^{d}g^{\prime}(d)$ for a $d\in{\rm supp}(g)$.
Let $\rho\in\Psi_{N}$ with $g=m(\rho)$.

\begin{enumerate}

\item\label{prp:hstepdown.1}
Let $d<s(f)$ and $h=h^{d}(f;a)$.
Then $h_{d}=g_{d}$ and $h<^{d}g^{\prime}(d)$.

\item\label{prp:hstepdown.2}
If $b<d$, then
$f_{b}=(h^{b}(g;a))_{b}$,
$f<^{b}(h^{b}(g;a))^{\prime}(b)$.

\item\label{prp:hstepdown.4}
Let $b<b_{0}<d$, 
$a_{0}, a_{1}<a< \Lambda$, $g_{0}=h^{b_{0}}(g;a_{0})*f^{b_{0}+1}$, and $k=h^{b}(g_{0};a_{1})$
Then
$k_{b}=(h^{b}(g;a))_{b}$ and
$k<^{b}(h^{b}(g;a))^{\prime}(b)$.

\end{enumerate}
\eprp
\bprf
\ref{prp:hstepdown}.\ref{prp:hstepdown.1}.
We have $h_{d}=f_{d}=g_{d}$.
We show $h(d)<g^{\prime}(d)$. 
Let $\{d=d_{0}<d_{1}<\cdots<d_{n}=s(f)\}=\{d,s(f)\}\cup\left((d,s(f))\cap {\rm supp}(f)\right)$ for $n>0$.
Let $\alpha_{i}$ be ordinals defined by
$\alpha_{n}=\alpha+a$ with $f(s(f))=\alpha+\Lambda$.
$\alpha_{i}=f(d_{i})+\tilde{\theta}_{c_{i}}(\alpha_{i+1})$ for
$c_{i}=d_{i+1}-d_{i}$.
Then $h(d)=\alpha_{0}+\Lambda$.

On the other, let $\mu_{0}$ be a part of $g^{\prime}(d_{0})$ such that $f(d_{0})<\mu_{0}$ and
$f<^{d_{1}}\tilde{\theta}_{-c_{0}}(tl(\mu_{0}))$.
For $i<n$, let $\mu_{i+1}$ be a part of $\tilde{\theta}_{-c_{i}}(tl(\mu_{i}))$ such that $f(d_{i+1})<\mu_{i+1}$.

For $f(d_{n})=\alpha+\Lambda$ we have $\alpha_{n}=\alpha+a<\alpha+\Lambda<\mu_{n}$, and
$\tilde{\theta}_{c_{n-1}}(\alpha+a)<\tilde{\theta}_{c_{n-1}}(\mu_{n})\leq\tilde{\theta}_{c_{n-1}}(\tilde{\theta}_{-c_{n-1}}(tl(\mu_{n-1})))\leq
tl(\mu_{n-1})$ by Proposition \ref{prp:tht4}.\ref{prp:tht4.2}.
We obtain $\alpha_{n-1}=f(d_{n-1})+\tilde{\theta}_{c_{n-1}}(\alpha_{n})<\mu_{n-1}$
by $f(d_{n-1})<\mu_{n-1}$.
We see inductively that $\alpha_{i}<\mu_{i}$, and $\alpha_{0}<\mu_{0}$.
Hence we obtain $h(d)=\alpha_{0}+\Lambda<\mu_{0}\leq g^{\prime}(d_{0})$.
\\
\ref{prp:hstepdown}.\ref{prp:hstepdown.2}.
Let $h=h^{b}(g;a)$.
We have $h_{b}=g_{b}=f_{b}$.
Let $b+x\in{\rm supp}(f)\cap d\subset{\rm supp}(g)\cap d$.
Then $f(b+x)= g(b+x)<\tilde{\theta}_{-x}(h^{\prime}(b))$ and
$g^{\prime}(d)< \tilde{\theta}_{-(d-b)}(h^{\prime}(b))$.
Proposition \ref{prp:idless} yields the proposition.
\\
\ref{prp:hstepdown}.\ref{prp:hstepdown.4}.
We have $h_{b}=g_{b}=(h^{b}(g;a))_{b}$.
We obtain $h^{\prime}(b)=(h^{b}(g;a_{0}))^{\prime}(b)<(h^{b}(g;a))^{\prime}(b)$ by $a_{0}<a$.
As in Proposition \ref{prp:hstepdown}.\ref{prp:hstepdown.2} we see that
$f(b+x)= g(b+x)<\tilde{\theta}_{-x}(h^{\prime}(b))$ and
$g^{\prime}(d)< \tilde{\theta}_{-(d-b)}(h^{\prime}(b))$ for $b+x\in{\rm supp}(f)\cap d\subset{\rm supp}(g)\cap d$.
\eprf

\subsection{Operator controlled derivations with caps}\label{subsec:derivationcap}
Our cut-elimination procedure goes roughly as follows, cf.\,the beginning of subsection \ref{subsec:*derivation}.
The initial sequents $({\rm stbl})$ in a given $*$-derivation
are replaced by inferences $({\rm rfl}(\rho,e,f))$
by putting a \textit{cap} $\rho$ on formulas, cf.\,Capping \ref{lem:capping}.
Our main task is to eliminate inferences $({\rm rfl}(\rho,e,f))$ from a resulting derivation $\mathcal{D}_{1}$.
Although a capped formula $A^{(\rho)}$ in Definition \ref{df:cap}.\ref{df:cap.1}
 is intended to denote the formula $A^{[\rho/\mathbb{S}]}$, 
we need to distinguish $A^{(\rho)}$ from $A^{[\rho/\mathbb{S}]}$.
The cap $\rho$ in $A^{(\rho)}$ is a temporary one, and the formula $A$ could
put on a smaller cap $A^{(\kappa)}$.
Let $\rho<\mathbb{S}$ be an ordinal for which inferences
$({\rm rfl}(\rho,e,f))$ occur in $\mathcal{D}_{1}$.
In Recapping \ref{mlem:singlemainl.1} the caps $A^{(\rho)}$
 are lowered by substituting a smaller
ordinal $\kappa$ for $\rho$, and simultaneously the ranks 
${\rm rk}(B)$ of formulas $B^{(\kappa)}$ to be reflected
are lowered.
In this process new inferences $({\rm rfl}(\sigma,e_{1},f_{1}))$ arise with $\sigma<\kappa$,
whose ranks might not be smaller.

Iterating this process,
we arrive at a derivation $\mathcal{D}_{2}$ such that
every formula $A$ occurring in it is in $\Sigma_{1}(\mathbb{S})\cup\Pi_{1}(\mathbb{S})$.
Then caps play no r\^{o}le, i.e.,
$A^{(\rho)}$ is `equivalent' to $A$, and
 inferences $({\rm rfl}(\rho,e,f))$ are removed from $\mathcal{D}_{2}$ by replacing 
these by a series of $(cut)$'s, cf.\,Lemmas \ref{lem:predcereg_S} and \ref{lem:main}.
See the beginning of subsection \ref{subsec:elimrfl} for more on an elimination procedure.

In what follows strongly critical numbers
$\Lambda_{0}=\beta<\Lambda=\Gamma(\beta)<\mathbb{K}$ will be fixed, which are ordinals in Collapsing \ref{lem:Kcollpase*}.
Also $\gamma_{0}$ denotes a fixed ordinal.

\bdf\label{df:cap}
{\rm
\begin{enumerate}
\item\label{df:cap.1}
By a \textit{capped formula} we mean a pair $(A,\rho)$ of $RS$-sentence $A$
and an ordinal 
$\rho<\mathbb{S}$ such that
$\mathsf{k}(A)\subset E^{\mathbb{S}}_{\rho}$.
Such a pair is denoted by $A^{(\rho)}$.
Sometimes it is convenient for us to regard \textit{uncapped formulas} $A$ as capped formulas
$A^{(\mathtt{ u})}$ with its cap $\mathtt{ u}$.
A \textit{sequent} is a finite set of capped or uncapped formulas, denoted by
$\Gamma_{0}^{(\rho_{0})},\ldots,\Gamma_{n}^{(\rho_{n})}$,
where each formula in the set $\Gamma_{i}^{(\rho_{i})}$ puts on the cap 
$\rho_{i}\in\mathbb{S}\cup\{\mathtt{ u}\}$.
When we write $\Gamma^{(\rho)}$, we tacitly assume that
$\mathsf{k}(\Gamma)\subset E^{\mathbb{S}}_{\rho}$, where
$E^{\mathbb{S}}_{\mathtt{ u}}=OT_{N}$.
A capped formula $A^{(\rho)}$ is said to be a $\Sigma(\pi)$-formula if
$A\in\Sigma(\pi)$.
Let $\mathsf{k}(A^{(\rho)}):=\mathsf{k}(A)$.

\item\label{df:cap.2}
Let $\Lambda<\mathbb{K}$ be a strongly critical number.
 A non-empty finite set $\mathtt{Q}$ of ordinals is said to be a \textit{finite family} 
(for ordinals $\Lambda,\gamma_{0}$) if each $\rho\in\mathtt{Q}$ is an ordinal 
$\rho=\psi_{\sigma}^{f}(a)\in\Psi_{N}$ such that
$m(\rho):\Lambda \to\Gamma(\Lambda)$\footnote{Actually $m(\rho):\Lambda \to\varphi_{\Lambda}(\Lambda\cdot\omega)$  suffices.} is special,
$\gamma_{0}\leq b(\rho):=a<\gamma_{0}+\mathbb{S}$, and
$\mathcal{H}_{\gamma_{0}}(\rho)\cap\mathbb{S}\subset\rho$.

For  a finite family $\mathtt{Q}$
let 
\[
\rho_{\mathtt{Q}}:=\min\mathtt{Q}
.\]

\end{enumerate}
}
\edf

\bdf\label{df:resolvent}
{\rm
$H_{\rho}(f,\gamma_{0},\Theta)$
 denotes the \textit{resolvent class} defined by
$\kappa\in H_{\rho}(f,\gamma_{0},\Theta)$ iff
$\kappa\in \Psi_{N}\cap\rho$,
$\{\rho,\kappa\}$ is a finite family,
$\Theta\cap E^{\mathbb{S}}_{\rho}
\subset E^{\mathbb{S}}_{\kappa}$,
and
$f\leq m(\kappa)$,
where $f\leq g :\Leftrightarrow\forall i(f(i)\leq g(i))$ for finite functions $f,g$.

}
\edf

We define another derivability relation 
$(\mathcal{H}_{\gamma},\Theta, \mathtt{ Q}
)\vdash^{a}_{c,d,m}\Gamma$, 
where
$c$ is a bound of ranks of cut formulas, 
$d$ is a bound of ranks of minor formulas in inferences $({\rm rfl}(\rho,e,f))$ in a witnessed derivation.
The relation depends on ordinals $\Lambda_{0},\gamma_{0}$, and
should be written as $(\mathcal{H}_{\gamma},\Theta,\mathtt{ Q})\vdash^{a}_{c,d,,m,\Lambda_{0},\gamma_{0}} \Gamma$.
The ordinals $\Lambda_{0},\gamma_{0}$ will be fixed.
So let us omit these.
Note that if $\mathsf{k}(\iota)\subset E^{\mathbb{S}}_{\rho}$, then $E_{m}(\iota)\subset E^{\mathbb{S}}_{\rho}$ by
Proposition \ref{prp:Emk}.\ref{prp:Emk.6}.

\bdf\label{df:controldercollapspi11}
{\rm
Let $\Theta$ be a finite set of ordinals,
$\gamma\leq\gamma_{0}$ and 
$c,d\leq\Lambda_{0}$. Let
$\mathtt{ Q}$  be
a finite family.
Let
$\Gamma=\bigcup\{\Gamma_{\sigma}^{(\sigma)}:\sigma\in\mathtt{Q}\}\subset\Delta_{0}(\mathbb{K})$
a set of formulas
such that 
$\mathsf{k}(\Gamma_{\sigma})\subset E^{\mathbb{S}}_{\sigma}$ 
for each $\sigma\in\mathtt{ Q}$.

$(\mathcal{H}_{\gamma},\Theta, \mathtt{Q})\vdash^{a}_{c,d,m} \Gamma$ holds
if

\begin{equation}
\label{eq:controlder_cap1}
\mathtt{Q}\subset\Theta
\end{equation}

\begin{equation}
\label{eq:controlder_cap2}
\forall\sigma\in\mathtt{Q}\left(
\mathsf{k}(\Gamma_{\sigma})\cup E_{m}(\Gamma_{\sigma}) \subset\mathcal{H}_{\gamma}[
\Theta\cap E^{\mathbb{S}}_{\sigma}
]
\right)
\end{equation}

\begin{equation}
\label{eq:controlder_cap3}
\{\gamma_{0},\gamma,\Lambda_{0},a,c,d\}\cup E_{m}(\{\gamma_{0},\gamma,\Lambda_{0},a,c,d\})
\subset
\mathcal{H}_{\gamma}[\Theta\cap E^{\mathbb{S}}_{\rho_{\mathtt{Q}}}]
\end{equation}

and one of the following cases holds:

\begin{description}

\item[$(\bigvee)$]
There exist 
$A\simeq\bigvee(A_{\iota})_{\iota\in J}$, 
an ordinal
$a(\iota)<a$, $A^{(\rho)}\in\Gamma$ with a cap $\rho\in\mathtt{ Q}$,  and an 
$\iota\in J$ 
 and
$(\mathcal{H}_{\gamma},\Theta,\mathtt{ Q})\vdash^{a(\iota)}_{c,d,m}\Gamma,
\left(A_{\iota}\right)^{(\rho)}$.

\item[$(\bigwedge)$]
There exist 
$A\simeq\bigwedge(A_{\iota})_{ \iota\in J}$,  a cap $\rho\in\mathtt{Q}$, 
ordinals $a(\iota)<a$ such that 
$A^{(\rho)}\in\Gamma$ and
$(\mathcal{H}_{\gamma},\Theta\cup\mathsf{k}(\iota)\cup E_{m}(\iota),\mathtt{Q}
)
\vdash^{a(\iota)}_{c,d,m}\Gamma,
A_{\iota}^{(\rho)}$
for each $\iota\in [\rho]J$.

\item[$(cut)$]
There exist an ordinal $a_{0}<a$ and a capped formula $C^{(\rho)}$
such that $\rho=\rho_{\mathtt{Q}}$, $\mathrm{rk}(C)<c\leq\Lambda_{0}$,
$(\mathcal{H}_{\gamma},\Theta,\mathtt{ Q})\vdash^{a_{0}}_{c,d,m}\Gamma,\lnot C^{(\rho)}$
and
$(\mathcal{H}_{\gamma},\Theta,\mathtt{ Q})\vdash^{a_{0}}_{c,d,m}C^{(\rho)},\Gamma$.

\item[$(\Sigma(\pi)\mbox{{\rm -rfl}})$]
There exist ordinals
$c>\pi\in\{\Omega\}\cup\{\Omega_{\mathbb{S}+k}:0<k<N-m\}$,
$a_{\ell}, a_{r}<a$, 
and a formula 
$C\equiv(\forall x\in u\, B(x))\in\Sigma(\pi)$ for a $u\in Tm(\pi)$ and a $B(\mathsf{L}_{0})\in\Sigma_{1}(\pi)$
such that 
$(\mathcal{H}_{\gamma},\Theta,\mathtt{ Q}
)\vdash^{a_{\ell}}_{c,d,m}\Gamma,C^{(\rho)}$
and
$(\mathcal{H}_{\gamma},\Theta,\mathtt{ Q}
)\vdash^{a_{r}}_{c,d,m}
\left(\lnot \exists x<\pi\,C^{(x,\pi)}\right)^{(\rho)}, \Gamma$ for $\rho=\rho_{\mathtt{Q}}$.

\item[$({\rm rfl}(\rho,e,f))$]
For $\rho=\rho_{\mathtt{Q}}$, there exist a finite function $f:\Lambda\to\Gamma(\Lambda)$,
ordinals $\Omega_{\mathbb{S}+N-1-m}+1\leq e\in {\rm supp}(m(\rho))$, 
$a_{0}<a$,
a set $\Xi\subset\Gamma$ of formulas such that
$A\in\Sigma_{1}(\mathbb{S})$ if $A^{(\rho)}\in\Xi$, 
and a
finite set $\Delta$ of uncapped formulas 
enjoying the following conditions (r1), (r2), (r3), (r4) and (r5).
Let $\mathtt{ Q}^{\sigma}=\mathtt{ Q}\cup\{\sigma\}$ and
$\mathrm{rk}(\Delta)=\max(\{0\}\cup\{\mathrm{rk}(\delta):\delta\in\Delta\})$.

\begin{enumerate}

\item[(r1)]
$\Delta\subset\bigvee(e):\Leftrightarrow 
\forall\delta \in\Delta
\left[
(
\delta \mbox{ {\rm is a}}
\bigvee\mbox{{\rm -formula}}
)
\,\&\,
\mathrm{rk}(\delta)<e
\right] $, 
and
$\mathrm{rk}(\Delta)<d\leq\Lambda_{0}$.

\item[(r2)]
$
f_{e}=g_{e} \,\&\, f<^{e}g(e)
$ for the finite function  $g=m(\rho)$, and
$SC_{\Lambda}(f)
\subset\mathcal{H}_{\gamma}[\Theta\cap 
E^{\mathbb{S}}_{\rho}]$.

 \item[(r3)]

For each $\delta\in\Delta$,
 $(\mathcal{H}_{\gamma},\Theta, \mathtt{Q}
 )\vdash^{a_{0}}_{c,d,m}\Gamma, \lnot\delta^{(\rho)}$ holds.

\item[(r4)]
$\max\{s(\rho),s(f)\}>\Omega_{\mathbb{S}+N-1-m}+1$: Then
$
(\mathcal{H}_{\gamma},\Theta\cup \{\sigma\},\mathtt{Q}^{\sigma}
)\vdash^{a_{0}}_{c,d,m}\Delta^{(\sigma)}, \Xi$ holds for each $\sigma\in H_{\rho}(f,\gamma_{0},\Theta)$.

 \item[(r5)]
$\max\{s(\rho),s(f)\}\leq\Omega_{\mathbb{S}+N-1-m}+1$:
 $
(\mathcal{H}_{\gamma},\Theta\cup\{\sigma\},
\mathtt{Q}^{\sigma}
)\vdash^{a_{0}}_{c,d,m}\Delta^{(\sigma)}, \Xi$ holds
for every $\sigma\in H_{\rho}(f,\gamma_{0},\Theta)$ such that
 $m(\sigma)=f$.
The case (r5) is said to be \textit{degenerated}.
{\small
\[
\hspace{-15mm}
\infer[({\rm rfl}(\rho,e,f))]{
(\mathcal{H}_{\gamma},\Theta,
\mathtt{Q}
)\vdash^{a}_{c,d,m}\Gamma
}
{
\{
(\mathcal{H}_{\gamma},\Theta,
\mathtt{Q}
)\vdash^{a_{0}}_{c,d,m}\Gamma, \lnot\delta^{(\rho)}
\}_{\delta}
&
\{
(\mathcal{H}_{\gamma},\Theta\cup\{\sigma\},
\mathtt{Q}^{\sigma}
)\vdash^{a_{0}}_{c,d,m}\Delta^{(\sigma)}, \Xi
\}_{\sigma}
}
\]
}
 
\end{enumerate}
\end{description}
}
\edf

We see from (r1) that $e$ in $({\rm rfl}(\rho,e,f))$ 
as well as $d$
is a bound of ranks of
the formulas $\delta$ to be reflected.
The conditions (\ref{eq:controlder_cap2}) and (\ref{eq:controlder_cap3})
ensure us (\ref{eq:notationsystem.11}) of Definition \ref{df:notationsystem}
in Lemma \ref{mlem:singlemainl.1}.

Note that the cap $\rho$ of cut formulas $C^{(\rho)}$ in inferences $(cut)$ as well as the cap of minor formulas in $(\Sigma(\pi)\mbox{{\rm -rfl}})$, and
ordinals $\rho$ in inferences $({\rm rfl}(\rho,e,f))$ are restricted to the case $\rho=\rho_{\mathtt{Q}}$.
Also note that the side formulas $A^{(\rho)}\in\Xi\subset\Gamma$ in the right upper sequents of a $({\rm rfl}(\rho,e,f))$
has to be $\Sigma_{1}(\mathbb{S})$-formula, cf.\,Reduction \ref{lem:predcereg}.

In this subsection the ordinals $\Lambda_{0}, \gamma_{0}$ 
will be fixed, and
we write $\vdash^{a}_{c,d,m}$ for $\vdash^{a}_{c,d,m,\Lambda_{0},\gamma_{0}}$.

\blem\label{lem:weakening}{\rm (Weakening)}
Let 
$(\mathcal{H}_{\gamma},\Theta,\mathtt{Q})\vdash^{a}_{c,d,m}\Gamma$ and
$\mathsf{k}(A)\cup E_{m}(A)\subset\mathcal{H}_{\gamma}[\Theta\cap E^{\mathbb{S}}_{\tau}]$ with
$\tau\in\mathtt{Q}$.
Then
$(\mathcal{H}_{\gamma},\Theta,\mathtt{Q})\vdash^{a}_{c,d,m}\Gamma,A^{(\tau)}$ holds.
\elem
\bprf
By induction on $a$. By the assumption $\mathsf{k}(A)\cup E_{m}(A)\subset\mathcal{H}_{\gamma}[\Theta\cap E^{\mathbb{S}}_{\tau}]$,
(\ref{eq:controlder_cap2}) is enjoyed.
In inferences ${(\rm rfl}(\rho,e,f))$, add the formula $A^{(\tau)}$ only to the left upper sequents $\Gamma, \lnot\delta^{(\rho)}$.
\eprf

\blem\label{lem:tautology.cap}
{\rm (Tautology)}
Let $\{\gamma\}\cup E_{m}(\gamma)\cup \mathsf{k}(A)\cup E_{m}(A) \subset\mathcal{H}_{\gamma}[\Theta\cap E^{\mathbb{S}}_{\rho}]$, 
$\mathsf{k}(A)\subset\Lambda_{m}$, and
$\mathtt{Q}$ be a finite family such that $\mathtt{Q}\subset\Theta$ and
 $\rho=\rho_{\mathtt{Q}}$.
Let $d=\mathrm{rk}(A)$.

\benu
\item\label{lem:tautology.cap1}
$(\mathcal{H}_{\gamma},\Theta,\mathtt{Q}
)\vdash^{2d}_{0,0,m}
\lnot A^{(\rho)}, A^{(\rho)}$ holds.

\item\label{lem:tautology.cap3}
Assume $\mathsf{k}(A)\subset \mathcal{H}_{\gamma}(\psi_{\Omega_{\mathbb{S}+N-m}}(\gamma))$,
$\gamma\geq\mathbb{S}$, and $\rho>\sigma$ be an ordinal such that $\mathtt{Q}^{\sigma}=\mathtt{Q}\cup\{\sigma\}$ is a finite family,
and $\Theta\cap E^{\mathbb{S}}_{\rho}\subset E^{\mathbb{S}}_{\sigma}$.
Then
$(\mathcal{H}_{\gamma},\Theta_{\sigma},\mathtt{Q}^{\sigma}
)\vdash^{2d}_{0,0,m}
\lnot A^{(\sigma)}, (A^{[\sigma/\mathbb{S}]})^{(\rho)}$ holds for $\Theta_{\sigma}=\Theta\cup\{\sigma\}$.

\eenu
\elem
\bprf
By induction on $d$. We have 
$ \mathsf{k}(A)\cup E_{m}(A) \subset\mathcal{H}_{\gamma}[\Theta\cap E^{\mathbb{S}}_{\rho}]$ with  $\rho=\rho_{\mathtt{Q}}$.
We obtain
$\{d\}\cup E_{m}(d)\subset\mathcal{H}_{\gamma}[\Theta\cap E^{\mathbb{S}}_{\rho}]$ by Proposition \ref{lem:rank}.\ref{lem:rank1} for (\ref{eq:controlder_cap3}).
In the proof let us write $\vdash^{a}_{0}$ for $\vdash^{a}_{0,0,m}$.
\\
\ref{lem:tautology.cap}.\ref{lem:tautology.cap1}.
Let $A\simeq\bigvee(A_{\iota})_{\iota\in J}$ and $\iota\in
[\rho]J$.
Let $d_{\iota}=\mathrm{rk}(A_{\iota})$.
We obtain
$\mathsf{k}(\iota)\cup E_{m}(\iota) \subset E^{\mathbb{S}}_{\rho}$,
and $(\mathsf{k}(\iota)\cup E_{m}(\iota))\cap E^{\mathbb{S}}_{\rho}= 
\mathsf{k}(\iota)\cup E_{m}(\iota)\subset\mathcal{H}_{\gamma}[(\Theta\cup\mathsf{k}(\iota)\cup E_{m}(\iota))\cap E^{\mathbb{S}}_{\rho}]$.
IH yields
$(\mathcal{H}_{\gamma},\Theta\cup\mathsf{k}(\iota)\cup E_{m}(\iota),\mathtt{Q}
)\vdash^{2d_{\iota}}_{0}
\lnot A_{\iota}^{(\rho)}, A_{\iota}^{(\rho)}$
for $d_{\iota}<d$.
A $(\bigvee)$ followed by a $(\bigwedge)$ yields the lemma.
\\
\ref{lem:tautology.cap}.\ref{lem:tautology.cap3}.
We have $\rho_{\mathtt{Q}^{\sigma}}=\sigma$.
We see 
$\{\gamma,d\}\cup E_{m}(\gamma,d)\cup\mathsf{k}(A)\cup E_{m}(A)
\subset\mathcal{H}_{\gamma}[\Theta\cap E^{\mathbb{S}}_{\sigma}]$
from $\Theta\cap E^{\mathbb{S}}_{\rho}\subset E^{\mathbb{S}}_{\sigma}$.
Hence (\ref{eq:controlder_cap3}) 
is enjoyed.

We have $\rho>\sigma\in\Theta_{\sigma}$.
We obtain
$\mathsf{k}(A^{[\sigma/\mathbb{S}]})\cup E_{m}(A^{[\sigma/\mathbb{S}]})\subset\mathcal{H}_{\gamma}[\Theta_{\sigma}\cap E^{\mathbb{S}}_{\rho}]$ by
Proposition \ref{prp:Emk}.\ref{prp:Emk.234} and $\mathsf{k}(A)\subset \mathcal{H}_{\gamma}(\psi_{\Omega_{\mathbb{S}+N-m}}(\gamma))$.
Hence (\ref{eq:controlder_cap2}) is enjoyed.

Let $A\simeq\bigvee(A_{\iota})_{\iota\in J}$ or  $A\simeq\bigwedge(A_{\iota})_{\iota\in J}$.
By Proposition \ref{lem:assigncollaps} we obtain
 $A^{[\sigma/\mathbb{S}]}\simeq\bigvee(A^{[\sigma/\mathbb{S}]}_{\iota})_{\iota\in [\sigma]J}$ or  
 $A^{[\sigma/\mathbb{S}]}\simeq\bigwedge(A^{[\sigma/\mathbb{S}]}_{\iota})_{\iota\in [\sigma]J}$.
 Let $I=\{\iota^{[\sigma/\mathbb{S}]}: \iota\in[\sigma]J\}$.
 We see 
 $\iota\in 
 [\sigma]J \Leftrightarrow \iota^{[\sigma/\mathbb{S}]}\in[\rho]I
 =I$
 from $\sigma<\rho$.
 Note that when $\mathrm{rk}(A)<\mathbb{S}$, we have $A^{[\sigma/\mathbb{S}]}\equiv A$ and $J=[\sigma]J=I$.
 
For each $\iota\in [\sigma]J$, we have $d_{\iota}=\mathrm{rk}(A_{\iota})<d$ by Proposition \ref{lem:rank}.\ref{lem:rank2}.
On the other hand we have $\mathsf{k}(\iota)\subset E^{\mathbb{S}}_{\sigma}\cap \mathcal{H}_{\gamma}(\psi_{\Omega_{\mathbb{S}+N-m}}(\gamma))$
 by Proposition \ref{prp:Emk}.\ref{prp:Emk.6}, 
$\iota\in[\sigma]J$ and $\mathsf{k}(A)\subset \mathcal{H}_{\gamma}(\psi_{\Omega_{\mathbb{S}+N-m}}(\gamma))$.
IH yields
$(\mathcal{H}_{\gamma},\Theta_{\sigma}\cup\mathsf{k}(\iota)\cup E_{m}(\iota),\mathtt{Q}^{\sigma}
)\vdash^{2d_{\iota}}_{0}
\lnot A_{\iota}^{(\sigma)}, (A_{\iota}^{[\sigma/\mathbb{S}]})^{(\rho)}$.
On the other hand we have
 $\mathcal{H}_{\gamma}[(\Theta\cup\mathsf{k}(\iota)\cup E_{m}(\iota))\cap E^{\mathbb{S}}_{\tau}]=
 \mathcal{H}_{\gamma}[(\Theta\cup\mathsf{k}(\iota^{[\sigma/\mathbb{S}]})\cup  E_{m}(\iota^{[\sigma/\mathbb{S}]}))\cap E^{\mathbb{S}}_{\tau}]$ 
 for any $\tau$ by Proposition \ref{prp:Emk}.\ref{prp:Emk.234}.
 We obtain
 $(\mathcal{H}_{\gamma},\Theta_{\sigma}\cup\mathsf{k}(\iota^{[\sigma/\mathbb{S}]})\cup  E_{m}(\iota^{[\sigma/\mathbb{S}]}),\mathtt{Q}^{\sigma}
)\vdash^{2d_{\iota}}_{0}
\lnot A_{\iota}^{(\sigma)}, (A_{\iota}^{[\sigma/\mathbb{S}]})^{(\rho)}$.
 A $(\bigvee)$ followed by a $(\bigwedge)$ yields 
 $(\mathcal{H}_{\gamma},\Theta_{\sigma},\mathtt{Q}^{\sigma}
)\vdash^{2d}_{0}
\lnot A^{(\sigma)}, (A^{[\sigma/\mathbb{S}]})^{(\rho)}$.
\eprf

\blem\label{lem:inversionpi11cap}
{\rm (Inversion)}
Let 
$(\mathcal{H}_{\gamma},\Theta,\mathtt{Q})\vdash^{a}_{c,d,m}\Gamma,A^{(\rho)}$,
$A\simeq \bigwedge(A_{\iota})_{\iota\in J}$, and
$\iota\in [\rho]J$.
Then
$(\mathcal{H}_{\gamma},\Theta\cup\mathsf{k}(\iota)\cup E_{m}(\iota),\mathtt{Q}
)
\vdash^{a}_{c,d,m}
\Gamma,\left(A_{\iota}\right)^{(\rho)}$ holds.

Furthermore if $\mathsf{k}(\iota)\cup E_{m}(\iota)\subset \mathcal{H}_{\gamma}[\Theta\cap E^{\mathbb{S}}_{\rho}]$,
 then
$(\mathcal{H}_{\gamma},\Theta,\mathtt{Q}
)
\vdash^{a}_{c,d,m}
\Gamma,\left(A_{\iota}\right)^{(\rho)}$ holds.
\elem
\bprf
This is seen as in Inversion \ref{lem:inversionreg*}.
Since 
$A$ is a $\bigwedge$-formula, $A$ is not 
a $\Sigma_{1}(\mathbb{S})$-side formula in a right upper sequent of a $({\rm rfl}(\tau,e,f))$.
We need to prune some branches at $({\rm rfl}(\tau,e,f))$ 
since
$\kappa\in H_{\tau}(f,\gamma_{0},\Theta\cup\mathsf{k}(\iota)\cup E_{m}(\iota))\subset H_{\tau}(f,\gamma_{0},\Theta)$
such that 
$(\mathsf{k}(\iota)\cup E_{m}(\iota))\cap E^{\mathbb{S}}_{\tau}\subset E^{\mathbb{S}}_{\kappa}$.
\eprf

\blem\label{lem:predcereg}{\rm (Reduction)}
Let
$(\mathcal{H}_{\gamma},\Theta,\mathtt{ Q})\vdash^{a_{0}}_{c,d,m}\Gamma_{0},\lnot C_{k}^{(\rho)}$
for $1\leq k\leq n$,
and
$(\mathcal{H}_{\gamma},\Theta, \mathtt{ Q})\vdash^{a_{1}}_{c,d,m}\Pi^{(\rho)},\Gamma_{1}$,
where
$\Pi^{(\rho)}=\{C_{1}^{(\rho)},\ldots,C_{n}^{(\rho)}\}$,
$\rho=\rho_{\mathtt{Q}}$, and
$C_{k}\simeq\bigvee(C_{k,\iota})_{\iota\in J_{k}}$. Let
$\max\{\mathrm{rk}(C): C\in\Pi\}\leq c+b$ with $c>\mathbb{S}$, and 
$\{c,b\}\cup E_{m}(c, b)\subset
\mathcal{H}_{\gamma}[\Theta\cap E^{\mathbb{S}}_{\rho}]$.
Then
$(\mathcal{H}_{\gamma},\Theta,\mathtt{ Q})\vdash^{\varphi_{b}(a_{0}+a_{1})}_{c,d,m}\Gamma_{0},\Gamma_{1}$
holds.
\elem
\bprf
By main induction on $b$ with subsidiary induction on $a_{1}$.
\\
\textbf{Case 1}. $a_{1}=0$:
Then 
$(\mathcal{H}_{\gamma},\Theta, \mathtt{ Q})\vdash^{a_{1}}_{c,d,m}\Pi^{(\rho)},\Gamma_{1}$
follows
by a void $(\bigwedge)$ with a major formula in $\Gamma_{1}$.
 In what follows assume $a_{1}>0$.
 \\
 \textbf{Case 2}. The last inference in $(\mathcal{H}_{\gamma},\Theta, \mathtt{ Q})\vdash^{a_{1}}_{c,d,m}\Pi^{(\rho)},\Gamma_{1}$
 is a $(\bigvee)$ with a major formula $C_{k}^{(\rho)}$:
 We have
 $(\mathcal{H}_{\gamma},\Theta, \mathtt{ Q})\vdash^{a_{2}}_{c,d,m}\Pi^{(\rho)},C_{k,\iota}^{(\rho)},\Gamma_{1}$
 for an $\iota\in J_{k}$ and an $a_{2}<a_{1}$.
 
  SIH yields
 $(\mathcal{H}_{\gamma},\Theta, \mathtt{ Q})\vdash^{\varphi_{b}(a_{0}+a_{2})}_{c,d,m}C_{k,\iota}^{(\rho)},\Gamma_{0},\Gamma_{1}$.
 Assuming $\mathsf{k}(\iota)\subset\mathsf{k}(C_{k,\iota})$,
 we obtain $\mathsf{k}(\iota)\subset \mathcal{H}_{\gamma}[\Theta\cap E^{\mathbb{S}}_{\rho}]$ by (\ref{eq:controlder_cap2}).
Hence
$\iota\in[\rho]J_{k}$.
We obtain
 $(\mathcal{H}_{\gamma},\Theta,\mathtt{ Q})\vdash^{a_{0}}_{c,d,m}\Gamma_{0},\lnot C_{k,\iota}^{(\rho)}$
 by Inversion \ref{lem:inversionpi11cap}, and
 $(\mathcal{H}_{\gamma},\Theta,\mathtt{ Q})\vdash^{\varphi_{b}(a_{0}+a_{2})}_{c,d,m}\Gamma_{0},\lnot C_{k,\iota}^{(\rho)}$
 for $a_{0}\leq \varphi_{b}(a_{0}+a_{2})$.
 If $\mathrm{rk}(C_{k,\iota})<c$, then a $(cut)$ yields the lemma.
 Let $\mathrm{rk}(C_{k,\iota})=c+b_{1}$
 with $b_{1}<b$.
 We obtain
$\mathsf{k}(c+b_{1})\cup E_{m}(c+b_{1})\subset
\mathcal{H}_{\gamma}[\Theta\cap E^{\mathbb{S}}_{\rho}]$ by
(\ref{eq:controlder_cap2}).
 This yields $\mathsf{k}(b_{1})\cup E_{m}(b_{1})\subset \mathcal{H}_{\gamma}[\Theta\cap E^{\mathbb{S}}_{\rho}]$.
 MIH then yields
  $(\mathcal{H}_{\gamma},\Theta,\mathtt{ Q})\vdash^{\varphi_{b}(a_{0}+a_{1})}_{c,d,m}\Gamma_{0},\Gamma_{1}$
 for $\varphi_{b_{1}}(\varphi_{b}(a_{0}+a_{2})\cdot 2)<\varphi_{b}(a_{0}+a_{1})$ by $b_{1}<b$ and $a_{2}<a_{1}$.
  \\
 \textbf{Case 3}. The last inference in $(\mathcal{H}_{\gamma},\Theta, \mathtt{ Q})\vdash^{a_{1}}_{c,d,m}\Pi^{(\rho)},\Gamma_{1}$
 is a $({\rm rfl}(\tau,e,f))$:
 Then $\tau=\rho_{\mathtt{Q}}=\rho$.
 We have an ordinal $a_{2}<a_{1}$, 
sets $\Xi^{(\rho)}$ and $\Xi_{1}$ of capped formulas, and a set $\Delta$ of uncapped formulas such that
$\Xi^{(\rho)}\subset \Pi^{(\rho)}$, $\Xi_{1}\subset\Gamma_{1}$, $\forall C_{k}^{(\rho)}\in\Xi^{(\rho)}(C_{k}\in\Sigma_{1}(\mathbb{S}))$,
  $
(\mathcal{H}_{\gamma},\Theta,
\mathtt{ Q}
)\vdash^{a_{2}}_{c,d,m}\Pi^{(\rho)},\Gamma_{1}, \lnot\delta^{(\rho)}$ for each $\delta\in\Delta$, and
 $
(\mathcal{H}_{\gamma},\Theta\cup\{\sigma\},
\mathtt{ Q}^{\sigma}
)\vdash^{a_{2}}_{c,d,m}\Delta^{(\sigma)}, \Xi^{(\rho)},\Xi_{1}$ holds for each $\sigma\in H_{\rho}(f,\gamma_{0},\Theta)$.
SIH yields
 $
(\mathcal{H}_{\gamma},\Theta,
\mathtt{ Q}
)\vdash^{\varphi_{b}(a_{0}+a_{2})}_{c,d,m}\Gamma_{0},\Gamma_{1}, \lnot\delta^{(\rho)}$ for each $\delta\in\Delta$.
We obtain
 $(\mathcal{H}_{\gamma},\Theta,\mathtt{ Q} )\vdash^{\varphi_{b}(a_{0}+a_{2})+1}_{c,d,m}\Gamma_{0},\Gamma_{1}, \Xi^{(\rho)}$
 by a $({\rm rfl}(\rho,e,f))$.
For each $C_{k}^{(\rho)}\in\Xi^{(\rho)}$, we have $\mathrm{rk}(C_{k})=\mathbb{S}<c$.
$(\mathcal{H}_{\gamma},\Theta,\mathtt{ Q} )\vdash^{\varphi_{b}(a_{0}+a_{1})}_{c,d,m}\Gamma_{0},\Gamma_{1}$
follows by a series of $(cut)$'s.
\eprf

\blem\label{lem:CE}{\rm (Cut-elimination)}
Let
$(\mathcal{H}_{\gamma},\Theta,\mathtt{ Q} )\vdash^{a}_{c+b,d,m}\Gamma$
with $c>\mathbb{S}$ and
$\mathsf{k}(c)\cup E_{m}(c)\subset
\mathcal{H}_{\gamma}[\Theta\cap E^{\mathbb{S}}_{\rho_{\mathtt{Q}}}]$.
Then
$(\mathcal{H}_{\gamma},\Theta,\mathtt{ Q})\vdash^{\varphi_{b}(a)}_{c,d,m}\Gamma$
holds.
\elem
\bprf
By induction on $a$. We may assume $b>0$.
Consider the case when the last inference is a $(cut)$.
We have an $a_{0}<a$ and a $\bigvee$-formula $C^{(\rho)}$ with $\rho=\rho_{\mathtt{Q}}$ such that
$(\mathcal{H}_{\gamma},\Theta,\mathtt{ Q} )\vdash^{a_{0}}_{c+b,d,m}\Gamma,\lnot C^{(\rho)}$
and $(\mathcal{H}_{\gamma},\Theta,\mathtt{ Q} )\vdash^{a_{0}}_{c+b,d,m}C^{(\rho)},\Gamma$.
IH yields
$(\mathcal{H}_{\gamma},\Theta,\mathtt{ Q} )\vdash^{\varphi_{b}(a_{0})}_{c,d,m}\Gamma,\lnot C^{(\rho)}$
and
$(\mathcal{H}_{\gamma},\Theta,\mathtt{ Q} )\vdash^{\varphi_{b}(a_{0})}_{c,d,m}C^{(\rho)},\Gamma$.
Let $c+b_{1}=\max\{c,\mathrm{rk}(C)\}<c+b$.
Then $\mathsf{k}(b_{1})\cup E_{m}(b_{1})\subset
\mathcal{H}_{\gamma}[\Theta\cap E^{\mathbb{S}}_{\rho_{\mathtt{Q}}}]$ by (\ref{eq:controlder_cap2}).
We obtain
$(\mathcal{H}_{\gamma},\Theta,\mathtt{ Q})\vdash^{\varphi_{b_{1}}(\varphi_{b}(a_{0})\cdot 2)}_{c,d,m}\Gamma$ by
Reduction \ref{lem:predcereg}, where
$\varphi_{b_{1}}(\varphi_{b}(a_{0})\cdot 2)<\varphi_{b}(a)$.
\eprf

\blem\label{lem:picollpase}
{\rm (Collapsing)}
Let $\pi=\Omega_{\mathbb{S}+N-1-m}$ for $m<N-1$.
Assume
$\Theta\subset\mathcal{H}_{\gamma}(\psi_{\pi}(\gamma))$
 and
$
(\mathcal{H}_{\gamma},\Theta,\mathtt{Q}
)\vdash^{a}_{\pi+1,\pi+1,m}\Gamma
$ with $\Gamma\subset\Sigma(\pi)$ and $\gamma\geq\mathbb{S}$.
Then 
$(\mathcal{H}_{\hat{a}+1},\Theta,\mathtt{Q})
\vdash^{\beta}_{\beta,\beta,m+1}
\Gamma^{(\beta,\pi)}$ holds for $\hat{a}=\gamma+\omega^{a}$ and $\beta=\psi_{\pi}(\hat{a})$, where
$(A^{(\sigma)})^{(\beta,\pi)}:\equiv (A^{(\beta,\pi)})^{(\sigma)}$.
\elem
\bprf
By induction on $a$.
Let $\rho\in\mathtt{Q}$, $A^{(\rho)}\in\Gamma$ and $\alpha\in\mathsf{k}(A)\subset E^{\mathbb{S}}_{\rho}$.
We have $\{\alpha\}\cup E_{m}(\alpha) \subset\mathcal{H}_{\gamma}[\Theta\cap E^{\mathbb{S}}_{\rho}]$
by (\ref{eq:controlder_cap2}).
The assumption $\Theta\subset\mathcal{H}_{\gamma}(\psi_{\pi}(\gamma))$ yields $\alpha\in\mathcal{H}_{\gamma}(\psi_{\pi}(\gamma))$.
On the other hand we have
$E_{m}(\alpha)\cup E_{m+1}(\alpha)\subset E^{\mathbb{S}}_{\rho}$ by Proposition \ref{prp:Em}.\ref{prp:Em.6}.
We obtain
$E_{m+1}(\alpha)\subset \mathcal{H}_{\gamma}[E_{m+1}(\alpha)]=
\mathcal{H}_{\gamma}[E_{m}(\alpha)]\subset \mathcal{H}_{\gamma}[\Theta\cap E^{\mathbb{S}}_{\rho}]$
by Proposition \ref{prp:Em}.\ref{prp:Em.7}.
Hence (\ref{eq:controlder_cap2}) is enjoyed in 
$(\mathcal{H}_{\hat{a}+1},\Theta,\mathtt{Q})\vdash^{\beta}_{\beta,\beta,m+1}\Gamma^{(\beta,\pi)}$.
(\ref{eq:controlder_cap3}) is similarly seen from $E_{m+1}(\beta)\subset E_{m+1}(\gamma,a)$.
\\
\textbf{Case 1}.
First consider the case when the last inference is a $(\bigwedge)$:
A formula $A^{(\rho)}\in\Sigma(\pi)$ with $A\simeq\bigwedge\left(A_{\iota}\right)_{\iota\in J}$ 
is introduced in $\Gamma$.
For every
$\iota\in [\rho]J$ there exists an $a(\iota)<a$ 
such that
$(\mathcal{H}_{\gamma},\Theta\cup\mathsf{k}(\iota)\cup E_{m}(\iota),\mathtt{Q}
)
\vdash^{a(\iota)}_{\pi+1,\pi+1,m}\Gamma, A_{\iota}^{(\rho)}$.
We see $\mathsf{k}(\iota)\subset\psi_{\pi}(\gamma)$ from $A\in\Sigma(\pi)$ and $\mathsf{k}(A)\subset\mathcal{H}_{\gamma}[\Theta]$ 
as in Collapsing \ref{lem:Kcollpase*}.
Proposition \ref{prp:Em}.\ref{prp:Em.7} with $\gamma\geq\mathbb{S}$ then yields
$\mathcal{H}_{\gamma}[(\Theta\cup\mathsf{k}(\iota)\cup E_{m}(\iota))\cap E^{\mathbb{S}}_{\tau}]=
\mathcal{H}_{\gamma}[(\Theta\cup\mathsf{k}(\iota)\cup E_{m+1}(\iota))\cap E^{\mathbb{S}}_{\tau}]$ for every $\tau$.
Hence $(\mathcal{H}_{\gamma},\Theta\cup\mathsf{k}(\iota)\cup E_{m+1}(\iota),\mathtt{Q}
)
\vdash^{a(\iota)}_{\pi+1,\pi+1,m}\Gamma, A_{\iota}^{(\rho)}$.
On the other hand we have
$E_{m+1}(\iota)\subset\Omega_{\mathbb{S}+N-m-2}<\psi_{\pi}(\gamma)$ by Proposition \ref{prp:Em}.\ref{prp:Em.1}.
Therefore IH yields
$(\mathcal{H}_{\hat{a}+1},\Theta\cup\mathsf{k}(\iota)\cup E_{m+1}(\iota),\mathtt{Q}
)
\vdash^{\beta(\iota)}_{\beta,\beta,m+1}\Gamma^{(\beta,\pi)}, (A_{\iota}^{(\beta,\pi)})^{(\rho)}$
for every $\iota\in[\rho]J$, where $\beta(\iota)=\psi_{\pi}(\widehat{a(\iota)})$ with $\widehat{a(\iota)}=\gamma+\omega^{a(\iota)}$.
We obtain
$(\mathcal{H}_{\hat{a}+1},\Theta,\mathtt{Q})
\vdash^{\beta}_{\beta,\beta,m+1}
\Gamma^{(\beta,\pi)}$
by a $(\bigwedge)$.
\\
\textbf{Case 2}.
Second consider the case when the last inference is a $(\Sigma(\pi)\mbox{{\rm -rfl}})$:
Let $\rho=\rho_{\mathtt{Q}}$.
There exist ordinals 
$a_{\ell}, a_{r}<a$ and a formula $C\in\Sigma(\pi)$
such that 
$(\mathcal{H}_{\gamma},\Theta,\mathtt{Q})
\vdash^{a_{\ell}}_{\pi+1,\pi+1,m}\Gamma,C^{(\rho)}$,
and
$(\mathcal{H}_{\gamma},\Theta,\mathtt{Q}
)\vdash^{a_{r}}_{\pi+1,\pi+1,m}
(\lnot \exists x<\pi\, C^{(x,\pi)})^{(\rho)},\Gamma$.
As in Collapsing \ref{lem:Kcollpase*} we see from IH that
$(\mathcal{H}_{\widehat{a}+1},\Theta,\mathtt{Q})
\vdash^{\beta_{\ell}}_{\beta,\beta,m+1}
\Gamma^{(\beta,\mathbb{K})}, (C^{(\beta_{\ell},\pi)})^{(\rho)}$ and
$(\mathcal{H}_{\hat{a}+1},\Theta,\mathtt{Q})
\vdash^{\beta_{r}}_{\beta,\beta,m+1}
(\lnot C^{(\beta_{\ell},\pi)})^{(\rho)},\Gamma^{(\beta,\pi)}$, where
$\beta_{\ell}=\psi_{\pi}(\widehat{a_{\ell}})$, 
$\widehat{a_{\ell}}=\gamma+\omega^{a_{\ell}}$,
$\beta_{\ell}\in\mathcal{H}_{\widehat{a_{\ell}}+1}[\Theta]$, 
$\beta_{r}=\psi_{\pi}(\widehat{a_{r}})<\beta$ and 
$\widehat{a_{r}}=\widehat{a_{\ell}}+1+\omega^{a_{r}}<\hat{a}$.

A $(cut)$ with the cap $\rho=\rho_{\mathtt{Q}}$ yields
$(\mathcal{H}_{\hat{a}+1},\Theta,\mathtt{Q})
\vdash^{\beta}_{\beta,\beta,m+1}
\Gamma^{(\beta,\pi)}$.
\\
\textbf{Case 3}.
Third consider the case when the last inference is a $({\rm rfl}(\rho,e,f))$ with $\rho=\rho_{\mathtt{Q}}$:
We have sets $\Delta$ and $\Xi\subset\Gamma$ of formulas and $a_{0}<a$ such that
$\Delta\subset\bigvee(\pi+1)\subset\Sigma_{1}(\pi)\cup\Delta_{0}(\pi)$,
 $(\mathcal{H}_{\gamma},\Theta, \mathtt{Q}
 )\vdash^{a_{0}}_{\pi+1,\pi+1,m}\Gamma, \lnot\delta^{(\rho)}$ for each $\delta\in\Delta$,
and
$
(\mathcal{H}_{\gamma},\Theta\cup \{\sigma\},\mathtt{Q}^{\sigma}
)\vdash^{a_{0}}_{\pi+1,\pi+1,m}\Delta^{(\sigma)}, \Xi$ for each $\sigma\in H_{\rho}(f,\gamma_{0},\Theta)$.
IH with $\sigma<\mathbb{S}<\psi_{\pi}(\gamma)$ yields
$
(\mathcal{H}_{\hat{a}+1},\Theta\cup \{\sigma\},\mathtt{Q}^{\sigma}
)\vdash^{\beta_{0}}_{\beta,\beta,m+1}(\Delta^{(\beta_{0},\pi)})^{(\sigma)}, \Xi^{(\beta,\pi)}$
for $\beta_{0}=\psi_{\pi}(\widehat{a_{0}})$, $\widehat{a_{0}}=\gamma+\omega^{a_{0}}$ and
 $\mathrm{rk}(\Delta^{(\beta_{0},\pi)})<\beta$.

Let $\delta\in\Sigma_{1}(\pi)$ with $\delta\simeq\bigvee(\delta_{\iota})_{\iota\in J}$, and $\iota\in[\rho]J$ with $\mathsf{k}(\iota)\subset\beta_{0}$.
On the other hand we have
 $(\mathcal{H}_{\gamma},\Theta\cup\mathsf{k}(\iota)\cup E_{m}(\iota), \mathtt{Q}
 )\vdash^{a_{0}}_{\pi+1,\pi+1,m}\Gamma, \lnot\delta_{\iota}^{(\rho)}$
by Inversion \ref{lem:inversionpi11cap}.
As in \textbf{Case 1} we obtian
 $(\mathcal{H}_{\widehat{a_{0}}},\Theta\cup\mathsf{k}(\iota)\cup E_{m+1}(\iota), \mathtt{Q}
 )\vdash^{a_{0}}_{\pi+1,\pi+1,m}\Gamma, \lnot\delta_{\iota}^{(\rho)}$ for $E_{m+1}(\iota)\subset\Omega_{\mathbb{S}+N-m-2}$.
IH yields
$(\mathcal{H}_{\hat{a}+1},\Theta\cup\mathsf{k}(\iota)\cup E_{m+1}(\iota), \mathtt{Q}
 )\vdash^{\beta_{1}}_{\beta,\beta,m+1}\Gamma^{(\beta,\pi)}, \lnot\delta_{\iota}^{(\rho)}$
 for $\delta_{\iota}\in\Delta_{0}(\pi)$, $\beta_{1}=\psi_{\pi}(\widehat{a_{1}})$ and $\widehat{a_{1}}=\widehat{a_{0}}+\omega^{a_{0}}=\gamma+\omega^{a_{0}}\cdot 2$.
 We obtain
 $(\mathcal{H}_{\hat{a}+1},\Theta, \mathtt{Q}
 )\vdash^{\beta_{1}+1}_{\beta,\beta,m+1}\Gamma^{(\beta,\pi)}, (\lnot\delta^{(\beta_{0},\pi)})^{(\rho)}$
 by a $(\bigwedge)$.
 When $\delta\in\Delta_{0}(\pi)$, this follows from IH.
A $({\rm rfl}(\rho,e,f))$ yields  $(\mathcal{H}_{\hat{a}+1},\Theta, \mathtt{Q}
 )\vdash^{\beta}_{\beta,\beta,m+1}\Gamma^{(\beta,\pi)}$.

Other case are seen from IH.
We have $\mathrm{rk}(C)\leq\pi$ for each cut formula $C^{(\rho)}$, and the case when the last inference is a $(cut)$ is similar to \textbf{Case 3}. 
\eprf
\\

Let us embed the derivability relation $\vdash^{*}$ in $\vdash$.
Let $\Lambda_{0}=\beta=\psi_{\mathbb{K}}(\delta)$ be the fixed ordinal in Collapsing \ref{lem:Kcollpase*}
with $\delta=\hat{a}$.

\blem\label{lem:capping}{\rm (Capping)}\\
Let $\Omega_{\mathbb{S}+N-1}<\Lambda_{0}=\psi_{\mathbb{K}}(\delta)<\mathbb{K}=\Omega_{\mathbb{S}+N}$ be a strongly critical number
with $\delta\geq\mathbb{S}$.
Let
$\Gamma\subset\Delta_{0}(\mathbb{K})$ be a set of uncapped formulas.
Suppose
$
(\mathcal{H}_{\gamma},\Theta
)\vdash^{* a}_{\Lambda_{0}}\Gamma
$ with $\delta\leq\gamma\leq\gamma_{0}$.

Let
$\rho=\psi_{\mathbb{S}}^{g}(\gamma_{\rho})$
be an ordinal such that 
$\Theta\subset E^{\mathbb{S}}_{\rho}$, and
$g=m(\rho)$ is a special finite function such that
$\mathrm{supp}(g)=\{\Lambda_{0}\}$ with $g(\Lambda_{0})=\Lambda\cdot 3$ and
$\Lambda\cdot 3\leq\gamma_{0}\leq\gamma_{\rho}<\gamma_{0}+\mathbb{S}$
with 
$\Lambda=\Gamma(\Lambda_{0})$ and
$\gamma_{\rho}\in\mathcal{H}_{\gamma}[\Theta]$.
Let
$\mathtt{Q}=\{\rho\}$ be a finite family for ordinals $\Lambda,\gamma_{0}$.

Then 
$(\mathcal{H}_{\gamma},\Theta\cup\{\rho\}, \mathtt{ Q})
\vdash^{\Lambda_{0}+a}_{\Lambda_{0},\Lambda_{0},0,\Lambda_{0},\gamma_{0}}
\Gamma^{(\rho)}$ holds for $c=d=\Lambda_{0}$ and $m=0$.

\elem
\bprf   
By induction on $a$.  
We have $s(\rho)=\Lambda_{0}>\Omega_{\mathbb{S}+N-1}$, $\mathcal{H}_{\gamma_{0}}(\rho)\cap\mathbb{S}\subset\rho$, 
$\{\gamma,a,\Lambda_{0}\}\cup E_{0}(\{\gamma,a,\Lambda_{0}\})\cup\mathsf{k}(\Gamma)\cup E_{0}(\Gamma)\subset\mathcal{H}_{\gamma}[\Theta]$ 
by (\ref{eq:controlder*1}).
We obtain $\rho=\rho_{\mathtt{Q}}$.
Also 
$\Theta\cap E^{\mathbb{S}}_{\rho}=\Theta\subset E^{\mathbb{S}}_{\rho}$
by the assumption.
We obtain $E_{0}(\Gamma)\subset E^{\mathbb{S}}_{\rho}$ by Proposition \ref{prp:Emk}.\ref{prp:Emk.6}.
Hence
each of (\ref{eq:controlder_cap1}), (\ref{eq:controlder_cap2}) and (\ref{eq:controlder_cap3}) is enjoyed in
$(\mathcal{H}_{\gamma},\Theta\cup\{\rho\}, \mathtt{ Q})
\vdash^{a}_{\Lambda_{0},\Lambda_{0},0,\Lambda_{0},\gamma_{0}}
\Gamma^{(\rho)}$.
In the proof we write $(\mathcal{H}_{\gamma}, \Theta\cup\{\rho\}, \mathtt{ Q})\vdash^{a}_{\Lambda_{0}}$ for $(\mathcal{H}_{\gamma}, \Theta\cup\{\rho\}, \mathtt{ Q})\vdash^{a}_{\Lambda_{0},\Lambda_{0},0,\Lambda_{0},\gamma_{0}}$.
\\
\textbf{Case 1}.
First consider the case when the last inference is a $({\rm stbl})$: 
We have 
a $\bigwedge$-formula
$B(\mathsf{L}_{0})\in\Delta_{0}(\mathbb{S})$, and
a term $u\in Tm(\mathbb{K})$
such that
$\{\lnot B(u),\exists x\in\mathsf{L}_{\mathbb{S}}B(x)\}\subset\Gamma$, 
$\mathbb{S}\leq d=\mathrm{rk}(B(u))<\Lambda_{0}$ and $\mathsf{k}(B(u))\subset\mathcal{H}_{\gamma}[\Theta]=\mathcal{H}_{\gamma}[\Theta\cap E^{\mathbb{S}}_{\rho}]$.
We obtain 
\beqn\label{eq:capping.1}
(\mathcal{H}_{\gamma},\Theta\cup\{\rho\},\mathtt{Q})\vdash^{2d}_{0}\lnot B(u)^{(\rho)},B(u)^{(\rho)}
\eeqn
by Tautology \ref{lem:tautology.cap}.\ref{lem:tautology.cap1}.
Let $f$ be a special finite function such that $\mathrm{supp}(f)=\{\Lambda_{0}\}$ and
$f(\Lambda_{0})=\Lambda$.
Then $f_{\Lambda_{0}}=g_{\Lambda_{0}}=\emptyset$ and $f<^{\Lambda_{0}}g^{\prime}(\Lambda_{0})$
by $f(\Lambda_{0})=\Lambda<\Lambda\cdot 2=g^{\prime}(\Lambda_{0})$.
Let $\sigma\in H_{\rho}(f,\gamma_{0},\Theta\cup\{\rho\})$ with
$\Theta=(\Theta\cup\{\rho\})\cap E^{\mathbb{S}}_{\rho}\subset E^{\mathbb{S}}_{\sigma}$.
We obtain $d<\Lambda_{0}=s(f)\leq s(\sigma)$ and $\mathsf{k}(B(u))\cup\{\Lambda_{0}\}\subset E^{\mathbb{S}}_{\sigma}$,
cf.\,(\ref{eq:notationsystem.11}).
We have $\rho_{\mathtt{Q}^{\sigma}}=\sigma<\rho$ for $\mathtt{Q}^{\sigma}=\mathtt{Q}\cup\{\sigma\}$.

We have $\mathsf{k}(B(u))\subset\Lambda_{0}=\psi_{\mathbb{K}}(\delta)$ with $\mathbb{S}\leq\delta\leq\gamma$.
Hence $\mathsf{k}(B(u))\subset\mathcal{H}_{\gamma}(\psi_{\mathbb{K}}(\gamma))$.
$
(\mathcal{H}_{\gamma},\Theta\cup \{\rho,\sigma\},\mathtt{Q}^{\sigma})\vdash^{2d}_{0}\lnot B(u)^{(\sigma)},B(u^{[\sigma/\mathbb{S}]})^{(\rho)}
$ follows by Tautology \ref{lem:tautology.cap}.\ref{lem:tautology.cap3}.
A $(\bigvee)$ with $u^{[\sigma/\mathbb{S}]}\in Tm(\mathbb{S})$ yields
\beqn\label{eq:capping.2}
(\mathcal{H}_{\gamma},\Theta\cup \{\rho,\sigma\},\mathtt{Q}^{\sigma})\vdash^{2d+1}_{\Lambda_{0}}\lnot B(u)^{(\sigma)},
(\exists x\in\mathsf{L}_{\mathbb{S}}B(x))^{(\rho)}
\eeqn

An inference $({\rm rfl}(\rho,\Lambda_{0},f))$ with 
$\mathrm{rk}(B(u))<\Lambda_{0}\in\mathrm{supp}(m(\rho))$, (\ref{eq:capping.1}) and (\ref{eq:capping.2})
yields
$
(\mathcal{H}_{\gamma},\Theta\cup\{\rho\},\mathtt{ Q})
\vdash^{\Lambda_{0}}_{\Lambda_{0}}\Gamma^{(\rho)}
$.
\\
\textbf{Case 2}.  
Second the last inference $(\bigvee)$ introduces
a $\bigvee$-formula $A\in\Gamma$ with $A\simeq\bigvee\left(A_{\iota}\right)_{\iota\in J}$:
There are an $\iota\in J$ an ordinal
 $a(\iota)<a$
such that
$(\mathcal{H}_{\gamma},\Theta)\vdash^{* a(\iota)}_{\Lambda_{0}}
\Gamma,A_{\iota}$.
Assume 
$\mathsf{k}(\iota)\subset\mathsf{k}(A_{\iota})$.
We obtain $\mathcal{H}_{\gamma}(\rho)\cap\mathbb{S}\subset\rho$ 
by $\gamma\leq\gamma_{0}\leq\gamma_{\rho}$, and hence
$\mathsf{k}(\iota)\subset\mathcal{H}_{\gamma}[\Theta]\subset 
E^{\mathbb{S}}_{\rho}$ by (\ref{eq:controlder*1}), $\Theta\subset E^{\mathbb{S}}_{\rho}$
and Proposition \ref{prp:EK2}.
Hence $\iota\in[\rho]J$.
IH yields $(\mathcal{H}_{\gamma},\Theta\cup\{\rho\}, \mathtt{ Q})
\vdash^{a(\iota)}_{\Lambda_{0}}
\Gamma^{(\rho)},\left(A_{\iota} \right)^{(\rho)}$.
$(\mathcal{H}_{\gamma},\Theta\cup\{\rho\}, \mathtt{ Q})
\vdash^{a}_{\Lambda_{0}}
\Gamma^{(\rho)}$ 
follows from a $(\bigvee)$.
\\
\textbf{Case 3}.
Third the last inference $(\bigwedge)$ introduces a $\bigwedge$-formula $A\in\Gamma$ with
$A\simeq\bigwedge\left(A_{\iota}\right)_{\iota\in J}$:
For every
$\iota\in J$ there exists an $a(\iota)<a$ 
such that
$(\mathcal{H}_{\gamma},\Theta\cup\mathsf{k}(\iota)\cup E_{0}(\iota)
)
\vdash^{* a(\iota)}_{\Lambda_{0}}
\Gamma,A_{\iota}$.
IH yields 
$(\mathcal{H}_{\gamma}, \Theta\cup\{\rho\}\cup\mathsf{k}(\iota)\cup E_{0}(\iota), \mathtt{ Q}
)
\vdash^{a(\iota)}_{\Lambda_{0}}
\Gamma^{(\rho)},\left(A_{\iota}\right)^{(\rho)}
$
for each $\iota\in [\rho]J$,
where $\mathsf{k}(\iota)\cup E_{0}(\iota)\subset E^{\mathbb{S}}_{\rho}$.
We obtain
$(\mathcal{H}_{\gamma},\Theta\cup\{\rho\}, \mathtt{ Q})\vdash^{a}_{\Lambda_{0}}
\Gamma^{(\rho)}$ by a $(\bigwedge)$.

Other cases 
$(cut)$ or $(\Sigma(\pi)\mbox{{\rm -rfl}})\, (\pi<\Lambda_{0})$ are seen from IH.
Each uncapped cut formula $C$ as well as
formulas $C$ and $\lnot \exists x<\pi\,C^{(x,\pi)}$ in $(\Sigma(\pi)\mbox{{\rm -rfl}})$
puts on the cap $\rho=\rho_{\mathtt{Q}}$ with $\mathrm{rk}(C)<\Lambda_{0}$.
\eprf

\subsection{Reducing ranks}\label{subsec:elimrfl}

In this subsection, ranks in inferences $({\rm rfl}(\rho,e,f))$ 
are lowered to $\mathbb{S}$ in operator controlled derivations $\mathcal{D}_{1}$ of $\Sigma_{1}$-sentences 
$\theta^{L_{\Omega}}$ over $\Omega$.
Let $\mathcal{D}_{2}$ be a derivation such that every formula occurring in it is in 
$\Sigma_{1}(\mathbb{S})\cup\Pi_{1}(\mathbb{S})$.
We see in Lemma \ref{lem:predcereg_S} that inferences $({\rm rfl}(\rho,e,f))$ 
are removed from $\mathcal{D}_{2}$, where each capped formula $A^{(\rho)}$
becomes the uncapped formula $A$ in Lemma \ref{lem:main}.
To have $\mathtt{ Q}\subset\mathcal{H}_{\gamma_{1}}[\Theta]$ for finite families $\mathtt{ Q}$,
we break through the threshold $\gamma_{0}$ in the sense that
 $\gamma_{1}\geq\gamma_{0}+\mathbb{S}$.
 We need the condition (\ref{eq:controlder_cap3}) to be enjoyed in Recapping \ref{mlem:singlemainl.1}.
Everything has to be done inside $\mathcal{H}_{\gamma_{0}}[\Theta]$
except ordinals in $\mathtt{ Q}$ until the rank is lowered to $\mathbb{S}$.
Our goal in this subsection is to transform derivations $\mathcal{D}_{1}$ to $\mathcal{D}_{2}$.
For this $\mathcal{D}_{1}$ is first transformed to a derivation $\mathcal{E}$ in which
every capped formula is in $\Sigma(\Omega_{\mathbb{S}+N-m})$.
Then Collapsing \ref{lem:picollpase} yields
a derivation in rank less than $\Omega_{\mathbb{S}+N-m}$.
Iterating this process, we arrive at a derivation $\mathcal{D}_{2}$.

In the following Recapping \ref{mlem:singlemainl.1} we show that inferences $({\rm rfl}(\rho,e,f))$
can be replaced by a series of $(cut)$'s.
Here caps $\rho$ are replaced by smaller caps
$\kappa\in H_{\rho}(g_{1},\gamma_{0},\Theta)$ for $g_{1}=h^{b}(m(\rho);a_{1})$
with an ordinal $a_{1}$, and other inferences
$({\rm rfl}(\kappa,b,h))$ are introduced for smaller ranks $b<e$.

When $b\leq \mathbb{S}+1$, we obtain 
$\delta^{(\sigma)}\in\Sigma_{1}(\mathbb{S})$
for the formula $\delta$ to be reflected in inferences $({\rm rfl}(\kappa,b,h))$.
However it may be the case $s(\kappa)\leq\mathbb{S}+1<s(h)\leq s(\sigma)$
for $\sigma$ in the resolvent class of $({\rm rfl}(\rho,e,f))$.
Therefore we need to replace inferences $({\rm rfl}(\sigma,e,g))$ in higher ranks by a series of $(cut)$'s.
We have to iterate the replacement inside $\mathcal{H}_{\gamma_{0}}[\Theta]$, and
an induction on $\sigma<\rho$ must be avoided.

For the ordinals $on_{m}(\mathtt{Q})$ in Definition \ref{df:onQ}
we obtain $on_{m}(\mathtt{Q})\in\mathcal{H}_{\gamma}[\Theta]$ 
if $SC_{\Lambda}(f)\subset\mathcal{H}_{\gamma}[\Theta]$,
and we see $on_{m}(\mathtt{Q}^{\lambda})<on_{m}(\mathtt{Q})$ from Proposition \ref{prp:hstepdown}.\ref{prp:hstepdown.1},
where $\lambda$ is a replacement for $\sigma$ and $\mathtt{Q}^{\lambda}=\mathtt{Q}\cup\{\lambda\}$.
By induction on the ordinals $on_{m}(\mathtt{Q})$ we see  in Lemma \ref{mlem:singlemainl_S} that
the rank is lowered to $\Omega_{\mathbb{S}+N-m}$.
This ends a rough sketch of the removals of inferences $({\rm rfl}(\rho,e,f))$,
and the details follow.

\blem\label{mlem:singlemainl.1}
{\rm (Recapping)}\\
Suppose
$
(\mathcal{H}_{\gamma},\Theta_{\rho}, \mathtt{ Q})
\vdash^{a}_{c,d,m,\Lambda_{0},\gamma_{0}}
\Pi,\Gamma^{(\rho)}
$, where $\Theta_{\rho}=\Theta\cup\{\rho\}$, $a<\Lambda$, $c=\Omega_{\mathbb{S}+N-m-1}+1<s(\rho)$, 
$c<d\leq\Lambda_{0}<\Lambda$,
$\rho\in\mathtt{Q}$ with $\Theta\cap E^{\mathbb{S}}_{\rho}\subset E^{\mathbb{S}}_{\rho_{\mathtt{Q}}}$,
and
$\Pi$ is a set of formulas such that $\tau\neq \rho$ for each $A^{(\tau)}\in\Pi$.

Let $b$ be an ordinal 
such that $c\leq b<s(\rho)\leq \Lambda_{0}$ and
$\Gamma\subset\bigvee(b)$.
Assume that $b\in\mathcal{H}_{\gamma}[\Theta\cap E^{\mathbb{S}}_{\rho}]$ and
$SC_{\Lambda}(g)\subset\mathcal{H}_{\gamma}[\Theta\cap E^{\mathbb{S}}_{\rho}]$ for $g=m(\rho)$.

Let 
$\kappa\in H_{\rho}(g_{1},\gamma_{0},\Theta)$ with 
$g_{1}=h^{b}(g;\varphi_{d}(a))$, and
$\mathtt{Q}^{[\kappa/\rho]}=\{\tau\in\mathtt{Q}:\tau\neq\rho\}\cup\{\kappa\}$.
Then 
\begin{equation}\label{eq:1.1.f}
(\mathcal{H}_{\gamma},\Theta_{\kappa}, \mathtt{Q}^{[\kappa/\rho]}
)
\vdash^{\varphi_{d}(a)}_{c,d,m,\Lambda_{0},\gamma_{0}}
\Pi,\Gamma^{(\kappa)}
\end{equation}
holds, where $\Theta_{\kappa}=\Theta\cup\{\kappa\}$.

\elem
\bprf
By induction on $a$. 
We have $\varphi_{d}(a)<\Lambda$ and $g_{1}:\Lambda\to\Gamma(\Lambda)$ with $SC_{\Lambda}(g_{1})\subset\Lambda$.
Let
$\kappa\in H_{\rho}(g_{1},\gamma_{0},\Theta)$.
By Definition \ref{df:resolvent} we have 
$\Theta\cap E^{\mathbb{S}}_{\rho}
\subset E^{\mathbb{S}}_{\kappa}$.
On the other hand we have
$\{a,d,b\}\cup\mathsf{k}(\Gamma)\subset\mathcal{H}_{\gamma}[\Theta\cap E^{\mathbb{S}}_{\rho}]$
by $\rho\geq\rho_{\mathtt{Q}}$, (\ref{eq:controlder_cap2}), (\ref{eq:controlder_cap3}) and the assumption.
Moreover $\rho_{\mathtt{Q}^{[\kappa/\rho]}}=\kappa$ if $\rho=\rho_{\mathtt{Q}}$.
Otherwise $\rho_{\mathtt{Q}^{[\kappa/\rho]}}=\rho_{\mathtt{Q}}$, cf.\,Definition \ref{df:cap}.\ref{df:cap.2}.
In each case we obtain $\Theta\cap E^{\mathbb{S}}_{\rho}\subset E^{\mathbb{S}}_{\rho_{\mathtt{Q}^{[\kappa/\rho]}}}$. 
Hence each of (\ref{eq:controlder_cap2}) and (\ref{eq:controlder_cap3}) 
is enjoyed in (\ref{eq:1.1.f}).

Also $\{b,d,a\}\cup SC_{\Lambda}(g)\subset E^{\mathbb{S}}_{\kappa}$ by $\{b,d,a\}\cup SC_{\Lambda}(g) \subset\mathcal{H}_{\gamma}[\Theta\cap E^{\mathbb{S}}_{\rho}]$,
Proposition \ref{prp:EK2} with $\gamma\leq\gamma_{0}$.
Hence 
$SC_{\Lambda}(h^{b}(g;\varphi_{d}(a)))\subset E^{\mathbb{S}}_{\kappa}$, cf.\,(\ref{eq:notationsystem.11}).

In the proof let us suppress the fourth and fifth subscripts, and we write
$
(\mathcal{H}_{\gamma},\Theta, \mathtt{ Q})
\vdash^{a}_{c,d,m}
\Pi,\Gamma^{(\rho)}
$
for
$
(\mathcal{H}_{\gamma},\Theta, \mathtt{ Q})
\vdash^{a}_{c,d,m,\Lambda_{0},\gamma_{0}}
\Pi,\Gamma^{(\rho)}
$.
\\
\textbf{Case 1}. 
First consider the case when the last inference is a $({\rm rfl}(\tau,e,f))$:
Then $\tau=\rho_{\mathtt{Q}}$.
If $\tau<\rho$, then (\ref{eq:1.1.f}) follows from IH.

In what follows assume $\rho_{\mathtt{Q}}=\tau=\rho$.
Then $\rho_{\mathtt{Q}^{[\kappa/\rho]}}=\kappa$.
We have 
a finite set $\Delta\subset\bigvee(e)$ of formulas with $b_{0}=\max\{c,\mathrm{rk}(\Delta)\}<d$
and an ordinal $a_{0}<a$ such that 
\[
(\mathcal{H}_{\gamma},\Theta_{\rho},
\mathtt{ Q}
)
\vdash^{a_{0}}_{c,d,m}\Pi,\Gamma^{(\rho)}, 
\lnot\delta^{(\rho)}
\]
for each $\delta\in\Delta$.
Inversion \ref{lem:inversionpi11cap} yields
\begin{equation}\label{eq:Case1left_inv}
(\mathcal{H}_{\gamma},\Theta_{\rho}\cup\mathsf{k}(\iota)\cup E_{m}(\iota),
\mathtt{ Q}
)
\vdash^{a_{0}}_{c,d,m}\Pi,\Gamma^{(\rho)}, 
\lnot\delta_{\iota}^{(\rho)}
\end{equation}
for each $\iota\in
[\rho]J$.
We have $(\lnot\delta_{\iota})\in\bigvee(b_{0})$.
Then $b_{0}< e\leq s(\rho)$ and $b_{0}<d$.
On the other hand we have 
\begin{equation}\label{eq:KppiNlowerCase1b}
(\mathcal{H}_{\gamma},\Theta_{\rho\sigma},
\mathtt{Q}^{\sigma}
)
\vdash^{a_{0}}_{c,d,m}
\Delta^{(\sigma)}, \Xi^{(\rho)},\Pi_{1}
\end{equation}
where $\Xi^{(\rho)}\in\Gamma^{(\rho)}\cap\Sigma_{1}(\mathbb{S})$, $\Pi_{1}\subset\Pi$, $\Theta_{\rho\sigma}=\Theta\cup \{\rho,\sigma\}$ and
$\sigma\in H_{\rho}(f,\gamma_{0},\Theta)$.
We have $\Theta\cap E^{\mathbb{S}}_{\rho}\subset E^{\mathbb{S}}_{\sigma}$ for $\sigma=\rho_{\mathtt{Q}^{\sigma}}$.
When $\max\{s(\rho),s(f)\}\leq\Omega_{\mathbb{S}+N-1-m}+1$, we have (\ref{eq:KppiNlowerCase1b}) for every $\sigma$ such that $m(\sigma)=f$.

$f$ is a finite function such that $e\in{\rm supp}(g)$ and
\beqn\label{eq:recapping_Pi11_r2}
f_{e}=g_{e} \,\&\, f<^{e}g^{\prime}(e) \,\&\, 
SC_{\Lambda}(f)\subset\mathcal{H}_{\gamma}[\Theta\cap E^{\mathbb{S}}_{\rho}]
\eeqn
\textbf{Case 1.1}. $b_{0}\leq b$: 
Let $\delta\in\Delta$, and $\iota\in J$ for $\delta\simeq\bigvee(\delta_{\iota})_{\iota\in J}$.
We have $\mathrm{rk}(\lnot\delta_{\iota})<\mathrm{rk}(\lnot\delta)\leq b$, and 
$\mathrm{rk}(\Gamma\cup\{\lnot\delta_{\iota}\})< b$.
Let $\iota\in[\kappa]J$.
We obtain $\iota\in[\rho]J$ by
$\kappa\leq\rho$, and $\mathsf{k}(\iota)\subset E^{\mathbb{S}}_{\kappa}$ by $\iota\in[\kappa]J$.
Hence $\kappa\in H_{\rho}(g_{1},\gamma_{0},\Theta\cup\mathsf{k}(\iota)\cup E_{m}(\iota))$.
By IH with (\ref{eq:Case1left_inv}) 
we obtain 
$(\mathcal{H}_{\gamma},\Theta_{\kappa}\cup\mathsf{k}(\iota)\cup E_{m}(\iota),
\mathtt{Q}^{[\kappa/\rho]}
)
\vdash^{\varphi_{d}(a_{0})}_{c,d,m}
\Pi,\Gamma^{(\kappa)},\lnot\delta_{\iota}^{(\kappa)}$.
A $(\bigwedge)$ yields
\begin{equation}\label{eq:Case1cpi11}
(\mathcal{H}_{\gamma},\Theta_{\kappa},
\mathtt{ Q}^{[\kappa/\rho]}
)
\vdash^{\varphi_{d}(a_{0})+1}_{c,d,m}
\Pi,\Gamma^{(\kappa)},\lnot\delta^{(\kappa)}
\end{equation}

Let $d_{1}=\min\{b,e\}$.
We claim that
\begin{equation}\label{eq:zigpi11.111}
f_{d_{1}}=(g_{1})_{d_{1}} \,\&\,
f<^{d_{1}}g_{1}^{\prime}(d_{1}) 
\end{equation}
We have $g_{1}=h^{b}(g;\varphi_{d}(a))$.
If $d_{1}=e\leq b$, then $(g_{1})_{e}=g_{e}=f_{e}$ and $g^{\prime}(e)\leq g_{1}^{\prime}(e)$.
We obtain the claim by Proposition \ref{prp:idless}.
If $d_{1}=b< e$, then the claim follows from Proposition \ref{prp:hstepdown}.\ref{prp:hstepdown.2}.

Let $\sigma\in H_{\kappa}(f,\gamma_{0},\Theta)$.
Then $\sigma\in H_{\rho}(f,\gamma_{0},\Theta)$ by 
$\kappa<\rho$, $\Theta\cap E^{\mathbb{S}}_{\rho}\subset E^{\mathbb{S}}_{\kappa}$.
By IH with (\ref{eq:KppiNlowerCase1b}) we obtain for $\Theta_{\kappa\sigma}=\Theta\cup\{\kappa,\sigma\}$
\begin{equation}\label{eq:1.1.fl}
(\mathcal{H}_{\gamma},\Theta_{\kappa\sigma},
\mathtt{Q}^{[\kappa/\rho]}\cup\{\sigma\}
)
\vdash^{\varphi_{d}(a_{0})}_{c,d,m}
\Delta^{(\sigma)}, \Xi^{(\kappa)},\Pi_{1}
\end{equation}

(\ref{eq:1.1.f})
follows by an inference $({\rm rfl}(\kappa,d_{1},f))$ with 
(\ref{eq:zigpi11.111}), (\ref{eq:Case1cpi11}) and (\ref{eq:1.1.fl}).

If $\max\{s(\kappa), s(f)\}\leq\Omega_{\mathbb{S}+N-m-1}+1$, then 
the inference is degenerated, and it suffices to have (\ref{eq:1.1.fl}) for ordinals $\sigma$ with $m(\sigma)=f$.
\\
\textbf{Case 1.2}. 
$b< b_{0}$: 
We obtain $\mathrm{rk}(\Delta)\in\mathcal{H}_{\gamma}[\Theta\cap E^{\mathbb{S}}_{\rho}]$  by
(\ref{eq:controlder_cap2}).
We have $b< b_{0}< e\leq s(\rho)$ and $b_{0}<d$.
Let $g_{0}=h^{b_{0}}(g;\varphi_{d}(a_{0}))*f^{b_{0}+1}$, and
$\sigma\in H_{\rho}(g_{0},\gamma_{0},\Theta)$.
We have 
$g_{0}\leq m(\sigma)$ 
by Definition \ref{df:resolvent}.
We see $SC_{\Lambda}(g_{0})\subset\mathcal{H}_{\gamma}[\Theta\cap E^{\mathbb{S}}_{\kappa}]$
from $\{b_{0},d,a_{0}\}\cup SC_{\Lambda}(g)\cup SC_{\Lambda}(f)\subset\mathcal{H}_{\gamma}[\Theta\cap E^{\mathbb{S}}_{\kappa}]$.
As in \textbf{Case 1.1}, we obtain 
$(\mathcal{H}_{\gamma},\Theta_{\sigma},
\mathtt{Q}^{[\sigma/\rho]}
)
\vdash^{\varphi_{d}(a_{0})+1}_{c,d,m}
\Pi,\Gamma^{(\sigma)},\lnot\delta^{(\sigma)}$
by IH and (\ref{eq:Case1left_inv}).

Let $\sigma\in H_{\kappa}(g_{0},\gamma_{0},\Theta)=H_{\rho}(g_{0},\gamma_{0},\Theta)\cap\kappa$. 
Then for $\mathtt{Q}^{[\kappa/\rho]}\cup\{\sigma\}=\mathtt{Q}^{[\sigma/\rho]}\cup\{\kappa\}$, we have
$\sigma=\rho_{\mathtt{Q}^{[\sigma/\rho]}}=\rho_{\mathtt{Q}^{[\sigma/\rho]}\cup\{\kappa\}}=\rho_{\mathtt{Q}^{[\kappa/\rho]}\cup\{\sigma\}}$.
Hence (\ref{eq:controlder_cap3}) holds by adding the ordinal $\kappa$ to $\mathtt{Q}^{[\sigma/\rho]}$, and we obtain
\begin{equation}\label{eq:L4.10case1rxi}
(\mathcal{H}_{\gamma},\Theta_{\kappa\sigma},
\mathtt{Q}^{[\kappa/\rho]}\cup\{\sigma\}
)
\vdash^{\varphi_{d}(a_{0})+1}_{c,d,m}
\Pi,\Gamma^{(\sigma)},\lnot\delta^{(\sigma)}
\end{equation}

We have
$f\leq g_{0}$.
Hence $\sigma\in H_{\rho}(f,\gamma_{0},\Theta)$.
We obtain by IH and (\ref{eq:KppiNlowerCase1b})
\begin{equation}\label{eq:KppiNlowerCase1b.12}
(\mathcal{H}_{\gamma},\Theta_{\kappa\sigma},
\mathtt{Q}^{[\kappa/\rho]}\cup\{\sigma\}
)
\vdash^{\varphi_{d}(a_{0})}_{c,d,m}
\Delta^{(\sigma)}, \Xi^{(\kappa)},\Pi_{1}
\end{equation}

We have $\sigma=\rho_{\mathtt{Q}^{[\kappa/\rho]}\cup\{\sigma\}}$, and
$\mathrm{rk}(\delta)\leq b_{0}\leq c+b_{0}$ with $\{c,b_{0}\}\cup E(c,b_{0})\subset\mathcal{H}_{\gamma}[\Theta\cap E^{\mathbb{S}}_{\sigma}]$
by $\Theta\cap E^{\mathbb{S}}_{\rho}\subset  E^{\mathbb{S}}_{\sigma}$.
Reduction \ref{lem:predcereg} with  (\ref{eq:L4.10case1rxi}) and (\ref{eq:KppiNlowerCase1b.12})
yields for $H_{\kappa}(g_{0},\gamma_{0},\Theta) =H_{\rho}(g_{0},\gamma_{0},\Theta)\cap\kappa $
\begin{equation}\label{eq:L4.10case134}
\forall\sigma\in H_{\kappa}(g_{0},\gamma_{0},\Theta) 
\left[
(\mathcal{H}_{\gamma},\Theta_{\kappa\sigma},
\mathtt{Q}^{[\kappa/\rho]}\cup\{\sigma\}
)
\vdash^{a_{1}}_{c,d,m}
\Pi,\Gamma^{(\sigma)}, \Xi^{(\kappa)}
\right]
\end{equation}
for
$2b\leq a_{1}=\varphi_{b_{0}}(\varphi_{d}(a_{0})\cdot 2)<\varphi_{d}(a)$
 by $b<b_{0}<d$ and $a_{0}<a$.

On the other, Tautology \ref{lem:tautology.cap}.\ref{lem:tautology.cap1} yields
for each $\theta\in\Gamma$ 
\begin{equation}\label{eq:L4.10case1.1b}
(\mathcal{H}_{\gamma},\Theta_{\kappa}, \mathtt{ Q}^{[\kappa/\rho]}
)\vdash^{2b}_{0,0,m}
\Gamma^{(\kappa)},\lnot \theta^{(\kappa)}
\end{equation}

For $g_{1}=h^{b}(g; \varphi_{d}(a))$, we obtain
$(g_{0})_{b}=g_{b}=(g_{1})_{b}$ and $g_{0}<^{b}g_{1}^{\prime}(b)$
by Proposition \ref{prp:hstepdown}.\ref{prp:hstepdown.4}.

By an inference rule 
$({\rm rfl}(\kappa,b,g_{0}))$ 
with 
its resolvent class 
$H_{\kappa}(g_{0},\gamma_{0},\Theta)$,
we conclude (\ref{eq:1.1.f}) by  (\ref{eq:L4.10case1.1b}), (\ref{eq:L4.10case134})
 with 
$\kappa=\rho_{\mathtt{Q}^{[\kappa/\rho]}}$ and $\mathrm{rk}(\Gamma)<b< d$.
\\
\textbf{Case 2}.
Second consider the case when the last inference $(\bigvee)$ introduces a $\bigvee$-formula
$B$:
Let
$B\equiv A^{(\rho)}\in\Gamma^{(\rho)}$
with $A\simeq\bigvee\left(A_{\iota}\right)_{\iota\in J}$.
We have 
$
(\mathcal{H}_{\gamma},\Theta, \mathtt{ Q}
)
\vdash^{a_{0}}_{c,d,m}
\Pi,\Gamma^{(\rho)},\left(A_{\iota}\right)^{(\rho)}
$,
where $\mathsf{k}(A_{\iota})\subset\mathcal{H}_{\gamma}[\Theta\cap E^{\mathbb{S}}_{\rho}]$.
Assuming $\mathsf{k}(\iota)\subset\mathsf{k}(A_{\iota})$, we obtain 
$\mathsf{k}(\iota)\subset E^{\mathbb{S}}_{\kappa}$, i.e., $\iota\in[\kappa]J$ by $\Theta\cap E^{\mathbb{S}}_{\rho}\subset E^{\mathbb{S}}_{\kappa}$.
IH 
yields
$
(\mathcal{H}_{\gamma},\Theta_{\kappa},
\mathtt{ Q}^{[\kappa/\rho]}
)\vdash^{\varphi_{d}(a_{0})}_{c,d,m}
\Pi,\Gamma^{(\kappa)},
\left(A_{\iota}\right)^{(\kappa)}
$ for $\mathrm{rk}(\Gamma\cup\{A_{\iota}\})=\mathrm{rk}(\Gamma)$.
We obtain (\ref{eq:1.1.f})
by a $(\bigvee)$.

Other cases are seen from IH.
Note that
for each capped cut formula $C^{(\rho_{\mathtt{Q}})}$, we have $\mathrm{rk}(C)<c=\Omega_{\mathbb{S}+N-1-m}+1\leq b$, and
for a minor formula $(\forall x\in u\,B(x))^{(\rho_{\mathtt{Q}})}$ of a $(\Sigma(\pi)\mathrm{-rfl})$,
$\mathrm{rk}(B(v))=\pi<\Omega_{\mathbb{S}+N-1-m}+1$ for $|v|<|u|<\pi$.
These formulas put on the cap $\rho_{\mathtt{ Q}^{[\kappa/\rho]}}$.
Also $[\kappa]J\subset [\rho]J$ for $(\bigwedge)$.
\eprf

\bdf\label{df:onQ}
{\rm
For a finite family $\mathtt{Q}$ with $\rho=\rho_{\mathtt{Q}}$ and $g=m(\rho)$, let
\[
on_{m}(\mathtt{Q})=g(\Omega_{\mathbb{S}+N-m-1}+1)
.\]
}
\edf

\blem\label{mlem:singlemainl_S}
Let
$
(\mathcal{H}_{\gamma},\Theta, \mathtt{ Q}
)
\vdash^{a}_{c,d,m,\Lambda_{0},\gamma_{0}}
\Pi
$, where $a<\Lambda$, $c=\Omega_{\mathbb{S}+N-m-1}+1<d<\Lambda$, and $\Pi\subset\Delta_{0}(\mathbb{K})$.
Assume
$s(\rho_{\mathtt{Q}})\leq\Omega_{\mathbb{S}+N-m-1}+1$ and
$SC_{\Lambda}(m(\rho_{\mathtt{Q}}))\subset \mathcal{H}_{\gamma}[\Theta\cap E^{\mathbb{S}}_{\rho_{\mathtt{Q}}}]$.
Then the following holds for $\hat{a}=\varphi_{d+\eta}(a)$ and $\eta=on_{m}(\mathtt{Q})$
\begin{equation}\label{eq:singlemainl_S1}
(\mathcal{H}_{\gamma},\Theta, \mathtt{Q}
)
\vdash^{\hat{a}}_{c,c,m,\Lambda_{0},\gamma_{0}}
\Pi
\end{equation}

\elem
\bprf
By main induction on $\eta=on_{m}(\mathtt{Q})$ with subsidiary induction on $a$.

We have 
$SC_{\Lambda}(m(\rho_{\mathtt{Q}}))\subset \mathcal{H}_{\gamma}[\Theta\cap E^{\mathbb{S}}_{\rho_{\mathtt{Q}}}]$ by the assumption,
and
 $\{a,d\}
\subset
\mathcal{H}_{\gamma}[\Theta\cap E^{\mathbb{S}}_{\rho_{\mathtt{Q}}}]$ by (\ref{eq:controlder_cap3}).
We obtain $\{\hat{a}\}\cup E(a,d,\eta)\subset
\mathcal{H}_{\gamma}[\Theta\cap E^{\mathbb{S}}_{\rho_{\mathtt{Q}}}]$ for (\ref{eq:controlder_cap3}).
In the proof let us omit the fourth and fifth subscripts in $\vdash^{a}_{c,d,m,\Lambda_{0},\gamma_{0}}$.

Consider the case when the last inference is a $({\rm rfl}(\rho,e,f))$ with $\rho=\rho_{\mathtt{Q}}$:
Let $g=m(\rho)$.
We have 
a finite set $\Delta\subset\bigvee(c)$ of formulas with $\mathrm{rk}(\Delta)<e\leq s(\rho)\leq c=\Omega_{\mathbb{S}+N-m-1}+1\leq e$
and an ordinal $a_{0}<a$ such that 
\begin{equation}\label{eq:Case1left_S}
(\mathcal{H}_{\gamma},\Theta,
\mathtt{Q}
)
\vdash^{a_{0}}_{c,d,m}\Pi,
\lnot\delta^{(\rho)}
\end{equation}
for each $\delta\in\Delta$.
$f$ is a finite function such that $c=e\in{\rm supp}(g)$ and
$f_{c}=g_{c}$,  $f<^{c}g^{\prime}(c)$ and
$SC_{\Lambda}(f)\subset\mathcal{H}_{\gamma}[\Theta\cap E^{\mathbb{S}}_{\rho}]$.

SIH yields for $\widehat{a_{0}}=\varphi_{d+\eta}(a_{0})$
\begin{equation}\label{eq:L4.10case1rxi_S}
(\mathcal{H}_{\gamma},\Theta,
\mathtt{Q}
)
\vdash^{\widehat{a_{0}}}_{c,c,m}
\Pi,\lnot\delta^{(\rho)}
\end{equation}

On the other hand we have
\begin{equation}\label{eq:singlemainl_SS}
(\mathcal{H}_{\gamma},\Theta_{\sigma},
\mathtt{ Q}^{\sigma}
)
\vdash^{a_{0}}_{c,d,m}\Pi_{1},
\Delta^{(\sigma)} 
\end{equation}
for every $\sigma\in H_{\rho}(f,\gamma_{0},\Theta)$
 with $m(\sigma)=f$,
$SC_{\Lambda}(m(\sigma))=SC_{\Lambda}(f)\subset\mathcal{H}_{\gamma}[\Theta\cap E^{\mathbb{S}}_{\rho}]\subset E^{\mathbb{S}}_{\rho}$, where
$\Pi_{1}\subset\Pi$ and $\forall A^{(\rho)}\in\Pi_{1}(A\in\Sigma_{1}(\mathbb{S}))$.

If $s(f)\leq c$, then SIH yields for each 
$\sigma\in H_{\rho}(f,\gamma_{0},\Theta)$
\begin{equation}\label{eq:KppiNlowerCase1b.12_S}
(\mathcal{H}_{\gamma},\Theta_{\sigma},
\mathtt{Q}^{\sigma}
)
\vdash^{\widehat{a_{0}}}_{c,c,m}
\Pi_{1}, \Delta^{(\sigma)}
\end{equation}
We obtain (\ref{eq:singlemainl_S1}) by a degenerated $({\rm rfl}(\rho,e,f))$ with (\ref{eq:L4.10case1rxi_S}) and (\ref{eq:KppiNlowerCase1b.12_S}).

Assume $s(\sigma)=s(f)>c$.
Let $f_{0}=h^{c}(f;\varphi_{d}(a_{0}))$, where
$SC_{\Lambda}(f_{0})\subset \mathcal{H}_{\gamma}[\Theta\cap E^{\mathbb{S}}_{\rho}]\subset E^{\mathbb{S}}_{\rho}$.
Let  $\lambda\in H_{\rho}(f_{0},\gamma_{0},\Theta)$ with $m(\lambda)=f_{0}$, and
$\sigma=\psi_{\rho}^{f}(\nu+\alpha)$ for $\nu=b(\rho)$ and
$\alpha=\max(\{\lambda\}\cup ((\Theta\cup SC_{\Lambda}(f))\cap E^{\mathbb{S}}_{\rho})\}$.
Then we see $\nu+\alpha<\gamma_{0}+\mathbb{S}$ from $\nu<\gamma_{0}+\mathbb{S}$ and $\alpha<\mathbb{S}$.
Hence
$\sigma\in H_{\rho}(f,\gamma_{0},\Theta)$ and
$\lambda\in H_{\sigma}(f_{0},\gamma_{0},\Theta)$.
Recapping \ref{mlem:singlemainl.1} with (\ref{eq:singlemainl_SS}) yields
\begin{equation}\label{eq:KppiNlowerCase1b.12_SS}
(\mathcal{H}_{\gamma},\Theta_{\lambda},
\mathtt{Q}^{\lambda}
)
\vdash^{\varphi_{d}(a_{0})}_{c,d,m}
\Pi_{1}, \Delta^{(\lambda)}
\end{equation}
where $SC_{\Lambda}(m(\lambda))=SC_{\Lambda}(f_{0})\subset\mathcal{H}_{\gamma}[\Theta\cap E^{\mathbb{S}}_{\rho}]$, 
$\mathtt{Q}^{\lambda}=\mathtt{Q}\cup\{\lambda\}$, $\rho_{\mathtt{Q}^{\lambda}}=\lambda<\sigma<\rho=\rho_{\mathtt{Q}}$, and
$s(f_{0})=c$.
We obtain 
$\xi=on_{m}(\mathtt{Q}^{\lambda})=f_{0}(c)<g(c)=on_{m}(\mathtt{Q})=\eta$
by Proposition \ref{prp:hstepdown}.\ref{prp:hstepdown.1}.

MIH then yields
\begin{equation}\label{eq:KppiNlowerCase1b.12_SSS}
(\mathcal{H}_{\gamma},\Theta_{\lambda},
\mathtt{Q}^{\lambda}
)
\vdash^{\widehat{a_{1}}}_{c,c,m}
\Pi_{1}, \Delta^{(\lambda)}
\end{equation}
where $\widehat{a_{1}}=\varphi_{d+\xi}(\varphi_{d}(a_{0}))<\varphi_{d+\eta}(a)$ by $\xi<\eta$ and $a_{0}<a$.

A degenerated $({\rm rfl}(\rho,c,f_{0}))$
with (\ref{eq:L4.10case1rxi_S}) and (\ref{eq:KppiNlowerCase1b.12_SSS})
yields (\ref{eq:singlemainl_S1}).

Other cases are seen from SIH.
\eprf

\subsection{A third calculus}

Combining Cut-elimination \ref{lem:CE}, Lemma \ref{mlem:singlemainl_S} and Collapsing \ref{lem:picollpase},
the rank of formulas in derivations is lowered to $\mathbb{S}$.
In other words, every formula occurring in derivations is in $\Sigma_{1}(\mathbb{S})\cup\Pi_{1}(\mathbb{S})\cup\Delta_{0}(\mathbb{S})$.
We are going to eliminate inferences $({\rm rfl}(\rho,e,f))$ and $\Sigma_{1}(\mathbb{S})$-formulas, and throw caps away.
In doing so, it is better to shift the calculus from $\vdash$ in subsection \ref{subsec:derivationcap} to a third one $\vdash^{\diamond}$.

An uncapped formula $A$ is denoted by $A^{(\mathtt{u})}$, and let $[\mathtt{u}]J=J$ for $E^{\mathbb{S}}_{\mathtt{u}}=OT_{N}$.
Let
$\mathsf{k}_{E}(\iota):=\bigcup\{\{\alpha\}\cup\ E_{N-1}(\alpha):\alpha\in \mathsf{k}(\iota)\}$ for $RS$-terms and $RS$-formulas $\iota$.
Also $\mathsf{k}_{E}(\Gamma)=\bigcup\{\mathsf{k}_{E}(A):A\in\Gamma\}$ for sets $\Gamma$ of formulas.

\bdf\label{df:controlder_diamond}
{\rm
Let $\mathtt{Q}$ be either $\mathtt{Q}=\emptyset$ or a finite family (for $\Lambda,\gamma_{0}$),
$\Theta$ a finite set of ordinals,
and $c\leq\mathbb{S}$. 
Let
$\Gamma=\bigcup\{\Gamma_{\sigma}^{(\sigma)}:\sigma\in\mathtt{Q}\cup\{\mathtt{u}\}\}\subset\Delta_{0}(\mathbb{K})$
a set of formulas
such that 
$\mathsf{k}(\Gamma_{\sigma})\subset E^{\mathbb{S}}_{\sigma}$ 
for each $\sigma\in\mathtt{ Q}\cup\{\mathtt{u}\}$.
$(\mathcal{H}_{\gamma},\Theta, \mathtt{Q})\vdash^{\diamond a}_{c} \Gamma$ holds
if

\begin{equation}
\label{eq:controlder_diamondS}
\{\gamma_{0},\gamma,\Lambda_{0},a,c\}\cup E_{N-1}(\gamma_{0},\gamma,\Lambda_{0},a,c)
 \cup\mathsf{k}_{E}(\Gamma)
\subset
\mathcal{H}_{\gamma}[\Theta]
\end{equation}

and one of the following cases holds:

\begin{description}

\item[$(\bigvee)$]
There exist 
$A\simeq\bigvee(A_{\iota})_{\iota\in J}$, 
an ordinal
$a(\iota)<a$, $A^{(\rho)}\in\Gamma$ with a cap $\rho\in\mathtt{ Q}\cup\{\mathtt{u}\}$,  and an 
$\iota\in [\rho]J$ 
 and
$(\mathcal{H}_{\gamma},\Theta,\mathtt{ Q})\vdash^{\diamond a(\iota)}_{c}\Gamma,
\left(A_{\iota}\right)^{(\rho)}$.

\item[$(\bigwedge)$]
There exist 
$A\simeq\bigwedge(A_{\iota})_{ \iota\in J}$,  a cap $\rho\in\mathtt{Q}\cup\{\mathtt{u}\}$, 
ordinals $a(\iota)<a$ such that 
$A^{(\rho)}\in\Gamma$ and
$(\mathcal{H}_{\gamma},\Theta\cup\mathsf{k}_{E}(\iota),\mathtt{Q}
)
\vdash^{\diamond a(\iota)}_{c}\Gamma,
A_{\iota}^{(\rho)}$
for each $\iota\in [\rho]J$.

\item[$(cut)$]
There exist $\rho\in\mathtt{Q}\cup\{\mathtt{u}\}$, an ordinal $a_{0}<a$ and an uncapped $\bigvee$-formula $C$
such that $\mathrm{rk}(C)<c$,
$(\mathcal{H}_{\gamma},\Theta,\mathtt{ Q})\vdash^{\diamond a_{0}}_{c}\Gamma,\lnot C^{(\rho)}$
and
$(\mathcal{H}_{\gamma},\Theta,\mathtt{ Q})\vdash^{\diamond a_{0}}_{c}C^{(\rho)},\Gamma$.

\item[$(\Sigma(\Omega)\mbox{{\rm -rfl}})$]
There exist ordinals
$a_{\ell}, a_{r}<a$, $\rho\in\mathtt{Q}\cup\{\mathtt{u}\}$
and a formula $C\in\Sigma(\Omega)$ 
such that 
$(\mathcal{H}_{\gamma},\Theta,\mathtt{ Q}
)\vdash^{\diamond a_{\ell}}_{c}\Gamma,C^{(\rho)}$
and
$(\mathcal{H}_{\gamma},\Theta,\mathtt{ Q}
)\vdash^{\diamond a_{r}}_{c}
\left(\lnot \exists x<\Omega\,C^{(x,\Omega)}\right)^{(\rho)}, \Gamma$
with $\Omega< c$.

\end{description}

}
\edf

\blem\label{lem:inversion_diamond}
{\rm (Inversion)}
Let 
$(\mathcal{H}_{\gamma},\Theta,\mathtt{Q})\vdash^{\diamond a}_{c}\Gamma,A^{(\rho)}$,
$A\simeq \bigwedge(A_{\iota})_{\iota\in J}$, 
$\iota\in [\rho]J$.
Then
$(\mathcal{H}_{\gamma},\Theta\cup \mathsf{k}_{E}(\iota),\mathtt{Q}
)
\vdash^{\diamond a}_{c}
\Gamma,\left(A_{\iota}\right)^{(\rho)}$ holds.
\elem
\bprf
This is seen as in Inversion \ref{lem:inversionpi11cap}.
\eprf

\blem\label{lem:reduction_diamond}{\rm (Reduction)}
Let
$(\mathcal{H}_{\gamma},\Theta,\mathtt{Q})\vdash^{\diamond a_{0}}_{c}\Gamma_{0},\lnot C^{(\rho)}$
and
$(\mathcal{H}_{\gamma},\Theta, \mathtt{Q})\vdash^{\diamond a_{1}}_{c} C^{(\rho)},\Gamma_{1}$,
where 
$C\simeq\bigvee(C_{\iota})_{\iota\in J}$ and $\mathrm{rk}(C)\leq c$.

Then
$(\mathcal{H}_{\gamma},\Theta,\mathtt{Q})\vdash^{\diamond a_{0}+a_{1}}_{c}\Gamma_{0},\Gamma_{1}$
holds.
\elem
\bprf
By induction on $a_{1}$ as in Reduction \ref{lem:predcereg}.
Consider the case when the last inference in $(\mathcal{H}_{\gamma},\Theta, \mathtt{ Q})\vdash^{\diamond a_{1}}_{c} C^{(\rho)},\Gamma_{1}$
 is a $(\bigvee)$ with a major formula $C^{(\rho)}$:
 We have
 $(\mathcal{H}_{\gamma},\Theta, \mathtt{ Q})\vdash^{\diamond a_{2}}_{\mathbb{S}}C^{(\rho)},C_{\iota}^{(\rho)},\Gamma_{1}$
 for an $\iota\in [\rho]J$ and an $a_{2}<a_{1}$.
 
  IH yields
$(\mathcal{H}_{\gamma},\Theta, \mathtt{ Q})\vdash^{\diamond a_{0}+a_{2}}_{\mathbb{S}}C_{\iota}^{(\rho)},\Gamma_{0},\Gamma_{1}$.
We obtain
 $(\mathcal{H}_{\gamma},\Theta\cup \mathsf{k}_{E}(\iota),\mathtt{ Q})\vdash^{\diamond a_{0}}_{c}\Gamma_{0},\lnot C_{\iota}^{(\rho)}$
 by Inversion \ref{lem:inversion_diamond}, and
 $(\mathcal{H}_{\gamma},\Theta,\mathtt{ Q})\vdash^{\diamond a_{0}+a_{2}}_{c}\Gamma_{0},\lnot C_{\iota}^{(\rho)}$
 for $a_{0}\leq a_{0}+a_{2}$ and $\mathsf{k}_{E}(\iota)\subset\mathcal{H}_{\gamma}[\Theta]$ by (\ref{eq:controlder_diamondS}).
 We have $\mathrm{rk}(C_{\iota})<\mathrm{rk}(C)\leq c$.
 A $(cut)$ yields the lemma.
\eprf

\blem\label{lem:ce_diamond}{\rm (Cut-elimination)}
Let
$(\mathcal{H}_{\gamma},\Theta,\mathtt{Q})\vdash^{\diamond a}_{c+b}\Gamma$
where $\{c\}\cup E_{N-1}(c)\subset \mathcal{H}_{\gamma}[\Theta]$, $c+b\leq\mathbb{S}$ and
$\lnot(c\leq\Omega<c+b)$.
Then
$(\mathcal{H}_{\gamma},\Theta,\mathtt{Q})\vdash^{\diamond \theta_{b}(a)}_{c}\Gamma$
holds.
\elem
\bprf
By main induction on $b$ with subsidiary induction on $a$ using Reduction \ref{lem:reduction_diamond}.
\eprf

\blem\label{lem:recap_diamond}
Let
$(\mathcal{H}_{\gamma},\Theta,\mathtt{Q}\cup\{\sigma\})\vdash^{\diamond a}_{c}\Gamma,\Pi^{(\sigma)}$, where
$c\leq\mathbb{S}$, $\sigma\in\Psi_{N}$,
$\tau\neq\sigma$ for each $A^{(\tau)}\in\Gamma$, and $\Pi\subset\bigvee(\mathbb{S}+1)$.
Then
$(\mathcal{H}_{\gamma},\Theta,\mathtt{Q})\vdash^{\diamond a}_{c}\Gamma,\Pi^{(\rho)}$ holds for $\sigma\leq\rho\in\mathtt{Q}$
and for $\rho=\mathtt{u}$.
\elem
\bprf
By induction on $a$.
We obtain $\mathsf{k}(\Pi)\subset E^{\mathbb{S}}_{\sigma}\subset E^{\mathbb{S}}_{\rho}$.
Let $A^{(\sigma)}\in\Pi^{(\sigma)}$ and $A\simeq\bigvee(A_{\iota})_{\iota\in J}$.
If $A\in\Delta_{0}(\mathbb{S})$, then
$[\rho]J=[\sigma]J=J=[\mathtt{u}]J$ holds by $\mathsf{k}(A)\subset  E^{\mathbb{S}}_{\sigma}\cap\mathbb{S}=\sigma\leq\rho$.
Let $A\in\Sigma_{1}(\mathbb{S})$. Then $[\sigma]J=Tm(\sigma)\subset Tm(\rho)=[\rho]J$ when $\rho\in\Psi_{N}$, and 
$[\sigma]J=Tm(\sigma)\subset Tm(\mathbb{S})=[\mathtt{u}]J=J$.
Note that each cut formula as well as minor formulas of $(\Sigma(\Omega)\mbox{{\rm -rfl}})$ is a $\Delta_{0}(\mathbb{S})$-formula
since $c\leq\mathbb{S}$.
\eprf

\blem\label{lem:predcereg_S}
Let
$(\mathcal{H}_{\gamma},\Theta, \mathtt{Q})\vdash^{a}_{\mathbb{S}+1,\mathbb{S}+1,N-1,\Lambda_{0},\gamma_{0}}\Gamma$
with $\gamma\leq\gamma_{0}$.

Then 
$(\mathcal{H}_{\gamma_{1}},\Theta(\mathtt{Q}), \mathtt{Q})\vdash^{\diamond \omega^{a}}_{\mathbb{S}}\Gamma$ holds for
$\gamma_{1}=\gamma_{0}+\mathbb{S}$ and
$\Theta(\mathtt{Q})=\Theta\cup E_{\mathbb{S}}(\Theta) \cup  b(\mathtt{Q})$, where
$E_{\mathbb{S}}(\Theta):=\bigcup\{E_{\mathbb{S}}(\alpha):\alpha\in\Theta\}$ and $b(\mathtt{Q})=\{b(\rho):\rho\in\mathtt{Q}\}$.
\elem
\bprf
By induction on $a$.
We have (\ref{eq:controlder_diamondS}) by (\ref{eq:controlder_cap3}).
In the proof let us write
$(\mathcal{H}_{\gamma},\Theta, \mathtt{Q})\vdash^{a}_{\mathbb{S}+1}\Gamma$
for $(\mathcal{H}_{\gamma},\Theta, \mathtt{Q})\vdash^{a}_{\mathbb{S}+1,\mathbb{S}+1,N-1,\Lambda_{0},\gamma_{0}}\Gamma$.
\\
\textbf{Case 1}.
Consider first the case when the last inference is a $(cut)$:
We have
$(\mathcal{H}_{\gamma},\Theta, \mathtt{ Q})\vdash^{a_{0}}_{\mathbb{S}+1}\Gamma,\lnot C^{(\rho)}$
and 
$(\mathcal{H}_{\gamma},\Theta, \mathtt{ Q})\vdash^{a_{0}}_{\mathbb{S}+1}\Gamma, C^{(\rho)}$
for $\rho=\rho_{\mathtt{Q}}$, $a_{0}<a$ and $\mathrm{rk}(C)\leq\mathbb{S}$.
IH yields
$(\mathcal{H}_{\gamma_{1}},\Theta(\mathtt{Q}), \mathtt{Q})\vdash^{\diamond \omega^{a_{0}}}_{\mathbb{S}}\Gamma,\lnot C^{(\rho)}$
and
$(\mathcal{H}_{\gamma_{1}},\Theta(\mathtt{Q}), \mathtt{Q})\vdash^{\diamond \omega^{a_{0}}}_{\mathbb{S}}\Gamma, C^{(\rho)}$.
By Reduction \ref{lem:reduction_diamond} we obtain
$(\mathcal{H}_{\gamma_{1}},\Theta(\mathtt{Q}), \mathtt{ Q})\vdash^{\diamond \omega^{a}}_{\mathbb{S}}\Gamma$
for $\omega^{a_{0}}\cdot 2<\omega^{a}$.
\\
\textbf{Case 2}.
Second consider the case when the last inference is a degenerated $({\rm rfl}(\rho,e,f))$:
We have $a_{0}<a$,  $\rho=\rho_{\mathtt{Q}}$, $SC_{\Lambda}(f)\subset \mathcal{H}_{\gamma}[\Theta\cap E^{\mathbb{S}}_{\rho}]$,
$(\mathcal{H}_{\gamma},\Theta, \mathtt{Q})\vdash^{a_{0}}_{\mathbb{S}+1}\Gamma,\lnot \delta^{(\rho)}$
for each $\delta\in\Delta$.
IH yields
\begin{equation}\label{eq:predcereg_S1}
(\mathcal{H}_{\gamma_{1}},\Theta(\mathtt{Q}), \mathtt{Q})\vdash^{\diamond \omega^{a_{0}}}_{\mathbb{S}}\Gamma,\lnot \delta^{(\rho)}
\end{equation}

On the other hand we have
$(\mathcal{H}_{\gamma},\Theta\cup\{\sigma\}, \mathtt{Q}\cup\{\sigma\})\vdash^{a_{0}}_{\mathbb{S}+1}\Xi ,\Delta^{(\sigma)}$
for $\Xi\subset\Gamma$ and $\Delta\subset\bigvee(\mathbb{S}+1)$, where
$\sigma$ ranges over ordinals such that
$\sigma\in H_{\rho}(f,\gamma_{0},\Theta)$ with $m(\sigma)=f$, and
$\tau\neq\sigma$ for each $A^{(\tau)}\in\Xi\subset\Gamma$.
We have $(\Theta\cup\{\sigma\})(\mathtt{Q}\cup\{\sigma\})=\Theta(\mathtt{Q})\cup\{b(\sigma),\sigma\}$.

We have $\gamma_{0}\leq\nu=b(\rho)<\gamma_{0}+\mathbb{S}$ by Definition \ref{df:cap}.\ref{df:cap.2}
and $SC_{\Lambda}(f)\subset \mathcal{H}_{\gamma}[\Theta\cap E^{\mathbb{S}}_{\rho}]\subset E^{\mathbb{S}}_{\rho}$.
Furthermore
$\{\rho,\nu\}\cup E_{\mathbb{S}}(\Theta)\subset\Theta(\mathtt{Q})$ by (\ref{eq:controlder_cap1}).
Let $\sigma=\psi_{\rho}^{f}(\nu+\alpha)$ for $\alpha=\max(\{1\}\cup (E_{\mathbb{S}}(\Theta)\cap E^{\mathbb{S}}_{\rho}))$.
We see $\nu+\alpha<\gamma_{0}+\mathbb{S}=\gamma_{1}$ from $\alpha<\mathbb{S}$.
Hence $\{b(\sigma),\sigma\}\subset\mathcal{H}_{\gamma_{1}}[\Theta(\mathtt{Q})]$, 
$\sigma\in H_{\rho}(f,\gamma_{0},\Theta)$ and $m(\sigma)=f$.

$(\mathcal{H}_{\gamma_{1}},\Theta(\mathtt{Q}), \mathtt{Q}\cup\{\sigma\})\vdash^{\diamond \omega^{a_{0}}}_{\mathbb{S}}\Xi, \Delta^{(\sigma)}$
follows from IH, and by Lemma \ref{lem:recap_diamond}
\begin{equation}\label{eq:predcereg_S2}
(\mathcal{H}_{\gamma_{1}},\Theta(\mathtt{Q}), \mathtt{Q})\vdash^{\diamond \omega^{a_{0}}}_{\mathbb{S}}\Xi ,\Delta^{(\rho)}
\end{equation}
Let $n=\#(\Delta)$ be the number of formulas in $\Delta$.
$(\mathcal{H}_{\gamma_{1}},\Theta(\mathtt{Q}), \mathtt{Q})\vdash^{\diamond \omega^{a_{0}}(n+1)}_{\mathbb{S}}\Gamma$
follows from (\ref{eq:predcereg_S1}), (\ref{eq:predcereg_S2}) and Reduction \ref{lem:reduction_diamond}, where
$\omega^{a_{0}}\cdot (n+1)<\omega^{a}$ by $a_{0}<a$.

Other cases are seen from IH.
In $(\bigvee)$ with a minor formula $(A_{\iota})^{(\rho)}$ and $A\simeq\bigvee(A_{\iota})_{\iota\in J}$,
we obtain $\iota\in[\rho]J$ by (\ref{eq:controlder_cap2}).
\eprf

\bdf\label{df:controlder_diamond_empty}
{\rm
For sets $\Gamma$ of uncapped formulas,
let
\[
(\mathcal{H}_{\gamma},\Theta)\vdash^{\diamond a}_{c}\Gamma:\Leftrightarrow
(\mathcal{H}_{\gamma},\Theta,\emptyset)\vdash^{\diamond a}_{c}\Gamma
\]
for the empty family $\mathtt{Q}=\emptyset$.
}
\edf

\blem\label{lem:main}{\rm (Uncapping)}
Suppose
$ 
(\mathcal{H}_{\gamma_{1}},\Theta,\mathtt{Q})
\vdash^{\diamond a}_{\mathbb{S}}
\Gamma^{(\cdot)}
$ for $\gamma_{1}=\gamma_{0}+\mathbb{S}$,
where 
$\Gamma\subset\Delta_{0}(\mathbb{S})$
is a set of uncapped formulas,
$\Gamma=\bigcup\{\Gamma_{\rho}:\rho\in\mathtt{Q}\cup\{\mathtt{u}\}\}$,
and
$\Gamma^{(\cdot)}=\bigcup\{\Gamma_{\rho}^{(\rho)}:\rho\in\mathtt{ Q}\cup\{\mathtt{u}\}\}$.
Then 
$(\mathcal{H}_{\gamma_{1}},\Theta)
\vdash^{\diamond a}_{\mathbb{S}}
\Gamma$ holds.
\elem
\bprf
This is seen from Lemma \ref{lem:recap_diamond}.
\eprf

\blem\label{lem:collpase_diamond}{\rm (Collapsing)}
Assume
$\Theta\subset
\mathcal{H}_{\gamma}(\psi_{\Omega}(\gamma))$
 and
$
(\mathcal{H}_{\gamma},\Theta
)\vdash^{\diamond a}_{\Omega+1}\Gamma
$ with $\Gamma\subset\Sigma(\Omega)$.
Then 
$(\mathcal{H}_{\hat{a}+1},\Theta)
\vdash^{\diamond \beta}_{\beta}
\Gamma^{(\beta,\Omega)}$ holds for $\hat{a}=\gamma+\omega^{a}$ and $\beta=\psi_{\Omega}(\hat{a})$.
\elem
\bprf
This is seen as in \cite{Buchholz}
by induction on $a$.
\eprf

\subsection{Proof of Theorem \ref{thm:2}}\label{subsec:final}

Let us prove Theorem \ref{thm:2}.
Let 
$T_{N}\vdash\theta^{\mathsf{L}_{\Omega}}$ for a $\Sigma$-sentence $\theta$.
By Lemma \ref{th:embedregthm} pick an $m$ so that 
 $(\mathcal{H}_{0},\emptyset)
 \vdash_{\mathbb{K}+1+m}^{* \mathbb{K}\cdot 2+m}
 \theta^{\mathsf{L}_{\Omega}} $.
Let $\gamma_{0}:=\omega_{m+2}(\mathbb{K}+1)\in\mathcal{H}_{0}[\emptyset]$.

Cut-elimination \ref{lem:predcereg*} yields
$(\mathcal{H}_{0},\emptyset)\vdash^{* a}_{\mathbb{K}+1} \theta^{\mathsf{L}_{\Omega}}$, where 
$a=\omega_{m}(\mathbb{K}\cdot 2+m)$.
Let $\hat{a}=\omega^{a}<\gamma_{0}$ and $\Lambda_{0}=\beta_{0}=\psi_{\mathbb{K}}(\hat{a})$.
We obtain 
$(\mathcal{H}_{\hat{a}+1},\emptyset)\vdash^{* \Lambda_{0}}_{\Lambda_{0}} \theta^{\mathsf{L}_{\Omega}}$
by Collapsing \ref{lem:Kcollpase*}.

In what follows each finite function is an $f:\Lambda\to\Gamma(\Lambda)$ with $\Lambda=\Gamma(\Lambda_{0})$.
Let
$\rho_{0}=\psi_{\mathbb{S}}^{g_{0}}(\gamma_{0})\in\mathcal{H}_{\gamma_{0}+\mathbb{S}}[\emptyset]$ with
$\mathrm{supp}(g_{0})=\{\Lambda_{0}\}$ with $s(\rho_{0})=\Lambda_{0}>\mathbb{S}+1$ and $g_{0}(\Lambda_{0})=\Lambda\cdot 3$.
Capping \ref{lem:capping} yields
$(\mathcal{H}_{\hat{a}+1},\{\rho_{0}\},\mathtt{Q}_{0})
\vdash_{\Lambda_{0},\Lambda_{0},0,\Lambda_{0},\gamma_{0}}^{\Lambda_{0}\cdot 2} (\theta^{\mathsf{L}_{\Omega}})^{(\rho_{0})}$
for $\mathtt{Q}_{0}=\{\rho_{0}\}$.
In the following we write $\vdash^{a}_{c,d,m}$ for $\vdash^{a}_{c,d,m,\Lambda_{0},\gamma_{0}}$, and
let $\pi_{m}=\Omega_{\mathbb{S}+N-m}$ for $m\leq N$ and $\mathbb{K}=\pi_{0}$.

We have
$(\mathcal{H}_{b_{0}+1},\{\rho_{0}\},\mathtt{Q}_{0})
\vdash_{\beta_{0},\beta_{0},0}^{a_{0}} (\theta^{\mathsf{L}_{\Omega}})^{(\rho_{0})}$ for $b_{0}=\hat{a}$,
$\beta_{0}=\Lambda_{0}$ and $a_{0}=\Lambda_{0}\cdot 2$.
We obtain
$(\mathcal{H}_{b_{0}+1},\{\rho_{0}\},\mathtt{Q}_{0})
\vdash_{\pi_{1}+1,\beta_{0},0}^{c} (\theta^{\mathsf{L}_{\Omega}})^{(\rho_{0})}$
for $c=\varphi_{\beta_{0}}(a_{0})<\Lambda$
by Cut-elimination \ref{lem:CE}.
We have $\pi_{1}+1<\Lambda_{0}=s(\rho_{0})$.
Let $g_{1}=h^{\pi_{1}+1}(g_{0};\varphi_{\beta_{0}}(c))$ and $\rho_{1}=\psi_{\rho_{0}}^{g_{1}}(\gamma_{0}+1)\in H_{\rho_{0}}(g_{1},\gamma_{0},\emptyset)$.
Also let $\mathtt{Q}_{1}=\mathtt{Q}^{[\rho_{1}/\rho_{0}]}=\{\rho_{1}\}$.
Recapping \ref{mlem:singlemainl.1} then yields  for $s(\rho_{1})=\pi_{1}+1$,
$(\mathcal{H}_{b_{0}+1},\{\rho_{1}\},\mathtt{Q}_{1})
\vdash_{\pi_{1}+1,\beta_{0},0}^{\varphi_{\beta_{0}}(c)} (\theta^{\mathsf{L}_{\Omega}})^{(\rho_{1})}$.
We obtain $\rho_{\mathtt{Q}_{1}}=\rho_{1}\in\mathcal{H}_{\gamma_{1}}[\emptyset]$ and
$\varphi_{\beta_{0}}(c)<\Lambda$.
Let $\eta_{1}=on_{0}(\mathtt{Q}_{1})=g_{1}(\pi_{1}+1)$.
Lemma \ref{mlem:singlemainl_S} yields
$(\mathcal{H}_{b_{0}+1},\{\rho_{1}\},\mathtt{Q}_{1})\vdash^{c_{1}}_{\pi_{1}+1,\pi_{1}+1,0}  (\theta^{\mathsf{L}_{\Omega}})^{(\rho_{1})}$
for $c_{1}=\varphi_{\beta_{0}+\eta_{1}}(\varphi_{\beta_{0}}(c))$.

Collapsing \ref{lem:picollpase} yields
$(\mathcal{H}_{b_{1}+1},\{\rho_{1}\},\mathtt{Q}_{1})\vdash^{a_{1}}_{\beta_{1},\beta_{1},1}  (\theta^{\mathsf{L}_{\Omega}})^{(\rho_{1})}$
for
$a_{1}=\beta_{1}=\psi_{\pi_{1}}(b_{1})$ and $b_{1}=b_{0}+\omega^{c_{1}}$.

Let $\beta_{0}=\Lambda_{0}$, $a_{0}=\Lambda_{0}\cdot 2$, $b_{0}=\omega^{a}$, and
$\mathtt{Q}_{0}=\{\rho_{0}\}$ with
$\rho_{0}=\psi_{\mathbb{S}}^{g_{0}}(\gamma_{0})$.
For $m< N$, let $g_{m+1}=h^{\pi_{m+1}+1}(g_{m};\varphi_{\beta_{m}}(\varphi_{\beta_{m}}(a_{m})))$,
$\rho_{m+1}=\psi_{\rho_{m}}^{g_{m+1}}(\gamma_{0}+m+1)$,
$\mathtt{Q}_{m+1}=\{\rho_{m+1}\}$, $\eta_{m+1}=on_{m}(\mathtt{Q}_{m+1})$, and
$c_{m+1}=\varphi_{\beta_{m}+\eta_{m+1}}(\varphi_{\beta_{m}}(\varphi_{\beta_{m}}(a_{m})))$.
Also let
$a_{m+1}=\beta_{m+1}=\psi_{\pi_{m+1}}(b_{m+1})$, and $b_{m+1}=b_{m}+\omega^{c_{m+1}}$ for $m<N-1$.

We see inductively that $a_{m},\beta_{m}<\Lambda$, $b_{m},c_{m}<\gamma_{0}$,
$\{a_{m},b_{m},c_{m},\beta_{m},\eta_{m},\rho_{m}\}\cup SC_{\Lambda}(g_{m})\subset\mathcal{H}_{\gamma_{0}+\mathbb{S}}[\emptyset]$,
 and for each $m<N-1$,
$(\mathcal{H}_{b_{m}+1},\{\rho_{m}\},\mathtt{Q}_{m})\vdash^{a_{m}}_{\beta_{m},\beta_{m},m}  (\theta^{\mathsf{L}_{\Omega}})^{(\rho_{m})}$.
We obtain 
$(\mathcal{H}_{b_{N-1}+1},\{\rho_{N}\},\mathtt{Q}_{N})\vdash^{c_{N}}_{\pi_{N}+1,\pi_{N}+1,N-1}  (\theta^{\mathsf{L}_{\Omega}})^{(\rho_{N})}$
for $\pi_{N}=\mathbb{S}$.

Let $\gamma_{1}=\gamma_{0}+\mathbb{S}$.
By Lemma \ref{lem:predcereg_S} we obtain
$(\mathcal{H}_{\gamma_{1}},\emptyset, \mathtt{Q}_{N} )\vdash^{\diamond c_{N}}_{\mathbb{S}}  (\theta^{\mathsf{L}_{\Omega}})^{(\rho_{N})}$, where
$\omega^{c_{N}}=c_{N}$ and  $\Theta(\mathtt{Q}_{N})=\{\rho_{N},\gamma_{0}+N\}\subset\mathcal{H}_{\gamma_{1}}[\emptyset]$ 
for $\Theta=\{\rho_{N}\}$.
Uncapping \ref{lem:main} yields 
$(\mathcal{H}_{\gamma_{1}},\emptyset,\emptyset)\vdash^{\diamond c_{N}}_{\mathbb{S}}  \theta^{\mathsf{L}_{\Omega}}$.
In what follows we write $(\mathcal{H}_{\gamma_{1}},\emptyset)$ for $(\mathcal{H}_{\gamma_{1}},\emptyset,\emptyset)$.

By Cut-elimination \ref{lem:ce_diamond} we obtain
$(\mathcal{H}_{\gamma_{1}},\emptyset)\vdash^{\diamond d}_{\Omega+1}  \theta^{\mathsf{L}_{\Omega}}$
for $d=\theta_{\mathbb{S}}(c_{N})<\mathbb{K}$.
Collapsing \ref{lem:collpase_diamond} then yields
$(\mathcal{H}_{\gamma_{1}+d+1},\emptyset)\vdash^{\diamond \delta}_{\delta}  \theta^{\mathsf{L}_{\delta}}$
for
$\delta=\psi_{\Omega}(\gamma_{1}+d)<\psi_{\Omega}(\omega_{m+2}(\mathbb{K}+1))$ with $\omega^{d}=d$.
We finally obtain
$(\mathcal{H}_{\gamma_{1}+d+1},\emptyset)\vdash^{\diamond \theta_{\delta}(\delta)}_{0}  \theta^{\mathsf{L}_{\delta}}$
by Cut-elimination \ref{lem:ce_diamond}.
We conclude
$L_{\delta}\models\theta$
by induction up to $\theta_{\delta}(\delta)$.

\bcor\label{cor:Pi2conservative}
${\sf KP}\ell^{r}+(M\prec_{\Sigma_{1}}V)$ is 
conservative over ${\rm I}\Sigma_{1}+\{TI(\alpha,\Sigma^{0}_{1}(\omega)): \alpha<\psi_{\Omega}(\Omega_{\mathbb{S}+\omega})\}$
with respect to $\Pi^{0}_{2}(\omega)$-arithmetic sentences, and 
each provably computable function in ${\sf KP}\ell^{r}+(M\prec_{\Sigma_{1}}V)$ is defined
by $\alpha$-recursion for an $\alpha<\psi_{\Omega}(\Omega_{\mathbb{S}+\omega})$.
\ecor
\bprf
Let $\theta$ be a $\Pi^{0}_{2}(\omega)$-arithmetic sentence on $\omega$, and assume
that ${\sf KP}\ell^{r}+(M\prec_{\Sigma_{1}}V)$ proves $\theta$.
Pick an $N>0$ such that $T_{N}\vdash\theta$.
Theorem \ref{thm:2} shows that $\theta$ is true.
The proof of Theorem \ref{thm:2} is seen to be formalizable in an intuitionistic fixed point theory
${\sf FiX}^{i}(T)$
over an extension 
$T={\sf PA}+\{TI(\alpha): \alpha<\psi_{\Omega}(\varepsilon_{\Omega_{\mathbb{S}+N}+1})\}$
of the first order arithmetic {\sf PA},
where transfinite induction schema 
applied to arithmetical formulas with fixed point predicates is available 
up to each $\alpha<\psi_{\Omega}(\varepsilon_{\Omega_{\mathbb{S}+N}+1})$ in the extension $T$.
From \cite{intfix} 
we see that ${\sf FiX}^{i}(T)$ is a conservative extension of $T$.
Therefore $T\vdash\theta$.
Since the ordinal $\psi_{\Omega}(\Omega_{\mathbb{S}+\omega})$ is an epsilon number,
we see that $\theta$ is provable in 
${\rm I}\Sigma_{1}+\{TI(\alpha,\Sigma^{0}_{1}(\omega)): \alpha<\psi_{\Omega}(\Omega_{\mathbb{S}+\omega})\}$.

Conversely we see that
${\sf KP}\ell^{r}+(M\prec_{\Sigma_{1}}V)$ proves
$TI(\alpha,\Sigma^{0}_{1}(\omega))$ up to each $\alpha<\psi_{\Omega}(\Omega_{\mathbb{S}+\omega})$
from Theorem \ref{th:wf} in \cite{singledistwfprf}.
\eprf

\end{document}

%% file: header.tex
\newcommand{\brem}{\begin{remark}}
\newcommand{\erem}{\end{remark}}
\newcommand{\blem}{\begin{lemma}}
\newcommand{\elem}{\end{lemma}}
\newcommand{\bth}{\begin{theorem}}
\newcommand{\ethm}{\end{theorem}}
\newcommand{\benu}{\begin{enumerate}}
\newcommand{\eenu}{\end{enumerate}}
\newcommand{\bdes}{\begin{description}}
\newcommand{\edes}{\end{description}}
\newcommand{\bdf}{\begin{definition}}
\newcommand{\edf}{\end{definition}}
\newcommand{\bcor}{\begin{cor}}
\newcommand{\ecor}{\end{cor}}
\newcommand{\bprp}{\begin{proposition}}
\newcommand{\eprp}{\end{proposition}}
\newcommand{\bmlem}{\begin{mlemma}}
\newcommand{\emlem}{\end{mlemma}}
\newcommand{\bclm}{\begin{claim}}
\newcommand{\eclm}{\end{claim}}
\newcommand{\bprf}{{\bf Proof}.\hspace{2mm}}

\newcommand{\eprf}{\hspace*{\fill} $\Box$}

\newcommand{\beqn}{\begin{equation}}
\newcommand{\eeqn}{\end{equation}}
\newcommand{\beqnarr}{\begin{eqnarray}}
\newcommand{\eeqnarr}{\end{eqnarray}}
\newcommand{\beqnarrs}{\begin{eqnarray*}}
\newcommand{\eeqnarrs}{\end{eqnarray*}}

\newtheorem{theorem}{Theorem}[section]
\newtheorem{definition}[theorem]{Definition}
\newtheorem{proposition}[theorem]{Proposition}
\newtheorem{lemma}[theorem]{Lemma}
\newtheorem{cor}[theorem]{Corollary}
\newtheorem{remark}[theorem]{Remark}
\newtheorem{mlemma}[theorem]{Main Lemma}
\newtheorem{claim}[theorem]{Claim}